\documentclass[12pt, a4paper]{amsart}

\usepackage[english]{babel}
\usepackage{amsfonts}
\usepackage{amsmath}
\usepackage{amssymb, latexsym}
\usepackage{graphicx}
\usepackage[usenames]{color}

\usepackage{fullpage}
\usepackage{amsmath}

\usepackage{mathrsfs}

\newtheorem{theorem}{Theorem}[section]
\newtheorem{proposition}[theorem]{Proposition}
\newtheorem{corollary}[theorem]{Corollary}
\newtheorem{lemma}[theorem]{Lemma}
\newtheorem{remark}[theorem]{Remark}
\newtheorem{definition}[theorem]{Definition}

\numberwithin{equation}{section} 

\numberwithin{equation}{section}

\newcommand{\R}{\mathbb{R}}

\newcommand{\ben}{\begin{eqnarray*}}
\newcommand{\enn}{\end{eqnarray*}}
\newcommand{\pa}{\partial}

\newcommand{\g}{\gamma}

\newcommand{\ve}{\varepsilon}

\newcommand{\al}{\alpha}
\newcommand{\la}{\lambda}

\newcommand{\ol}{\overline}

\newcommand{\half}{\frac{1}{2}}

\newcommand{\na}{\nabla}
\newcommand{\be}{\begin{equation}}
\newcommand{\ee}{\end{equation}}
\newcommand{\ba}{\begin{aligned}}
\newcommand{\ea}{\end{aligned}}
\newcommand{\lf}{\left}
\newcommand{\rt}{\right}

\newcommand{\wt}{\tilde{w}}

\def\9{{\infty}}
\def\a{{\alpha}}
\def\b{{\beta}}
\def\g{{\gamma}}

\def\lbb{{\lambda}}
\def\t{{\theta}}

\def\calo{{\mathcal{O}}}
\def\calp{{\mathcal{P}}}

\def\bbp{{\mathbb{P}}}
\def\bbr{{\mathbb{R}}}
\def\ve{{\varepsilon}}
\def\vf{{\varphi}}

\def\wt{\widetilde}
\def\wh{\widehat}
\def\ol{\overline}
\def\({\left(}
\def\){\right)}
\def\<{\langle}
\def\>{\rangle}

\begin{document}

\title{On the multi-bubble blow-up solutions to rough nonlinear Schr\"odinger equations}

\author{Yiming Su}
\address[Y.M. Su]{Department of mathematics,
Zhejiang University of Technology, 310014 Zhejiang, China.}
\author{Deng Zhang}
\address[D. Zhang]{School of mathematical sciences,
Key laboratory of sientific and engineering computing (Ministry of Education),
Shanghai Jiao Tong University, 200240 Shanghai, China.}
\begin{abstract}
  We are concerned with the multi-bubble blow-up solutions to rough nonlinear Schr\"odinger equations in
  the focusing mass-critical case.
  In both dimensions one and two,
  we construct the finite time multi-bubble solutions,
  which concentrate at $K$ distinct points, $1\leq K<\9$,
  and behave asymptotically like a sum of pseudo-conformal blow-up solutions
  in the pseudo-conformal space $\Sigma$ near the blow-up time.
  The upper bound of the asymptotic behavior
  is closely related to the flatness of noise at blow-up points.
  Moreover,
  we prove the conditional uniqueness of multi-bubble solutions
  in the case where the asymptotic behavior in the energy space $H^1$
  is of the order $(T-t)^{3+\zeta}$,
  $\zeta>0$.
  These results are also obtained for  nonlinear Schr\"odinger equations
  with lower order perturbations, particularly,
  in the absence of the classical pseudo-conformal symmetry
  and the conversation law of energy.
  The existence results are applicable to the canonical deterministic nonlinear Schr\"odinger equation
  and complement the previous work \cite{M90}.
  The conditional uniqueness results are new in both
  the stochastic and deterministic case.
\end{abstract}

\keywords{Blow up, mass-critical, multi-bubbles, nonlinear Schr\"odinger equation, rough path}
\subjclass[2010]{60H15, 35B44, 35B40, 35Q55}

\maketitle
\tableofcontents
\section{Introduction}    \label{Sec-Intro}

This work is devoted to the existence and uniqueness of blow-up solutions
at multiple points
to the rough nonlinear Schr\"{o}dinger equations
in the focusing mass-critical case.
Precisely,
we consider
\begin{align} \label{equa-X-rough}
    & idX = - \Delta X dt - |X|^{\frac 4d} X dt - i \mu X dt +i X dW(t), \\
    & X(0) = X_0 \in H^1(\bbr^d).   \nonumber
\end{align}
Here,
$W$ is the Wiener process of the form
$$W(t,x)=\sum_{k=1}^N i\phi_k(x)B_k(t),\ \ x\in \bbr^d,\ t\geq 0,$$
where
$\phi_k \in C_b^\9(\bbr^d, \bbr)$,
$B_k$ are the standard $N$-dimensional real valued Brownian motions
on a stochastic basis $(\Omega, \mathscr{F}, \{\mathscr{F}_t\}, \bbp)$,
$1\leq k\leq N$,
and
$\mu= \frac 12 \sum_{k=1}^N  \phi_k ^2.$
The last term $XdW(t)$ in \eqref{equa-X-rough} is taken in the sense of controlled rough path
(see Definition \ref{def-X-rough} below).
In particular,
the rough integration coincides with the usual It\^o integration
if the corresponding processes are $\{\mathscr{F}_t\}$-adapted
(see \cite[Chapter 5]{FH14}).
For simplicity we assume that $N<\9$,
but the arguments below can be also applied to the case when $N=\9$
under suitable summability conditions.

The noise here is mainly considered of {\it conservative type},
i.e., ${\rm Re} W=0$.
In this case,
the quantum system evolves on the unit ball
if the initial state is normalized $\|X_0\|_{L^2}=1$
and thus verifies the conservation of probability.

One significant model
of nonlinear Schr\"odinger equations with noise
arises from the molecular aggregates
with thermal fluctuations,
where the multiplicative noise corresponds to
scattering of exciton by phonons,
due to thermal vibrations of the molecules.
The noise effect on the coherence of the ground state solitary solution
was studied in \cite{BCIR94,BCIRG95} for
the two dimensional case with the critical cubic nonlinearity.
The influence of noise on collapse
in the one dimensional case with quintic nonlinearity
was  studied in  \cite{RGBC95}.

Another important application
is related to the open quantum systems \cite{BG09},
in which the stochastic perturbation $iX \phi_k dB_k$
represents a stochastic continuous measurement via
the pointwise quantum observable $R_k(X) = X \phi_k$,
while $B_k$ represents the output of continuous measurement, $1\leq k \leq N$.
We refer to \cite[Chapter 2]{BG09} for more physical interpretations.
See also \cite{SS99}
for other physical applications
of Schr\"odinger equations.

The local well-posedness of equation \eqref{equa-X-rough} is quite well known
when the stochastic integration is taken in the sense of It\^o or rough path.
See, e.g., \cite{BRZ16,BM13,BD03,SZ20}.

However,
the situation becomes much more delicate
for the large time behavior of solutions.
As a matter of fact,
solutions may formalize singularities in finite time in the focusing mass-(super)critical case.

When the input noise is of conservative type,
it was first proved by de Bouard and Debussche  \cite{BD02,BD05}
that the noise has the effect to accelerate blow-up with positive probability
in the  mass-supercritical case
(i.e., the exponent of nonlinearity is in the region $(1+\frac 4d, 1+ \frac{4}{(d-2)_+})$).
When the noise is of non-conservative type,
the explosion, however, can be prevented with high probability
as long as the strengthen of noise is large enough,
which  reflects the damped effect of non-conservative noise (\cite{BRZ17}).
We also refer to \cite{MR20} for the
global well-posedness below the threshold
in the mass-(super)critical case.

Moreover,
many numerical simulations have been  made to investigate the blow-up phenomena in the stochastic case.
It was observed in \cite{BDM01,DM02,DM02.2} that
the colored multiplicative noise has the effect to delay blow-up,
while the white noise may even prevent blow-up.
Such phenomena have been also confirmed by the recent numerical results in \cite{MRRY20,MRRY20.0}.
The noise effects on the energy, global well-posedness
and blow-up profiles are also studied in
\cite{MRRY20,MRRY20.0},
which partially confirm
the conjecture that,
in the mass-critical case
the stable blow-up solutions with slightly supercritical mass
shall have the $\log\log$ blow-up rate,
while in the mass-supercritical case the blow-up rate
is of a different polynomial type.

Recently,
the minimal mass blow-up solutions have been constrcuted
by the authors in \cite{SZ20},
and it is shown that the mass of the ground state
characterizes the threshold of global well-posedness
and blow-up in the stochastic case.
The log-log blow-up solutions have been also constructed in \cite{FSZ20}
in the stochastic case.
We also would like to refer to the recent work \cite{FM20}
for the log-log blow-up solutions with $L^2$-regularity
randomized initial data.

In this paper,
we are mainly interested in the blow-up dynamics
in the  large mass regime,
particulary, the existence and uniqueness
of multi-bubble blow-up solutions in the stochastic case.
It should be mentioned that,
the presence of noise destroys
the symmetries and several conservation laws,
which makes it rather difficult to obtain
the  quantitative descriptions
of blow-up dynamics.

Before stating the main results,
let us first review the existing results
for the deterministic nonlinear Schr\"odinger equation (NLS)
\be \label{equa-NLS}
\left\{ \begin{aligned}
\ &i\partial_tu+\Delta u+|u|^\frac{4}{d}u=0,\\
\ &u(0)=u_{0} \in H^1(\bbr^d),
\end{aligned}\right.
\ee
Equation \eqref{equa-NLS} admits a number of symmetries and conservation laws.
It is invariant under the translation, scaling, phase rotation and Galilean transform,
i.e.,
if $u$ solves (\ref{equa-NLS}),
then so does
\begin{align}\label{symmerty}
v(t,x)=\lambda_0^{-\frac{d}{2}}
u (\frac{t-t_0}{\lambda_0^2},\frac{x-x_0 }{\lambda_0}- \frac{\beta_0(t-t_0)}{\lbb_0})
e^{i\frac{\beta_0}{2}\cdot(x-x_0)-i\frac{|\beta_0|^2}{4}(t-t_0)+i\t_0},
\end{align}
with $v(t_0,x)=\lambda_0^{-\frac{d}{2}} u_0 (\frac{x-x_0}{\lambda_0}) e^{i\frac{\beta_0}{2}\cdot(x-x_0)+i\t_0}$,
where
$(\lambda_0, \beta_0, \theta_0) \in \bbr^+ \times \bbr^d \times \bbr$,
$x_0\in\R^d$, $t_0\in\R$.
In particular,
the $L^2$-norm of solutions is  preserved under the symmetries above,
and thus
\eqref{equa-NLS} is  called the {\it mass-critical} equation.
Another symmetry,
particularly important in the blow-up analysis,
is related to the
{\it pseudo-conformal transformation} in the pseudo-conformal space
$\Sigma:=\{u\in H^1(\bbr^d), \|x u\|_{L^2(\bbr^d)}<\9 \}$,
\begin{align}\label{pseudo}
 (-t)^{-\frac{d}{2}} u (\frac{1}{-t},\frac{x}{-t}) e^{-i\frac{|x|^{2}}{4t}}, \quad\ \ t\not =0.
\end{align}
The conservation laws related to \eqref{equa-NLS} contain
\begin{align}
   & {\rm Mass}:\ \ M(u)(t):=\int_{\R^d}|u(t)|^2dx=M(u_0).  \label{mass}  \\
   & {\rm Energy}:\ \ E(u)(t) := \half\int_{\R^d}|\nabla u(t)|^2dx-\frac{d}{2d+4}\int_{\R^d}|u(t)|^{2+\frac{4}{d}}dx=E(u_0). \label{energy} \\
   & {\rm Momentum}: \ \ Mom(u):={\rm Im} \int_{\R^d}\nabla u\bar{u}dx
       =  Mom(u_0). \label{momentum}
\end{align}

An important role here is played by the ground state $Q$,
which is the unique positive spherically symmetric solution to the elliptic equation
\begin{align} \label{equa-Q}
    \Delta Q - Q + Q^{1+\frac{4}{d}} =0.
\end{align}
It is well known that (see, e.g., \cite{Wenn}),
 the mass of the ground state characterizes the threshold for
the global well-posedness and blow-up of solutions to NLS.
More precisely,
solutions to \eqref{equa-NLS} exist globally
if the  initial data have the {\it subcritical mass}, i.e.,  $\|u_0\|_{L^2} < \|Q\|_{L^2}$,
while solutions may formolize singularities in finite time
in the {\it critical mass} case, i.e., $\|u_0\|_{L^2} = \|Q\|_{L^2}$.
In particular,
by virtue of the pseudo-conformal transformation \eqref{pseudo},
one may construct the so-called {\it pseudo-conformal blow-up solutions}
\begin{align}  \label{S}
S_{T}(t,x)=(\omega(T-t))^{-\frac d2}Q (\frac{x-\alpha}{\omega(T-t)}) e^{-\frac i4\frac{|x-\alpha|^{2}}{T-t}+\frac{i}{\omega^2(T-t)}+i\theta},\ \ T\in \bbr.
\end{align}
Note that,
$\|S_T\|_{L^2}=\|Q\|_{L^2}$,
and $S_T$ blows up at time $T$ with the blow-up rate $\sim (T-t)^{-1}$.
Thus, $S_T$ is the {\it minimal mass blow-up solution}.

When the mass of initial data is slightly above $\|Q\|_{L^2}$,
two different blow-up scenarios have been observed:
the pseudo-conformal blow-up rate $\sim (T-t)^{-1}$,
and the {\it log-log} blow-up rate $\sim\sqrt{\frac{log|log(T-t)|}{T-t}}$.
For the blow-up and classification results in this case
we refer to \cite{B-W,MR05.0,MR05,P}
and the references therein.

For even larger mass of initial data,
the complete characterization of the formation of singularity is still an open problem.
It is conjectured by Merle and Rapha\"{e}l \cite{MR05.0} that every $H^1$ blow-up solution
can be decomposed into a singular part and a $L^2$ residual,
and the singular part expands asymptotically as multiple bubbles concentrating at a finite number of points.
This conjecture is known as the blow-up version of the {\it soliton resolution conjecture}.

Thus,
an important step to understand the singularity formulation
in the large mass regime is
to construct multi-bubble blow-up solutions.

In the pioneering work \cite{M90},
Merle initiated the construction
of blow-up solutions concentrating at arbitrary $K(<\9)$ distinct points,
which behave asymptotically like a sum of $K$ pseudo-conformal blow-up solutions
and thus have the pseudo-conformal blow-up rate.
The multi-bubble solutions with the log-log blow-up rate
have been constructed by Fan \cite{F17}.
Moreover,
Martel and Raph\"{e}al \cite{MP18} constructed blow-up solutions with  multiple bubbles concentrating at exactly the same point.
Bubbling phenomena have been also exhibited in various other settings.
See, e.g.,
\cite{J17} for the energy-critical NLS,
and \cite{SG19} for the nonlinear Schr\"{o}dinger system.
We also would like to refer to
\cite{CM18,MMR15} for
the generalized Korteweg-de Vries equations (gKDV),
 \cite{J19,JL18,KST08} for the wave maps,
and \cite{CPM20,M92} for the nonlinear heat equations.

However,
one major challenge in the stochastic case is,
that the symmetries and several conservation laws
are  destroyed, because of the presence of noise.
Equation \eqref{equa-X-rough} is no longer invariant under the pseudo-conformal symmetry,
which, however, is the key ingredient in the classification of minimal mass blow-up solutions to NLS in \cite{M93}.
Moreover,
the failure of the conservation law of energy  creates a major problem
to understand the global behavior in the stochastic case,
which motivates the recent numerical tracking of the energy in the works \cite{MRRY20, MRRY20.0}.
Another important qualitative difference is,
that the perturbation order of profiles
is of merely polynomial type in the stochastic case,
which makes it rather intricate to decouple different bubbles,
particularly, the remainders in the corresponding geometrical decomposition.
This is different from the previous construction of multi-bubble solutions
in \cite{M90} for NLS,
where the interactions are exponentially small in time.

In the present work,
in both dimensions one and two,
we are able to construct the multiple bubble blow-up solutions
concentrating at $K$ distinct points
to the rough nonlinear Schr\"odinger equation \eqref{equa-X-rough}, $1\leq K<\9$.
The constructed multi-bubble solutions behave asymptotically like
a sum of pseudo-conformal blow-up solutions
in the pseudo-conformal space $\Sigma$
near the blow-up time
and so have the pseudo-conformal blow-up rate $\sim (T-t)^{-1}$.
The upper bound of the approximation
is also obtained and, interestingly,
is closely related to that of the flatness of noise at blow-up points.
As a matter of fact,
if the noise is more flat at the blow-up points,
the approximation can be even taken in the more regular space $H^\frac 32$,
and the perturbation orders of the corresponding geometrical parameters
and the remainder can be also improved.

Another novelty of this work is concerned with the uniqueness of multi-bubble solutions.
The uniqueness issue is of significant importance in the classification
of blow-up solutions to dispersive equations.
In the remarkable paper \cite{M93},
Merle obtained the uniqueness of minimal mass blow-up solutions to NLS,
which states that the pseudo-conformal blow-up solution
is indeed the unique minimal mass blow-up solution up to the symmetries of NLS.
Such strong rigidity results were also obtained
by Rapha\"el and Szeftel \cite{R-S} for the inhomogeneous nonlinear Schr\"odinger equation,
and by Martel, Merle, Rapha\"el \cite{MMR15} for the mass-critical gKdV equation.
We also refer to \cite{MRS13} for the conditional uniqueness
result for the Bourgain-Wang solutions,
and \cite{KK20} for the Chern-Simons-Schr\"odinger equation.

However,
to the best of our knowledge,
there are very few results on the uniqueness of multi-bubble blow-up solutions to dispersive equations.

We prove that,
two multi-bubble blow-up solutions to equation \eqref{equa-X-rough} are indeed the same
if they have the same asymptotic blow-up profile
within the order $(T-t)^{3+\zeta}$, $\zeta>0$,
in the energy space $H^1$.
Hence,
in this asymptotic regime,
the $H^1$ multi-bubble blow-up solution
is exactly the above constructed solution
and so lies in the more regular pseudo-conformal space $\Sigma$.
This conditional uniqueness result of multi-bubbles solutions
can be also viewed as similar to the local uniqueness results in the elliptic setting,
see, e.g., \cite{CLL15,CPY20,DLY15}.

The existence and conditional uniqueness results are also obtained for a class of nonlinear Schr\"odinger equations
with lower order perturbations (see equations \eqref{equa-u-RNLS} and \eqref{equa-v-NLS-perturb} below),
particularly,
in the absence of the pseudo-conformal symmetry and the conservation of energy.

In particular,
the obtained results
are applicable to the single bubble case
and give the existence and conditional uniqueness of minimal mass blow-up solutions
for both the stochastic equation \eqref{equa-X-rough}
and the deterministic equation \eqref{equa-v-NLS-perturb}.

We would like to mention that,
the existence result is also applicable to the canonical deterministic NLS.
The positive frequencies $\{\omega_j\}_{j=1}^K$ in the construction here
are allowed to be arbitrarily small,
and
the asymptotic behavior can be taken in the pseudo-conformal space $\Sigma$
instead of the space $L^{2+\frac 4d}$,
which complement the previous results in \cite{M90}.
Furthermore,
the conditional uniqueness results
are new in both the stochastic and deterministic case.

The strategy of proof is mainly based on the modulation method
developed in \cite{R-S} for the minimal mass blow-up solutions
to the inhomogeneous nonlinear Schr\"odinger equation.
See also the recent work \cite{SZ20} in the stochastic setting.
One major difference here
is,
that the study of multi-bubble solutions
requires a delicate localization procedure.
A great effort is dedicated to
the decoupling of different bubbles.
Particularly,
because of the low polynomial type perturbation orders,
the decoupling of the remainders in the geometrical decomposition
is quite delicate.
Moreover,
a new generalized energy is introduced here,
it incorporates the localized functions in an appropriate way
such that different bubbles can be decoupled
and, simultaneously,
the key monotonicity property
keeps still preserved.
Let us also mention that,
the proof of the conditional uniqueness result requires an iterated argument,
which is also different from the single bubble case.
We expect the arguments developed here would be  also of interest
in the further understanding of multi-bubble solutions of other dispersive equations.   \\

{\bf Notations.}
For any $x=(x_1,\cdots,x_d) \in \bbr^d$
and any multi-index $\nu=(\nu_1,\cdots, \nu_d)$,
let
$|\nu|= \sum_{j=1}^d \nu_j$,
$\<x\>=(1+|x|^2)^{1/2}$,
$\partial_x^\nu=\partial_{x_1}^{\nu_1}\cdots \partial_{x_d}^{\nu_d}$,
and
$\<\na\>=(I-\Delta)^{1/2}$.

We use the standard Sobolev spaces
$H^{s,p}(\bbr^d)$, $s\in \bbr, 1\leq p\leq\9$.
In particular,
$L^p := H^{0,p}(\bbr^d)$ is
the space of $p$-integrable (complex-valued) functions,
$L^2$ denotes the Hilbert space endowed with the scalar product
$\<v,w\> =\int_{\bbr^d} v(x) \ol w(x)dx$,
and $H^s:= H^{s,2}$.
Let $\Sigma$ denote the pseudo-conformal space,
i.e., $\Sigma:=\{u\in H^1, |x|u\in L^2\}$.
As usual,
if $B$ is a Banach space,
$L^q(0,T;B)$ means the space of all integrable $B$-valued functions $f:(0,T)\to B$ with the norm
$\|\cdot\|_{L^q(0,T;B)}$,
and $C([0,T];B)$ denotes the space of all $B$-valued continuous functions on $[0,T]$ with the sup norm over $t$.
The local smoothing spaces is defined by
$L^2(I;H^\a_{\beta})=\{u\in \mathscr{S}': \int_{I} \int \<x\>^{2\beta}|\<\na\>^{\a} u(t,x)|^2  dxdt <\9 \}$,
$\a, \beta \in \mathbb{R}$.

Throughout this paper,
the positive constants $C$ and $\delta$ may change from line to line.

\section{Formulation of main results}  \label{Sec-Main}

\subsection{Main results}

To begin with,
let us first present the precise definition
of solutions to \eqref{equa-X-rough},
in which the noise term is taken in the sense of the controlled rough path.
For more details on the theory of rough paths,
see \cite{FH14,G04}.

\begin{definition} \label{def-X-rough}
We say that $X$ is a solution to \eqref{equa-X-rough} on $[0,\tau^*)$,
where $\tau^*\in (0,\9]$ is a random variable,
if $\bbp$-a.s. for any $\vf\in C_c^\9$,
$t \mapsto \<X(t), \vf\>$ is continuous on $[0,\tau^*)$
and for any $0<s<t<\tau^*$,
\begin{align*}
   \<X(t)-X(s), \vf\>
   - \int_s^t \<i X, \Delta\vf\>  + \<i|X|^{\frac 4d} X, \vf\>  - \< \mu X, \vf\> dr
   = \sum\limits_{k=1}^N \int_s^t \<i\phi_k X, \vf\> dB_k(r).
\end{align*}
Here,
the integral $\int_s^t \<i\phi_k X, \vf\> d B_k(r)$
is taken in the sense of controlled rough path
with respect to the rough paths $(B, \mathbb{B})$,
where $\mathbb{B}=(\mathbb{B}_{jk})$, $\mathbb{B}_{jk,st}:= \int_s^t \delta B_{j,sr} dB_k(r)$
with the integration taken in the sense of It\^o.
That is,
$\<i\phi_k X, \vf\> \in C^\a([s,t])$,
\begin{align} \label{phikX-st}
   \delta (\<i\phi_k X, \vf\>)_{st}
   = -\sum\limits_{j=1}^N \<\phi_j\phi_k X(s), \vf\> \delta B_{j,st}
     + \delta R_{k,st},
\end{align}
and
$\|\<\phi_j\phi_k X, \vf\> \|_{\a, [s,t]} <\9, \ \
   \|R_k\|_{2\a, [s,t]} <\9.
$
\end{definition}

We mention that,
because of the backward propagation procedure in the construction below,
the solution to \eqref{equa-X-rough} is no longer adapted
to the filtration $\{\mathscr{F}_t\}$.
Hence, the stochastic integration in \eqref{equa-X-rough} should be
interpreted in the sense of the controlled rough path,
instead of the It\^o sense.

The theory of (controlled) rough paths,
and the recent development of the theory of regularity structures \cite{H14}
and the para-controlled calculus \cite{GIP15}
have led to significant progress in solving
singular parabolic stochastic partial differential equations
with white noises.
We refer the interested readers  to  the monograph \cite{FH14}
and the references therein for
more details on these topics.

Throughout this paper we assume that
\begin{enumerate}
   \item[(A0)] ({\it Asymptotical flatness})
  For any multi-index $\nu \not = 0$ and $1\leq k\leq N$,
\begin{align} \label{decay}
   \lim_{|x|\to \9} \<x\>^2 |\partial_x^\nu \phi_k(x)| =0.
\end{align}

  \item[(A1)] ({\it Flatness at blow-up points})
  There exists $\nu_*\in \mathbb{N}$
  such that for every $1\leq k\leq N$ and $1\leq j\leq K$,
\begin{align} \label{degeneracy}
   \partial_x^\nu \phi_k(x_j)=0,\ \  \forall\ 0\leq|\nu|\leq \nu_*.
\end{align}
\end{enumerate}

\begin{remark}
The asymptotical flatness condition guarantees the
local well-posedness of equation \eqref{equa-X-rough} (see \cite{BRZ14,BRZ16}),
while the flatness condition at blow-up points
is mainly for the construction of multi-bubble solutions.
More interestingly,
the order $\nu_*$ is closely related to that of the
asymptotic behavior of solutions near the blow-up time.
See \eqref{X-Sj-H1} and \eqref{u-Sj-H1} below.
\end{remark}

For the frequencies $\omega_j>0$
and the blow-up points $x_j\in \bbr^d$,
$1\leq j\leq K$,
we mainly consider two cases below:

{\rm \bf Case (I).}
$\{x_j\}_{j=1}^K$ are arbitrary distinct points in $\bbr^d$,
and $\{\omega_j\}_{j=1}^K (\subseteq \bbr^+)$ satisfy that
for some $\omega>0$, $|\omega_j - \omega| \leq \ve$ for any $1\leq j \leq K$,
where $\ve >0$;

{\rm \bf Case (II).}
$\{\omega_j\}_{j=1}^K$ are arbitrary points in $\bbr^+$,
and $\{x_j\}_{j=1}^K (\subseteq \bbr^d)$ satisfy that
$|x_j-x_l| \geq \ve^{-1}$ for any $1\leq j \neq l \leq K$,
where $\ve>0$.

The existence of multiple bubble solutions is formulated in Theorem \ref{Thm-Multi-Blowup} below.

\begin{theorem} {\bf (Existence)} \label{Thm-Multi-Blowup}
Consider $d=1,2$.
Assume $(A0)$ and $(A1)$ with $\nu_*\geq 5$.
For every $1\leq K<\9$,
let
$\{\vartheta_{j}\} \subseteq \bbr$,
$\{x_j\}_{j=1}^K\subseteq \bbr^d$
be distinct points,
$\{\omega_j\}_{j=1}^K \subseteq \bbr^+$,
satisfying {\rm Case (I)} or {\rm Case (II)}.

Then,
for $\bbp-a.e.$ $\omega$
there exists $\ve^*(\omega) >0$ sufficiently small
such that for any $\ve\in (0,\ve^*]$,
there exists $\tau^*>0$ small enough
such that
for any $T\in (0,\tau^*(\omega)]$,
there exist $X_0(\omega)\in \Sigma$
and a corresponding blow-up solution $X(\omega)$ to \eqref{equa-X-rough},
satisfying that
for some $C>0$, $\zeta\in (0,1)$,
\begin{align} \label{X-Sj-H1}
\|e^{-W(t,\omega)} X(t,\omega)-\sum_{j=1}^KS_j(t)\|_{\Sigma}\leq C(T-t)^{\frac 12(\nu_*-5)+\zeta},\ \ t\in [0,T),
\end{align}
where $S_j$, $1\leq j\leq K$, are the pseudo-conformal blow-up solutions
\begin{align}  \label{Sj-blowup}
S_{j}(t,x)=(\omega_j(T-t))^{-\frac d2}Q (\frac{x-x_j}{\omega_j(T-t)})
             e^{ - \frac i 4 \frac{|x-x_j|^2}{T-t} + \frac{i}{\omega_j^2(T-t)} + i\vartheta_j}, \ \ t\in (0,T).
\end{align}
\end{theorem}

\begin{remark}
The asymptotic behavior \eqref{X-Sj-H1} yields that
the blow-up solution concentrates at the given $K$ points,
i.e.,
\begin{align}
|X(t,\omega)|^2\rightharpoonup \sum_{j=1}^K\|Q\|^2_{L^2}\delta_{x=x_j}, \ \ as\ t\to T.
\end{align}
In particular,
$\|X(t,\omega)\|_{L^2} = K\|Q\|_{L^2}$.
Hence the multi-bubble solutions
are constructed in the large mass regime,
which is different from the minimal mass case in \cite{SZ20}.
Moreover,
the asymptotic can be taken in the pseudo-conformal space $\Sigma$,
which improves the $H^1$-approximation result in \cite{SZ20}.
\end{remark}

\begin{remark}
The estimate \eqref{X-Sj-H1} also shows that
the order of approximation can be improved
if the noise is more flat at the blow-up points.
In the case $\nu_*\geq 6$,
the approximation \eqref{X-Sj-H1} can be even taken in the more regular space $H^\frac 32$
(see Proposition \ref{Prop-Rn-H23-bdd} below).
One can also improve the perturbation orders of the geometrical parameters
and the remainder with more flat noise,
see estimates \eqref{R-Tt}-\eqref{thetan-Tt}
and Theorem \ref{Thm-u-Unibdd} below.
Let us also mention that,
such asymptotic behavior \eqref{X-Sj-H1} is exhibited
only after the stochastic solutions are
rescaled by the random transformation $e^{-W}$.
\end{remark}

\begin{remark}\label{rek-exi}
The multi-bubble blow-up solutions were first constructed by Merle
in the pioneering work \cite{M90} for NLS in any dimensions,
the main blow-up profile in \cite{M90} is built on any functions
that decay exponentially,
while the frequencies $\{\omega_j\}_{j=1}^K$ are required to
have a uniform positive lower bound
and the asymptotic behavior is taken in the space $L^{2+\frac 4d}$.
In Theorem \ref{Thm-Multi-Blowup},
the multi-bubble solutions are constructed in dimensions one and two
and the blow-up profile is built on the ground state,
because the corresponding linearized operators are used in the construction.
The gain here is that,
in {\rm Case (I)} the frequencies are allowed to be arbitrarily small,
and in \eqref{X-Sj-H1} the approximation can be taken in the energy space $H^1$,
which is important in the proof of uniqueness result below.
\end{remark}

Our next main result is concerned with
the conditional uniqueness of multi-bubble solutions,
which is the content of Theorem \ref{Thm-Uniq-Blowup} below.

\begin{theorem} \label{Thm-Uniq-Blowup} {\bf (Conditional uniqueness)}
Consider $d=1,2$.
Assume $(A0)$ and $(A1)$ with $\nu_*\geq 5$.
For any $1\leq K<\9$,
let
$\{\vartheta_{j}\} \subseteq \bbr$,
$\{x_j\}_{j=1}^K\subseteq \bbr^d$
be distinct points,
$\{\omega_j\}_{j=1}^K \subseteq \bbr^+$,
satisfying {\rm Case (I)} or {\rm Case (II)}.

Then,
for $\bbp-a.e.$ $\omega$
there exists $\ve^*(\omega) >0$ sufficiently small
such that for any $\ve\in (0,\ve^*]$,
there exists $\tau^*>0$ small enough
such that
for any $T\in (0,\tau^*(\omega)]$,
there exists a unique blow-up solution $X(\omega)$ to \eqref{equa-X-rough}
satisfying
\begin{align} \label{X-Sj-Uniq}
\|e^{-W(t,\omega)} X(t,\omega)-\sum_{j=1}^KS_j(t)\|_{H^1}\leq C(T-t)^{3+\zeta},\ \ t\in[0,T),
\end{align}
where $C>0$, $\zeta\in (0,1)$,
and $S_j$ are the  pseudo-conformal blow-up solutions as in \eqref{Sj-blowup},
$1\leq j \leq K$.
\end{theorem}

\begin{remark}
Theorem \ref{Thm-Uniq-Blowup} states that
two multi-bubble solutions are the same if they both
have the asymptotic behavior \eqref{X-Sj-Uniq}
in the energy space $H^1$.
Moreover,
it also yields that any $H^1$ solution satisfying (2.7)
is exactly the constructed solution in Theorem 2.3,
which lies in the more regular $\Sigma$ space.
\end{remark}

\begin{remark}
The conditional uniqueness result also holds for the (deterministic) nonlinear Schr\"odinger equations
with lower order perturbations (see Theorem \ref{Thm-Uniq-Blowup-RNLS} and Remark \ref{Rem-deter-NLS-pert} below),
which include the canonical NLS equation.
Let us also mention that,
these conditional uniqueness results are new in both the stochastic and deterministic case.
\end{remark}

In particular,
in the special single bubble case (i.e., $K=1$),
we have the following
existence and conditional uniqueness of
minimal mass blow-up solutions.

\begin{theorem} \label{Thm-Minimass-Blowup}
Consider $d=1,2$.
Assume $(A0)$ and $(A1)$ with $\nu_*\geq 5$.
Let $x_*, \omega_*, \vartheta_*$
be any given points.
Then, for $\bbp$-a.e. $\omega$
there exists $\tau^*(\omega)>0$ sufficiently small,
such that
for any $T\in (0,\tau^*(\omega)]$
there exists a minimal mass blow-up solution $X(\omega)$  to \eqref{equa-X-rough}
satisfying that
\begin{align} \label{asymp-minimass}
\|e^{-W(t,\omega)} X(t,\omega)- S_*(t)\|_{\Sigma}\leq C(T-t)^{\frac 12(\nu_*-5)+\zeta},\ \ t\in [0,T),
\end{align}
where $C>0$, $\zeta\in (0,1)$,
and $S_*$ is as in \eqref{S}
with  $x_*, \omega_*, \vartheta_*$ replacing $\a_j, \omega_j, \vartheta_j$, respectively.

Moreover,
in the case where $\nu_*\geq 11$,
there exists a unique minimal mass blow-up solution $X(\omega)$ to \eqref{equa-X-rough}
satisfying that
\begin{align}
\|e^{-W(t,\omega)} X(t,\omega)- S_*(t)\|_{H^1}\leq C(T-t)^{3+\zeta},\ \ t\in [0,T).
\end{align}
\end{theorem}

\begin{remark}
The existence of minimal mass blow-up solutions are proved in the recent work \cite{SZ20},
but with the asymptotic \eqref{asymp-minimass} taken in the $H^1$ space.
The conditional uniqueness result is new in the stochastic case.
It should be mentioned that,
the strong uniqueness of minimal mass blow-up solutions to NLS in the deterministic case
was first obtained by Merle in the remarkable paper \cite{M93}.
Such strong rigidity results have been also obtained for the inhomogeneous NLS equations \cite{R-S}
and for the gKdV equations \cite{MMR15}.
For the stochastic equation \eqref{equa-X-rough}
the strong uniqueness of
minimal mass blow-up solutions
is at present still unclear,
due to the lack of the conservation law of energy.
\end{remark}

Equation \eqref{equa-X-rough} is indeed closely related to the
nonlinear Schr\"odinger equations with lower order perturbations.
More precisely,
by virtue of  the rescaling transformation
\begin{align} \label{X-u-rescal}
   u:=e^{-W}X,
\end{align}
we may reduce equation \eqref{equa-X-rough} to the
random equation below
\begin{align} \label{equa-u-RNLS}
   &i\partial_t u +\Delta u +|u |^{\frac{4}{d}}u +b \cdot \nabla u +c u =0,  \\
   &u(0)=u_{0},  \nonumber
\end{align}
where $b$ and $c$ are the coefficients of lower order perturbations
\begin{align}
 b(t,x)&=2\na W(t,x)= 2 i \sum\limits_{k=1}^N \na \phi_k(x)B_k(t),   \label{b} \\
 c(t,x)&= \sum\limits_{j=1}^d (\partial_j W(t,x))^2 + \Delta W(t,x) \nonumber \\
       &= - \sum\limits_{j=1}^d (\sum\limits_{k=1}^N \partial_j \phi_k(x)B_k(t))^2
          + i \sum\limits_{k=1}^N \Delta \phi_k(x) B_k(t). \label{c}
\end{align}
The key equivalent result has been proved in the recent work \cite{SZ20},
based on a delicate analysis of the temporal regularities.
Let us mention that,
such transformation is known as the Doss-Sussman transformation
in finite dimensional case,
and proves to be also very robust in the infinite dimensional spaces.
One main advantage is that,
from the viewpoint of analysis,
it enables one to
treat equation \eqref{equa-X-rough} as a random dynamic system
outside a uniform probability null set,
and thus to perform the sharp path-by-path analysis of stochastic solutions,
which is in general not easy by standard stochastic analysis.
Furthermore,
the rescaling approach also reveals the structure of stochastic equations,
which becomes more apparent
in the reduced nonlinear Schr\"odinger equations
with lower order perturbations.
See, for instance,
the stochastic logarithmic Schr\"odinger equations in \cite{BRZ17.0},
the damped effect of non-conservative noise in \cite{BRZ17},
and the scattering behavior in the stochastic setting in \cite{HRZ18}.
See also
the existence and geometrical characterization of optimal controllers
in  \cite{BRZ18,Z19}, related to the Ekeland's variational principle
and the theory of $U^p-V^p$ spaces.

The solutions to equation \eqref{equa-u-RNLS} are defined below.
\begin{definition} \label{def-u-RNLS}
We say that $u$ is a solution to \eqref{equa-u-RNLS} on $[0,\tau^*)$,
where $\tau^* \in (0,\9]$ is a random variable,
if $\bbp-a.s.$  $u\in C([0,\tau^*); H^1)$,
$|u|^{\frac{4}{d}}u\in L^1(0,\tau^*; H^{-1})$,
and $u$ satisfies
\begin{align} \label{equa-u-RNLS-def}
     u(t)= u(0) + \int_0^t i e^{-W(s)}\Delta(e^{W(s)}u(s))+i|u(s)|^{\frac{4}{d}}u(s) ds,  \ \
      t\in [0,\tau^*),
\end{align}
as an equation in $H^{-1}$.
\end{definition}

The key equivalent relationship
between equations \eqref{equa-X-rough} and \eqref{equa-u-RNLS}
is stated in Theorem \ref{Thm-Equiv-X-u} below.

\begin{theorem} (\cite[Theorem 2.10]{SZ20}) \label{Thm-Equiv-X-u}
$(i)$. Let $u$ be the solution to \eqref{equa-u-RNLS} on $[0,\tau^*)$
with $u(0) = u_0\in H^1$
in the sense of Definition \ref{def-u-RNLS},
where $\tau^* \in (0,\9]$ is a random variable.
Then, $\bbp$-a.s.,
$X := e^{W} u$ is the solution to equation \eqref{equa-X-rough} on $[0,\tau^*)$
with $X(0)= u_0$
in the sense of Definition \ref{def-X-rough}.

$(ii)$. Let $X$ be the solution to equation \eqref{equa-X-rough} on $[0,\tau^*)$ with $X(0)=X_0\in H^1$
in the sense of Definition \ref{def-X-rough},
satisfying that $\bbp$-a.s.
$\|X\|_{C([0,T]; H^1)} + \|X\|_{L^2(0,T; H^{\frac 32}_{-1})} <\9$,
$T\in (0,\tau^*)$,
and
$$\|e^{-it\Delta}e^{-W(t)}X(t) - e^{-is\Delta}e^{-W(s)}X(s)\|_{L^2} \leq C(t)(t-s),\ \ \forall 0\leq s<t<\tau^*.$$
Then,
$u:=e^{-W} X$ solves equation \eqref{equa-u-RNLS} on $[0,\tau^*)$ with $u(0)=X_0$
in the sense of Definition \ref{def-u-RNLS}.
\end{theorem}

Hence, by virtue of Theorem \ref{Thm-Equiv-X-u},
the proof of Theorems \ref{Thm-Multi-Blowup} and \ref{Thm-Uniq-Blowup}
is now reduced
to that of Theorems \ref{Thm-Multi-Blowup-RNLS} and \ref{Thm-Uniq-Blowup-RNLS} below
corresponding to the  equation \eqref{equa-u-RNLS}.

\begin{theorem} \label{Thm-Multi-Blowup-RNLS} {\bf (Existence)}
Consider $d=1,2$.
Assume $(A0)$ and $(A1)$ with $\nu_*\geq 5$.
For any $1\leq K <\9$,
let
$\{\vartheta_j\}_{j=1}^K \subseteq \bbr$,
$\{x_j\}_{j=1}^K\subseteq \bbr^d$
be distinct points,
and $\{\omega_j\}_{j=1}^K \subseteq \bbr^+$,
satisfying either {\rm Case (I)} or {\rm Case (II)}.

Then,
for $\bbp-a.e.$ $\omega$
there exists $\ve^*(\omega) >0$ sufficiently small
such that for any $\ve\in (0,\ve^*]$,
there exists $\tau^*>0$ small enough
such that
for any $T\in (0,\tau^*(\omega)]$,
there exist $u_0(\omega)\in \Sigma$ and a corresponding blow-up solution $u(\omega)$ to \eqref{equa-u-RNLS}
such that
\begin{align} \label{u-Sj-H1}
\|u(t,\omega)-\sum_{j=1}^KS_j(t)\|_{\Sigma}\leq C(T-t)^{\frac{1}{2}(\nu_*-5)+ \zeta},\ \ t\in [0,T),
\end{align}
where $C>0$, $\zeta\in (0, 1)$,
and $S_j$ are the pseudo-conformal blow-up solutions given by \eqref{Sj-blowup}, $1\leq j\leq K$.
\end{theorem}

\begin{theorem} \label{Thm-Uniq-Blowup-RNLS} {\bf (Conditional uniqueness)}
Consider the situations as in Theorem \ref{Thm-Multi-Blowup-RNLS}.
Then,
for $\bbp-a.e.$ $\omega$
there exists $\ve^*(\omega) >0$ sufficiently small
such that for any $\ve\in (0,\ve^*]$,
there exists $\tau^*>0$ small enough
such that
for any $T\in (0,\tau^*(\omega)]$,
there exists a unique blow-up solution $u(\omega)$ to \eqref{equa-u-RNLS}
satisfying
\begin{align} \label{u-Sj-Tt3-Uniq-RNLS}
\|u(t,\omega)-\sum_{j=1}^KS_j(t)\|_{H^1}\leq C(T-t)^{3+\zeta},\ \ t\in [0,T),
\end{align}
where $C>0$, $\zeta\in (0,1)$,
and $S_j$ are  as in \eqref{Sj-blowup}, $1\leq j\leq K$.
\end{theorem}

\begin{remark} \label{Rem-deter-NLS-pert}
The existence and conditional uniqueness results of multi-bubbles
also hold for the deterministic nonlinear Schr\"odinger equation
with lower order perturbations
if the Brownian motions $\{B_k\}$ in \eqref{equa-u-RNLS} are replaced by any deterministic continuous functions,
namely,
\begin{align} \label{equa-v-NLS-perturb}
    i\partial_t v +\Delta v +|v|^{\frac{4}{d}}v + a_1 \cdot \nabla u + a_2 u =0,
\end{align}
where
\begin{align*}
 a_1(t,x)= 2 i \sum\limits_{k=1}^N \na \phi_k(x)h_k(t),  \ \
 a_2(t,x)= - \sum\limits_{j=1}^d (\sum\limits_{k=1}^N \partial_j \phi_k(x)h_k(t))^2
          + i \sum\limits_{k=1}^N \Delta \phi_k(x) h_k(t),
\end{align*}
$\phi_k$ satisfy Assumptions $(A0)$ and $(A1)$
and $h_k \in C(\bbr^+; \bbr)$, $1\leq k\leq K$.
In particular,
these results are applicable to the canonical NLS equation,
in which $a_1, a_2 \equiv 0$.
Note that,
the standard pseudo-conformal symmetry and the conservation law of energy
are also destroyed in equation \eqref{equa-v-NLS-perturb}.
\end{remark}

\subsection{Strategy of proof}   \label{Subsec-strategy proof}

By virtue of the equivalent result  Theorem \ref{Thm-Equiv-X-u},
we shall mainly focus on the proof of Theorems \ref{Thm-Multi-Blowup-RNLS} and \ref{Thm-Uniq-Blowup-RNLS},
namely,
the existence and uniqueness of multi-bubble blow-up solutions
to nonlinear Schr\"odinger equations with lower order perturbations \eqref{equa-u-RNLS}.

As mentioned above,
the absence of specific symmetries  and the conservation law of energy
makes the blow-up analysis rather intricate.
A robust modulation method was developed by Rapha\"el and Szeftel \cite{R-S}
for the existence and uniqueness of minimal mass blow-up solutions
to inhomogeneous nonlinear Schr\"odinger equations,
which is a canonical model proposed by Merle \cite{M96} to break
the pseudo-conformal symmetry.
This modulation method has been recently applied in \cite{SZ20}
to construct minimal mass blow-up solutions
for both equations \eqref{equa-X-rough} and \eqref{equa-u-RNLS}.
The main strategy consists of
geometrical decompositions,
a bootstrap device
and backward propagation from the singularity.

Here, we use and extend the modulation method
to address the multi-bubble problem.
More specifically,
we first decompose the solution to \eqref{equa-u-RNLS}
into a main blow-up profile
and a remainder
\begin{align*}
    u(t,x)=  \sum_{j=1}^{K} \lbb_j^{-\frac d2} Q_j(t,\frac{x-\a_j}{\lbb_j}) e^{i\theta_j} +R(t,x), \ \
    with\ Q_j(t,y) = Q(y) e^{i(\beta_j\cdot y - \frac 14 \g_j|y|^2)},
\end{align*}
where $Q_j$ and $R$
satisfy the orthogonality conditions in \eqref{ortho-cond-Rn-wn} below,
which are related to the null space of the linearized operators around the ground state
and ensure the uniqueness of this decomposition.
Such geometrical decomposition enables us  to reduce the blow-up analysis
into those of the five finite-dimensional geometrical parameters
$(\lbb_j, \a_j, \beta_j, \g_j, \theta_j)$ and the infinite-dimensional remainder $R$.
As a first consequence,
the estimate of the modulation equations is derived,
which indeed captures the dynamics of the geometrical parameters.
This part is contained in Section \ref{Sec-Gem-Mod}.

The key unform estimates of the geometrical parameters
and the remainder are obtained by using a bootstrap device
and the propagation backward from the singularity.
The main efforts here
are dedicated to the analysis of
the localized mass,
the energy and
a new generalized functional.

It should be mentioned that,
unlike the single bubble case in \cite{R-S,SZ20},
the growth in the unstable direction $Q_j$
is analyzed via a localized mass,
instead of the usual whole mass.
Moreover,
we introduce a new generlized energy \eqref{def-I}
adapted to the multi-bubble setting.
It incorporates the localized functions in an appropriate way,
such that different bubbles
can be decoupled
while
the key monotonicity property keeps preserved.
One delicate problem here lies in the decoupling of the remainders,
due to the corresponding low polynomial type
perturbation orders,
and, actually,
extra decays have to be explored from the test functions.
Furthermore,
a refined estimate of the modulation parameter $\beta$
is derived from the coercivity
of the energy,
of which the proof requires a careful treatment
to balance the localized functions and the test functions
involved in the localized coercivity of the linearized operators.
These constitute the main part of Section \ref{Sec-Bootstrap}.

Then in Section \ref{Sec-Exit},
the construction of the multi-bubble blow-up solutions
follows from a compactness argument,
based on
the uniform estimates and integrating the flow backward from the blow-up time.
Let us mention that,
the uniform estimates can be also obtained in the pseudo-conformal space $\Sigma$,
which improves the previous $H^1$-estimate in \cite{SZ20}
and also simplifies the compactness argument.

Concerning the uniqueness part in Section \ref{Sec-Unique},
the key idea again relies on the monotonicity formula of the generalized energy
adapted to the difference of two multi-bubble solutions,
and is to show that the implied a priori bound of the difference
is exactly zero.

More precisely,
via the generalized energy \eqref{def-wtI} below,
we obtain the estimate (see Theorem \ref{Thm-Nt-Scal} below)
\begin{align*}
    \sup\limits_{t\leq s < T} D(s)
\leq C (\sum_{j=1}^{K} \sup\limits_{t\leq s < T} \frac{Scal_j(s)}{\lbb_j^2(s)}
+ \int_{t}^{T}\sum_{j=1}^{K}\frac{Scal_j(s)}{\lbb_j^3(s)}  + \ve^* \frac{D(s)}{T-s}ds).
\end{align*}
where
$D(t):=\|\nabla w(t)\|_{L^2}^2+\sum_{j=1}^{K} \lbb_j^{-2} \| w_j(t)\|_{L^2}^2$
is defined on the difference $w$, $w_j = w\Phi_j$ with the localized function $\Phi_j$,
and $Scal_j$ denotes the scalar products
of $w_j$ and the unstable directions in the null space.
This step requires a careful analysis of the differences between nonlinearities.

The next step is to control the unstable growth generated by the null space,
that is,
we  prove that (see Theorem \ref{Thm-Scal} below), for some $\zeta\in (0,1)$,
\begin{align*}
Scal_j(t)\leq C(T-t)^{2+\zeta} \sup\limits_{t\leq s<T}D(s).
\end{align*}
For this purpose,
a new renormalized variable is introduced.
It satisfies a neat formulation of equation
and enables us to obtain the estimates of the scalar products in $Scal_j$
in a simplified diagonalized form (see Proposition \ref{Prop-ej-Lkernel} below).

At this stage,
by virtue of the two estimates above,
we obtain the estimate of $D(t)$ in a closed form.
It should be mentioned that,
because of the localization functions in the multi-bubble case,
an extra error $ \ve^* \frac{D(t)}{T-t}$ is also involved here.
This requires an iterated argument
to show that $D(t)$ is exactly zero,
which is different from the single bubble case in \cite{R-S}.

The remainder of this paper is structured as follows.
In Section \ref{Sec-Pre} we first present some preliminaries
including the localization functions,
the coercivity of linearized operator,
and the Taylor expansion in the complex situation.
Then, Sections \ref{Sec-Gem-Mod}, \ref{Sec-Bootstrap} and \ref{Sec-Exit}
are mainly devoted to the proof of the existence of blow-up solutions
at multiple points.
In Section \ref{Sec-Unique} we prove the uniqueness of blow-up solutions.
Finally, some technical proofs are postponed to the Appendix, i.e., Section \ref{Sec-App}.

\section{Preliminaries}  \label{Sec-Pre}

\subsection{Localization}  \label{Subsec-Local}

We shall use the localization functions
in order to construct the blow-up profiles concentrating at
distinct blow-up points.

For this purpose,
we note that,
because equation \eqref{equa-u-RNLS} is invariant under orthogonal transforms,
we may take an orthonormal basis $\{{\bf v_j}\}_{j=1}^d$ of $\bbr^d$,
such that $(x_j-x_l)\cdot {\bf v_1}\neq0$ for any $1\leq j\neq l\leq K$.
Hence,
we may assume that $x_1\cdot {\bf v_1}<x_2\cdot {\bf v_1}<\cdots<x_K\cdot {\bf v_1}$.
Then,
\begin{align} \label{sep-xj-0}
\sigma :=\frac{1}{12}\min_{1\leq j \leq K-1}\{(x_{j+1}-x_j)\cdot {\bf v_1}\}> 0.
\end{align}

Let $\Phi(x)$ be a smooth function on $\R^d$ such that $0\leq \Phi(x)\leq 1$,
$|\na \Phi(x)| \leq C \sigma^{-1}$,
$\Phi(x)=1$ for $x\cdot {\bf v_1}\leq 4\sigma$
and $\Phi(x)=0$ for $x\cdot {\bf v_1} \geq 8\sigma$.
The localization functions $\Phi_j$ are defined by
namely,
\be\ba \label{phi-local}
&\Phi_1(x) :=\Phi(x-x_1), \ \ \Phi_K(x) :=1-\Phi(x-x_{K-1}),  \\
&\Phi_j(x) :=\Phi(x-x_{j})-\Phi(x-x_{j-1}),\ \ 2\leq j\leq K-1.
\ea\ee
In particular,
we have the  partition of unity $1= \sum_{j=1}^K \Phi_j$.

Lemma \ref{Lem-inter-est} below enables us to decouple different blow-up profiles
and will be used frequently throughout this paper.

\begin{lemma}(Interaction estimates) \label{Lem-inter-est}
Let $0<t^*<T_*<T<\9$.
For every $1\leq j\leq K$,
set
\begin{align} \label{Gj-gj-g}
    G_j(t,x)
    := \lbb_j^{-\frac d2} g_j(t,\frac{x-\a_j}{\lbb_j}) e^{i\theta_j}, \ \
    with\ \ g_{j}(t,y) := g(y) e^{i(\beta_{j}(t) \cdot y - \frac 14\g_{j}(t) |y|^2)},
\end{align}
where $g \in C_b^2(\bbr^d)$
decays exponentially fast at infinity
\begin{align*}
   |\partial^\nu g(y)| \leq C e^{-\delta |y|}, \ \ |\nu|\leq 2,
\end{align*}
with $C,\delta>0$,
for $1\leq j\leq K$,
$\calp_j:=(\lbb_j,\alpha_j,\beta_j,\gamma_j,\theta_j) \in C([t^*,T_*]; \bbr^{2d+3})$
 satisfies
\begin{align} \label{aj-xj}
  |(\alpha_j(t)-x_j)\cdot {\bf v_1}|\leq \sigma, \
  |x_j - \a_j(t)| \leq 1,\
  \half\leq\frac{\lambda_{j}(t)}{|\omega_{j}(T-t)|}\leq2, \ \ t\in[t^*,T_*],
\end{align}
and $|\beta_j|+|\g_j|\leq 1$,
\begin{align} \label{T-M-0}
    C T (1+ \max_{1\leq j\leq K} |x_j|) \leq  1,
\end{align}
where $C$ is sufficiently large but independent of $T$.
Then,
there exists $\delta>0$ such that
for any $1\leq j\not =l\leq K$,
$m\in \mathbb{N}$,
and for any multi-index $\nu$ with $|\nu|\leq 2$,
\begin{align} \label{Gj-Gl-decoup}
  \int\limits_{\bbr^d} |x-\a_l|^n |\partial^\nu G_l(t)| |x-\a_j|^m |G_j(t)| dx
   \leq Ce^{-\frac{\delta}{T-t}}, \ \ t\in [t^*, T_*].
\end{align}
Moreover, for any $h\in L^1$ or $L^2$,
$1\leq j\not = l\leq K$,
$m,n\in \mathbb{N}$,
and for any multi-index $\nu$ with $|\nu|\leq 2$,
\begin{align}  \label{Gj-hl-decoup}
   \int\limits_{\bbr^d} |x-\a_l|^n |\partial^\nu G_l(t)| |x-\a_j|^m |h| \Phi_jdx
   \leq Ce^{-\frac{\delta}{T-t}} \min\{\|h\|_{L^1}, \|h\|_{L^2}\}, \ \ t\in [t^*,T_*].
\end{align}
\end{lemma}

The proof is postponed to the Appendix for simplicity.
Lemma \ref{Lem-inter-est} actually shows that the interactions between $\{U_j\}$
and other profiles are very weak,
mainly due to the exponential decay of the ground state.

\subsection{Coercivity of linearized operators} \label{Subsec-Coer}

We denote  $Q$ the ground state that solves  the soliton equation \eqref{equa-Q}.
It follows from \cite[Theorem 8.1.1]{C} that
$Q$ is smooth and decays at infinity exponentially fast,
i.e., there exist $C, \delta>0$ such that for any multi-index $|\nu|\leq 3$,
\be\label{Q-decay}
|\partial_x^\nu Q(x)|\leq C e^{-\delta |x|}, \ \ x\in \bbr^d.
\ee

Let $L=(L_+,L_-)$ be the linearized operator around the ground state,
defined by
\begin{align} \label{L+-L-}
     L_{+}:= -\Delta + I -(1+{\frac{4}{d}})Q^{\frac{4}{d}}, \ \
    L_{-}:= -\Delta +I -Q^{\frac{4}{d}}.
\end{align}
The generalized null space of operator $L$ is
spanned by $\{Q, xQ, |x|^2 Q, \na Q, \Lambda Q, \rho\}$,
where
$\Lambda := \frac{d}{2}I + x\cdot \na$,
and $\rho$ is the unique $H^1$ spherically symmetric solution to the equation
\begin{align} \label{def-rho}
L_{+}\rho= - |x|^2Q,
\end{align}
which satisfies the exponential decay property (see, e.g., \cite{K-M-R, MP18}),
i.e., for some $C,\delta>0$,
\ben
|\rho(x)|+|\nabla \rho(x)|
\leq Ce^{-\delta|x|}.
\enn
Moreover, we have (see, e.g., \cite[(B.1), (B.10), (B.15)]{Wenm})
\be \ba \label{Q-kernel}
&L_+ \na Q =0,\ \ L_+ \Lambda Q = -2 Q,\ \ L_+ \rho = -|x|^2 Q, \\
&L_{-} Q =0,\ \ L_{-} xQ = -2 \na Q,\ \ L_{-} |x|^2 Q = - 4 \Lambda Q.
\ea\ee

For any complex valued $H^1$ function $f = f_1 + i f_2$
in terms of the real and imaginary parts,
we set
\be
(Lf,f) :=\int f_1L_+f_1dx+\int f_2L_-f_2dx.
\ee
Let $\mathcal{K}$ denote the
set of all complex valued $H^1$ functions $f=f_1+if_2$
satisfying the orthogonality conditions below
\be\ba\label{ortho-cond}
&\int Qf_1dx=0,\;\; \int xQf_1dx=0,\;\;\int |x|^2Qf_1dx=0,\\
&\int \nabla Qf_2dx=0,\;\;\int \Lambda Qf_2dx=0,\;\;
\int \rho f_2dx=0.
\ea\ee

The coercivity property below is crucial in the proof of main results.

\begin{lemma} {(\cite[Theorem 2.5]{Wenm})}  \label{Lem-coerc-L}
There exists $C>0$ such that
\begin{align}
(Lf,f)\geq C \|f\|_{H^1}^2,\ \ \forall f\in \mathcal{K}.
\end{align}
\end{lemma}

We define the scalar products along all the unstable directions
in the null space
\begin{align} \label{Scal-def}
Scal(f)=\<f_1,Q\>^2+\<f_1,xQ\>^2+\<f_1,|x|^2Q\>^2+\<f_2,\nabla Q\>^2+\<f_2,\Lambda Q\>^2+\<f_2,\rho\>^2,
\end{align}
where $f=f_1 + if_2 \in H^1$.
As a consequence of Lemma \ref{Lem-coerc-L} we have
\begin{corollary}(\cite[Corollary 3.2]{SZ20})   \label{Cor-coer-f-H1}
There exist positive constants $C_1, C_2>0$,  such that
\begin{align} \label{coer-f-H1}
(Lf,f)\geq&  C_1\|f\|_{H^1}^2
             -C_2Scal(f),\ \ \forall f \in H^1,
\end{align}
where $f_1$ and $f_2$ are the real and imaginary parts of $f$, respectively.
\end{corollary}

\begin{corollary}(Localized coercivity)   \label{Cor-coer-f-local}
Let $\phi$ be a positive smooth radial function on $\R^d$,
such that
$\phi(x) = 1$ for $|x|\leq 1$,
$\phi(x) = e^{-|x|}$ for $|x|\geq 2$,
$0<\phi \leq 1$,
and $\lf|\frac{\nabla\phi}{\phi}\rt|\leq C$ for some $C>0$.
Set $\phi_A(x) :=\phi\lf(\frac{x}{A}\rt)$, $A>0$.
Then,
for $A$ large enough we have
\begin{align} \label{coer-f-local}
\int (|f|^2+|\nabla f|^2)\phi_A -(1+\frac 4d)Q^{\frac4d}f_1^2-Q^{\frac4d}f_2^2dx\geq C_1\int(|\nabla f|^2+|f|^2)\phi_A dx-C_2Scal(f),
\end{align}
where $C_1, C_2>0$, and $f_1, f_2$ are the real and imaginary parts of $f$, respectively.
\end{corollary}

The proof of Corollary \ref{Cor-coer-f-local} is  similar to that of
\cite[Corollary 3.3]{SZ20}
and is postponed to the Appendix for simplicity.

\subsection{Expansion of the nonlinearity} \label{Subsec-Expan}

We shall use the notations that,
for any continuous differentiable function $g:\mathbb{C}\to \mathbb{C}$
and for any $v, R\in \mathbb{C}$,
\begin{align}
  g'(v,R) \cdot R
     :=&  R  \int_0^1 \partial_z g(v+sR) ds
           + \ol R  \int_0^1 \partial_{\ol z} g(v+sR) ds,  \label{f'R} \\
  g''(v,R) \cdot R^2
     :=&   R^2 \int_0^1 t  \int_0^1 \partial_{zz}g(v+st R) dsdt
      + 2|R|^2 \int_0^1 t \int_0^1 \partial_{z \ol z}g(v+st R) dsdt \nonumber \\
      &   + \ol R^2 \int_0^1 t  \int_0^1 \partial_{\ol z \ol z}g(v+st R) dsdt,    \label{f''R2}
\end{align}
where $\partial_z g$ and $\partial_{\ol{z}} g$
are the usual complex derivatives $\partial_z g= \frac 12(\partial_x g - i \partial_y g)$,
$\partial_{\ol z} g= \frac 12(\partial_x g + i \partial_y g)$,
respectively.
Then,
one has (see, e.g., \cite[(3.10)]{KV13})
\begin{align} \label{g-gz-expan}
   g(v+R) = g(v) + g'(v, R)\cdot R.
\end{align}
Moreover,
if  $\partial_z g$ and $\partial_{\ol{z}} g$  are also continuously differentiable,
we may expand $g$ up to the second order
\begin{align} \label{g-gzz-expan}
   g(v+R)
   =& g(v)
     +  \partial_z g(v) R
     + \partial_{\ol{z}} g(v)  \ol R
     + g''(v, R) \cdot R^2.
\end{align}

In particular,
for the complex function
$f(z)=|z|^{\frac 4d} z$ with $d=1,2$,
we have
\begin{align}
   |f'(v, R)\cdot R| \leq& C (|v|^{\frac 4d} + |R|^\frac 4d) |R|, \label{f'vR-bdd} \\
   |f''(v, R)\cdot R^2| \leq& C (|v|^{\frac 4d-1} + |R|^{\frac 4d-1}) |R|^2. \label{f''vR2-bdd}
\end{align}

It would be also useful to use the expansion,
for  $f(z)=|z|^{\frac 4d} z$ with $d=1,2$,
\begin{align} \label{f-Taylor}
   f(v+R) =&  f(v) + f'(v)\cdot R
         +  f''(v)\cdot R^2
         + \calo(\sum\limits_{k=3}^{1+\frac 4d} |v|^{1+\frac 4d-k} |R|^k),
\end{align}
where
\begin{align}
   f'(v)\cdot R :=&  \partial_z f(v) R+ \partial_{\ol{z}} f(v) \ol{R}
                 = (1+\frac 2d) |v|^{\frac 4d} R + \frac 2d |v|^{\frac 4d-2} v^2 \ol{R},   \label{f-linear} \\
   f''(v)\cdot R^2
                :=& \frac 12 \partial_{zz}f(v)R^2 + \partial_{z\ol{z}} f(v) |R|^2 + \frac 12 \partial_{\ol{z}\ol{z}} f(v) \ol{R}^2 \nonumber \\
                 =& \frac 1d(1+\frac 2d)|v|^{\frac 4d -2} \ol{v} R^2
                    + \frac 2d (1+\frac 2d) |v|^{\frac 4d -2} v |R|^2
                    +   \frac 1d(\frac 2d-1) |v|^{\frac 4d -4} v^3 \ol{R}^2.  \label{f-quadratic}
\end{align}

Similarly, for $F(z):= \frac{d}{2d+4} |z|^{2+ \frac 4d}$ with $d=1,2$,
we have the expansion
\begin{align} \label{F-Taylor*}
  F(u)=& F(v)
         + \frac 12 |v|^{\frac 4d}\ol{v} R
         + \frac 12  |v|^{\frac 4d} v \ol{R} \nonumber \\
       & + \frac{1}{2d} |v|^{\frac 4d -2}\ol{v}^2 R^2
         + \frac 12 (1+\frac 2d) |v|^{\frac 4d} |R|^2
         + \frac{1}{2d} |v|^{\frac 4d -2} v^2 \ol{R}^2
         + \calo(\sum\limits_{k=3}^{2+\frac 4d} |v|^{2+\frac 4d -k} |R|^k).
\end{align}

In most cases in this paper,
the high order terms in the expansion of nonlinearity
can be controlled by the Gagliardo-Nirenberg  inequality
contained in Lemma \ref{Lem-GN} below.

\begin{lemma} \label{Lem-GN}
(\cite[Theorem 1.3.7]{C})
Let $d\geq 1$ and $2\leq p< \infty$.
Then, there exists $C>0$ such that
\begin{align}  \label{G-N}
\|f\|_{L^p}\leq C \|f\|_{L^2}^{1 - d (\frac 12-\frac 1p)} \|\nabla f\|_{L^2}^{d(\frac 12-\frac 1p)}, \ \ \forall f\in H^1.
\end{align}
In particular,
for any $1<p<\9$,
\begin{align} \label{fLp-fH1}
   \|f\|_{L^p} \leq C \|f\|_{H^1},\ \ \forall f\in L^p.
\end{align}
\end{lemma}

We also have the product rule below.
\begin{lemma}  \label{Lem-Pro-rule} (Product rule \cite[p.105, Proposition 1.1]{T00})
For any $s>0$,
\begin{align} \label{fra-lei-law}
   \|uv\|_{H^{s,p}}
   \leq C (\|u\|_{L^{q_1}} \|v\|_{H^{s,q_2}}
         + \|v\|_{L^{r_1}} \|u\|_{H^{s,r_2}}),
\end{align}
where
$\frac 1p = \frac{1}{q_1}+\frac{1}{q_2}=\frac{1}{r_1}+\frac{1}{r_2}$,
$q_1,r_1 \in (1,\9]$,
$q_2, r_2 \in (1,\9)$.
\end{lemma}

As a consequence we have
\begin{lemma} \label{Lem-H12-H1}
For any complex functions $f,g,h$
and for any $l,m,n\in \mathbb{N}$,
we have
\begin{align} \label{fghLp-H1}
   \|f^l g^m h^n\|_{L^p}
   \leq C \|f\|_{H^1}^l \|g\|^m_{H^1} \|h\|^n_{H^1}.
\end{align}
Moreover, we also have
\begin{align} \label{fghH12-H1}
   \|f^l g^m h^n\|_{\dot H^{\frac 12}}
   \leq C  \|f\|_{H^1}^l \|g\|^m_{H^1} \|h\|^n_{H^1}.
\end{align}
\end{lemma}

{\bf Proof.}
\eqref{fghLp-H1} follows from H\"older's inequality and \eqref{fLp-fH1}.
Regarding \eqref{fghH12-H1},
by the product rule,
\begin{align*}
   \|f^l g^m h^n\|_{\dot H^{\frac 12}}
   \leq& C (\|f\|_{H^{\frac 12, p_{11}}} \|f\|^{l-1}_{L^{p_{12}}} \|g\|^m_{L^{q_1}} \|h\|^n_{L^{r_1}}
        + \|f\|^{l}_{L^{p_2}} \|g\|_{H^{\frac 12,q_{21}}} \|g\|^{m-1}_{L^{q_{22}}} \|h\|^n_{L^{r_2}}  \nonumber \\
        &\quad      + \|f\|^{l}_{L^{p_2}}   \|g\|^m_{L^{q_3}} \|h\|^n_{H^{\frac 12,r_{31}}} \|h\|^{n-1}_{L^{r_{32}}}),
\end{align*}
where
$\frac 12 = \frac {1}{p_{11}} + \frac{l-1}{p_{12}} + \frac{m}{q_1} + \frac{n}{r_1}
= \frac{l}{p_2} + \frac {1}{q_{21}} + \frac{m-1}{q_{22}} + \frac{n}{r_2}
= \frac{l}{p_2} + \frac {m}{q_3} + \frac{1}{r_{31}} + \frac{n-1}{r_{32}}$.
We then take $p_{11},q_{21},r_{31}$ close to $2$ such that
$H^1$ is imbedded into the Sobolev spaces $H^{\frac 12, p_{11}}$,
$H^{\frac 12, q_{21}}$ and $H^{\frac 12, r_{31}}$.
Then, taking into account \eqref{fLp-fH1}
we obtain \eqref{fghH12-H1} and finish the proof.
\hfill $\square$

\section{Geometrical decomposition and modulation equations}  \label{Sec-Gem-Mod}

\subsection{Geometrical decomposition}  \label{Subsec-Geom-Decomp}
For each $1\leq j\leq K$,
define the modulation parameters by
$\mathcal{P}_j:=(\lbb_j,\alpha_j,\beta_j,\gamma_j,\theta_j) \in \bbr^{2d+3}$,
where $\lbb_j \in \bbr$,
$\a_j \in \bbr^d$,
$\beta_j \in \bbr^d$, $\gamma_j \in \bbr$ and
$\theta_j \in \bbr$,
and set
$P_{j}:=   |\lbb_{j}|+|\a_{j}-x_j| + |\beta_{j}| + |\g_{j}|$,
where $x_j$ are the given blow-up points,
$1\leq j\leq K$.

We also set
$\calp:= (\calp_1, \cdots, \calp_K) \in \bbr^{(2d+3)K}$,
$P:= \sum_{j=1}^K P_j$.
Similarly,
let $\lbb := (\lbb_1, \cdots, \lbb_K) \in \bbr^K$
and
$|\lbb|:=\sum_{j=1}^{K}|\lbb_{j}|$.
Similar notations are also used
for the remaining parameters.

\begin{proposition} (Geometrical decomposition)   \label{Prop-dec-un}
Assume that $u\in C([\wt t, T_*]; H^1)$ for some $\wt t \in [0,T_*)$
and $u(T_*)=S_T(T_*)$.
Then,
for $T_*$ sufficiently close to $T$,
there exist $t^*<T_*$
and
unique modulation parameters
$\mathcal{P}
\in C^1((t^*,T_*); \bbr^{(2d+3)K})$,
such that
$u$
can be geometrically decomposed
into the main blow-up profile and the remainder
\begin{align} \label{u-dec}
    u(t,x)=U(t,x)+R(t,x),\ \ t\in [t^*, T_*],\ x\in \bbr^d,
\end{align}
where the main blow-up profile
\begin{align}  \label{U-Uj}
    U(t,x)
    =   \sum_{j=1}^{K} U_j(t,x),
\end{align}
with
\begin{align} \label{Uj-Qj-Q}
   U_j(t,x) = \lbb_j^{-\frac d2} Q_j(t,\frac{x-\a_j}{\lbb_j}) e^{i\theta_j},  \ \
   Q_j(t,y) = Q(y) e^{i(\beta_j\cdot y - \frac 14 \g_j|y|^2)},
\end{align}
and $R(T_*)=0$,  the modulation parameters satisfy
\begin{align} \label{PjT*}
   \calp_j(T_*)=(\omega_j (T-T_*), x_j, 0, \omega^2_j (T-T_*), \omega^{-2}_j (T-T_*)^{-1} +\vartheta_j).
\end{align}
Moreover,
for each $1\leq j\leq K$, the following orthogonality conditions hold on $[t^*,T_*] $:
\be\ba\label{ortho-cond-Rn-wn}
&{\rm Re}\int (x-\a_{j}) U_{j}(t)\ol{R}(t)dx=0,\ \
{\rm Re} \int |x-\a_{j}|^2 U_{j}(t) \ol{R}(t)dx=0,\\
&{\rm Im}\int \nabla U_{j}(t) \ol{R}(t)dx=0,\ \
{\rm Im}\int \Lambda U_{j}(t) \ol{R}(t)dx=0,\ \
{\rm Im}\int \varrho_{j}(t) \ol{R}(t)dx=0,
\ea\ee
where
\be \label{rhon}
   \varrho_{j}(t,x)= \lbb^{-\frac d2}_{j} \wt \rho_{j}(t,\frac{x-\a_{j}}{\lbb_n}) e^{i\theta_{j}}
 \ \ with\  \ \wt \rho_{j}(t,y) := \rho(y)^{i(\beta_{j}(t)\cdot y - \frac 14 \g_{j}(t) |y|^2)},
\ee
and $\rho$ is given by \eqref{def-rho}.
\end{proposition}

\begin{remark}
Proposition \ref{Prop-dec-un}
is actually a local version of the geometrical decomposition as $t^*$ may depend on $T_*$.
However,
as we shall see later,
by virtue of the bootstrap estimates in Theorem \ref{Thm-u-Boot} below
we indeed have the global geometrical decomposition
on the time interval $[0,T_*]\subseteq [0,T)$
if $T$ is sufficiently small.
See Theorem \ref{Thm-u-Unibdd} below.
\end{remark}

The proof of Proposition \ref{Prop-dec-un}
is quite similar to that of \cite[Proposition 4.1]{SZ20}
and is mainly based on the implicit function theorem.
We mention that,
the computations of the Jacobian of transformation
also include the interactions between different profiles $\{U_j\}_{j=1}^K$
which, however,
by Lemma \ref{Lem-inter-est},
only contribute exponentially small errors
due to the exponential decay of the ground state.
Thus the Jacobian is still non-zero.
The details are omitted here for simplicity.

\subsection{Modulation equations} \label{Subsec-Mod-Equa}
Let $\dot{g}:= \frac{d}{dt}g$  for any $C^1$ function $g$.
For each $1\leq j\leq K$,
define \emph{the vector of modulation equations} by
\begin{align} \label{Mod-def}
   Mod_{j}:= |\lambda_{j}\dot{\lambda}_{j}+\gamma_{j}|+|\lambda_{j}^2\dot{\gamma}_{j}+\gamma_{j}^2|
   +|\lambda_{j}\dot{\alpha}_{n,j}-2\beta_{j}|
                +|\lambda_{j}^2\dot{\beta}_{j}+\gamma_{j}\beta_{j}|
              +|\lambda_{j}^2\dot{\theta}_{j}-1-|\beta_{j}|^2|,
\end{align}
and set $Mod:=\sum_{j=1}^{K}Mod_{j}.$

The main result
of this subsection
is formulated in Proposition \ref{Prop-Mod-bdd}
below.

\begin{proposition} \label{Prop-Mod-bdd}
Assume that $u$ has the geometrical decomposition on $[t^*,T_*] \subseteq [0,T)$ as in \eqref{u-dec}
with the modulation parameters $\calp=(\lbb, \a, \beta, \g, \theta)\in \bbr^{(2d+3)K}$,
and
\begin{align*}
  C_1(T-t) \leq |\lbb(t)| \leq C_2(T-t), \ \ t\in [t^*, T_*],
\end{align*}
where $C_1, C_2>0$.
Then, for $T$ small enough
and for $t^*$ close to $T_*$,
we have for any $ t\in[t^*,T_*]$,
\begin{align} \label{Mod-bdd}
Mod(t)\leq
C( \sum_{j=1}^{K} |{\rm Re}\langle R_j, U_j\rangle|+  P^2(t) \|R(t)\|_{L^2}+\|R(t)\|_{L^2}^2
+\|R(t)\|_{H^1}^3   + P^{\nu_*+1}(t) + e^{-\frac{\delta}{T-t}}),
\end{align}
where $C>0$ and $\nu_*$ is the index of flatness in Assumption $(A1)$.
\end{proposition}

The proof relies mainly
on the analysis of the equation of remainder $R$,
the almost orthogonality of profiles $U_j$ and $R_j$,
and the decoupling of different blow-up profiles $U_j$ and $U_l$, $j\not = l$.

To be precise,
we use  the partition of unity  $1=\sum_{j=1}^{K} \Phi_j$
to get
\begin{align}  \label{R-Rj}
R=\sum_{j=1}^{K}R_{j}, \ \ with\ \ R_j:= R \Phi_j.
\end{align}
Define the renormalized remainder $\ve_{j}$ by
\begin{align} \label{Rj-ej}
R_{j}(t,x) =\lbb_{j}^{-\frac d2} \varepsilon_{j} (t,\frac{x-\a_{j}}{\lbb_{j}}) e^{i\theta_{j}}.
\end{align}
Then, by \eqref{equa-u-RNLS} and \eqref{u-dec},
the remainder $R$ satisfies the equation
\begin{align}\label{eq-U-R}
 &i\partial_tR
   +\sum_{l=1}^{K}
   (\Delta R_{l}+(1+\frac 2d)|U_{l}|^{\frac{4}{d}}R_{l}
   +\frac 2d|U_{l}|^{\frac{4}{d}-2}U_{l}^2 \ol{R_{l}}
   +  i\partial_tU_{l}+\Delta U_{l}+|U_{l}|^{\frac{4}{d}}U_{l}) \nonumber \\
=& -H_1-H_2-f''(U,R)\cdot R^2
   -  \sum\limits_{l=1}^K(b\cdot\nabla  (U_l+R_l)  +c  (U_l+R_l) )
\end{align}
Here,
$H_1, H_2$ contain the interactions between different blow-up profiles
\begin{align}
  H_1 : =& (1+\frac 2d)|U|^{\frac{4}{d}}R+\frac 2d|U|^{\frac{4}{d}-2}U^2\ol{R}
-\sum_{l=1}^{K}((1+\frac 2d)|U_{l}|^{\frac{4}{d}}R_{l}+\frac 2d|U_{l}|^{\frac{4}{d}-2}U_{l}^2\ol{{R}_{l}}), \label{G} \\
  H_2 : =&  |U|^{\frac{4}{d}}U-\sum_{l=1}^{K} |U_{l}|^{\frac{4}{d}}U_{l}, \label{H}
\end{align}
and $f''(U,R)\cdot R^2$ is defined as in \eqref{f''R2}
with $f,U$ replacing $g$ and $v$, respectively.

Using \eqref{Uj-Qj-Q}
we have
\begin{align}   \label{equa-Ut}
 i\partial_tU_{j}+&\Delta U_{j}+|U_{j}|^{\frac{4}{d}}U_{j}
=\frac{e^{i\theta_{j}}}{\lambda_{j}^{2+\frac d2}}
\big\{-(\lambda_{j}^2\dot{\theta}_{j}-1-|\beta_{j}|^2)Q_{j}
      -(\lambda_{j}^2\dot{\beta}_{j}+\gamma_{j}\beta_{j})\cdot yQ_{j}  \nonumber \\
   &   +\frac 14 (\lambda_{j}^2\dot{\gamma}_{j}+\gamma_{j}^2) |y|^2 Q_{j}
  -i(\lambda_{j}\dot{\alpha}_{j}-2\beta_{j})\cdot \nabla Q_{j}
     -i(\lambda_{j}\dot{\lambda}_{j}+\gamma_{j})\Lambda Q_{j} \big\}(t,\frac{x-\alpha_{j}}{\lambda_{j}}).
\end{align}

The important fact here is,
that the modulation equations show up on the right-hand side of \eqref{equa-Ut}
as the coefficients of the five directions
in the generalized null space of the operator $L$ defined in \eqref{L+-L-}.
This enables us to extract each modulation equation
by applying the almost orthogonality in Lemma \ref{Lem-almost-orth} below,
which in turn follows from  Lemma \ref{Lem-inter-est}
and the orthogonality conditions in \eqref{ortho-cond-Rn-wn}.

\begin{lemma} (Almost orthogonality)  \label{Lem-almost-orth}
Let $t^*$ be as in Proposition \ref{Prop-Mod-bdd}.
Then, for $t^*$ close to $T$,
there exists $\delta >0$ such that for any
$1\leq j\leq K$,
it holds on $[t^*,T]$ that
\be\ba \label{Orth-almost}
&|{\rm Re}\int (x-\a_{j}) U_{j} \ol{R_{j}}dx|
  + |{\rm Re} \int |x-\a_{j}|^2 U_{j} \ol{R_{j}}dx|\leq Ce^{-\frac{\delta}{T- t}}\|R\|_{L^2},  \\
& |{\rm Im}\int \nabla U_{j} \ol{R_{j}}dx|
  + |{\rm Im}\int \Lambda U_{j} \ol{R_{j}}dx|
  + |{\rm Im}\int  \varrho_{j} \ol{R_{j}}dx|\leq Ce^{-\frac{\delta}{T- t}}\|R\|_{L^2}.
\ea\ee
\end{lemma}

{\bf Proof.}
By the orthogonality condition (\ref{ortho-cond-Rn-wn}),
\begin{align}
{\rm Re}\int (x-\a_{j}) U_{j}(t)\ol{R_{j}}(t)dx=-\sum_{l\neq j}{\rm Re}\int (x-\a_{j}) U_{j}(t)\ol{R_{l}}(t)dx,
\end{align}
which along with Lemma \ref{Lem-inter-est} yields immediately that
for some $\delta>0$,
\begin{align}
 |{\rm Re}\int (x-\a_{j}) U_{j}(t)\ol{R_{j}}(t)dx|
\leq Ce^{-\frac{\delta}{T-t}}\|R\|_{L^2}.
\end{align}
The remaining four estimates in \eqref{Orth-almost}
can be proved similarly.
\hfill $\square$

We are now ready to prove Proposition \ref{Prop-Mod-bdd}.

{\bf Proof of Proposition \ref{Prop-Mod-bdd}.}
The proof is similar to that of \cite[Proposition 4.3]{SZ20}.

Below, we take the modulation equation $\lbb^2_j \dot{\g}_j + \g_j^2$,
corresponding to the direction $\Lambda U_j$,
for an example to illustrate the main arguments
and to show that the scalar
${\rm Re} \<R_j, U_j\>$ in the unstable direction $U_j$
is also required to bound the modulation equation.

Taking the inner product of \eqref{eq-U-R} with
$\Lambda  {U_{j}}$ and then taking the real part
we get
\begin{align}\label{eq-Ul-Rl}
   &-{\rm Im}\langle\partial_tR,\Lambda U_{j}\rangle
   +{\rm Re}\langle\Delta R_{j}+(1+\frac 2d)|U_{j}|^{\frac{4}{d}}R_{j}+\frac 2d|U_{j}|^{\frac{4}{d}-2}U_{j}^2\ol{{R}_{j}},\Lambda U_{j}\rangle  \nonumber \\
   &+{\rm Re}\langle i\partial_tU_{j}+\Delta U_{j}+|U_{j}|^{\frac{4}{d}}U_{j},\Lambda U_{j}\rangle  \nonumber \\
 = &-{\rm Re} \<\sum_{l\neq j} (\Delta R_{l}+(1+\frac 2d)|U_{l}|^{\frac{4}{d}}R_{l}+\frac 2d|U_{l}|^{\frac{4}{d}-2}U_{l}^2\ol{{R}_{l}}) + H_1,\Lambda U_{j} \>   \nonumber \\
   & -{\rm Re} \< \sum_{l\neq j}( i\partial_tU_{l}+\Delta U_{l}+|U_{l}|^{\frac{4}{d}}U_{l}) + H_2,\Lambda U_{j} \>  \nonumber \\
   &-{\rm Re}\langle f''(U,R)\cdot R^2,\Lambda U_{j}\rangle  \nonumber \\
   &- {\rm Re}\< \sum\limits_{l=1}^K (b\cdot \nabla (U_{l} + R_{l}) +c(U_{l} + R_{l})), \Lambda U_{j}\>,
\end{align}
where $H_1$ an $H_2$ are given by \eqref{G} and \eqref{H}, respectively.

As we shall see below that,
the right-hand side of equation \eqref{eq-Ul-Rl}
merely contribute negligibly small errors.

Actually,
we may take $t^*$ close to $T_*$
such that \eqref{aj-xj} and thus Lemma \ref{Lem-inter-est} hold.
By Lemma \ref{Lem-inter-est},
the interactions between different profiles are exponentially small,
and thus we infer that for some $\delta>0$,
\begin{align}  \label{RHS-Rj-Ul}
  | \<  \sum_{l\neq j} (\Delta R_{l}+(1+\frac 2d)|U_{l}|^{\frac{4}{d}}R_{l}+\frac 2d|U_{l}|^{\frac{4}{d}-2}U_{l}^2 \ol{{R}_{l}})
              + H_1,\Lambda U_{j} \>|
\leq C |\lbb|^{-2}e^{-\frac{\delta}{T-t}}\|R\|_{L^2},
\end{align}

Similarly,
by \eqref{equa-Ut},
\begin{align} \label{RHS-Uj-Ul}
    | \<\sum_{l\neq j}( i\partial_tU_{l}+\Delta U_{l}+|U_{l}|^{\frac{4}{d}}U_{l}) + H_2,\Lambda U_{j} \>|
\leq C  |\lbb|^{-2}e^{-\frac{\delta}{T-t}}(1+Mod)  ,
\end{align}

For the remainder $f''(U,R)\cdot R^2$ containing high order terms of $R$,
using \eqref{f''vR2-bdd} with $U$ replacing $v$,
\eqref{fLp-fH1} and $\|R\|_{H^1}\leq 1$,
we get
\begin{align} \label{RHS-Nf}
|\<f''(U,R)\cdot R^2, \Lambda U_{j} \>|
    \leq C |\lbb|^{-2}  ( \|R\|^2_{L^2}+   \|R\|_{H^1}^{3}).
\end{align}

Regarding the last term involving $b$ and $c$
on the right-hand side of  \eqref{eq-Ul-Rl},
using Lemma \ref{Lem-inter-est} and
\eqref{Uj-Qj-Q}, \eqref{Rj-ej} to rewrite it in the renormalized variables
we have
\begin{align} \label{bc-wtbwtc}
   &   {\rm Re}\<\sum\limits_{l=1}^K  (b\cdot \nabla (U_{l} + R_{l}) +c(U_{l} + R_{l})), \Lambda U_{j}\>   \nonumber \\
   = & {\rm Re} \< \lbb_{j}^{-1} \wt {b}\cdot\nabla (Q_{j}+\varepsilon_{j})
+ \wt {c}(Q_{j}+\varepsilon_{j}),\Lambda Q_{j} \>
       + \calo(e^{-\frac{\delta}{T-t}}(1+\|R\|_{L^2})),
\end{align}
where $\wt{b}(y):=b(\lbb_{j}y+\alpha_{j})$
and $\wt{c}(y):=c(\lbb_{j}y+\alpha_{j})$,
$y\in \bbr^d$.
Then,
using \eqref{b}, \eqref{c} and
integrating by parts formula
we get
\begin{align*}
  & {\rm Re} \<\lbb_j^{-1} \wt b \cdot \na (Q_j+\ve_j)
             +  \wt c (Q_j+\ve_j), \Lambda Q_j\>  \\
  =&2 {\rm Im} \sum\limits_{k=1}^N
     B_k \int \Delta \wt{\phi_k} (Q_j +\ve_j)  \ol{\Lambda Q_j} dy
     + 2 \lbb_j^{-1} {\rm Im}
      \sum\limits_{k=1}^N B_k  \int (Q_j + \ve_j)\na \wt{\phi_k} \cdot \na (\ol{\Lambda Q_j}) dy \\
  & -  \sum\limits_{l=1}^d {\rm Re}
      \int ( \sum\limits_{k=1}^N   \partial_l \wt{\phi_k} B_k)^2
           (Q_j+\ve_j) \ol{\Lambda Q_j}dy
    -{\rm Im}\sum\limits_{k=1}^N\int  \Delta \widetilde{\phi_k} B_k(Q_j+\ve_j)  \ol{\Lambda Q_j} dy,
\end{align*}
where $\partial^\nu \wt{\phi_k}(y) := (\partial^\nu \phi_k)(\lbb_j y + \a_j)$, $|\nu|\leq 2$.
Note that,
by the flatness condition (\ref{degeneracy})
and  the fact that  $\pa_y^\nu \phi_k \in L^\9$ for any multi-index $\nu$,
\begin{align}\label{phik-Taylor}
| \partial_y^\nu \wt{\phi_k} (y)|
\leq C(\lambda_{j} y+\alpha_{j}-x_j)^{\nu_*+1-|\nu|}
\leq C P^{\nu_*+1-|\nu|} \<y\>^{\nu_*+1}, \ \ 0\leq |\nu|\leq \nu_*.
\end{align}
This yields that
\begin{align} \label{RHS-b-c}
    |{\rm Re}\langle \lbb_{j}^{-1}\tilde{b}\cdot\nabla (Q_{j}+\varepsilon_{j})
    +\tilde{c}(Q_{j}+\varepsilon_{j}), \Lambda Q_{j}\rangle|
    \leq  C |\lbb|^{-2} P^{\nu_*+1}(1+\|R\|_{L^2}).
\end{align}

Thus,
we conclude from estimates \eqref{RHS-Rj-Ul}-\eqref{RHS-b-c} that
\begin{align} \label{Mod-RHS-bdd}
   {\rm R.H.S.\ of}\ \eqref{eq-Ul-Rl}
   \leq C |\lbb|^{-2}((e^{-\frac{\delta}{T-t}}+P^{\nu_*+1})(1+ \|R\|_{L^2})
              +e^{-\frac{\delta}{T-t}} Mod
               + \|R\|^2_{L^2}
               + \|R\|_{H^1}^{3} ).
\end{align}

Regarding the left-hand side,
by the orthogonality condition (\ref{ortho-cond-Rn-wn}), \eqref{equa-Ut} and Lemma \ref{Lem-inter-est},
\begin{align}  \label{ptR-Ul}
{\rm Im}\langle\partial_t R,\Lambda U_{j}\rangle
&={\rm Im}\langle \Lambda R,\partial_t U_{j}\rangle
={\rm Im}\langle \Lambda R_{j},\partial_t U_{j}\rangle+
   \sum_{l\neq j}{\rm Im}\langle \Lambda R_{l},\partial_t U_{j}\rangle \nonumber  \\
&={\rm Im}\langle \Lambda R_{j},\partial_t U_{j}\rangle
      + |\lbb|^{-2}\calo(Mod + e^{-\frac{\delta}{T-t}})\|R\|_{L^2}.
\end{align}
Then, using the identity \eqref{equa-Ut} and
the renormalized variables $Q_j, \ve_j$ in \eqref{Uj-Qj-Q} and \eqref{Rj-ej}, respectively,
we get
\begin{align}\label{ptR-x2Uj}
  \lbb_j^2{\rm Im}\langle \Lambda R_{j},\partial_t U_{j}\rangle
  =& -  {\rm Re}\langle \Lambda \varepsilon_{j},\Delta Q_{j} +|Q_j|^{\frac{4}{d}}Q_j\rangle
   + \calo(Mod\|R\|_{L^2}) \nonumber \\
=&  {\rm Re} \<\ve_j, \Lambda Q_j\>
    +  \g_j{\rm Im}\langle \Lambda \varepsilon_{j},\Lambda Q_{j}\rangle
   -2 \b_j{\rm Im}\langle \Lambda \varepsilon_{j},\nabla Q_{j}\rangle  \nonumber \\
   & + \calo((Mod+P^2)\|R\|_{L^2}),
\end{align}
where in the last step
we used the almost orthogonality \eqref{Orth-almost} and the identity
\begin{align} \label{DQj-Qj-fQj}
  \Delta Q_j - Q_j + |Q_j|^\frac 4d Q_j
  = |\beta_j - \frac{\g_j}{2}|^2 Q_j
    - i \g_j \Lambda Q_j
    + 2i \beta_j \cdot \na Q_j.
\end{align}

Furthermore,
using  \eqref{equa-Ut}, the identities
\begin{align}
   & \Lambda  Q_j = (\Lambda Q + i(\beta_j \cdot y - \frac 12 \g_j |y|^2) Q)e^{i(\beta_j\cdot y - \frac 14 \g_j |y|^2)}, \label{LaQj-LaQ} \\
   & \na  Q_j = (\na Q + i (\beta_j - \frac 12 \g_j y) Q)e^{i(\beta_j\cdot y - \frac 14 \g_j |y|^2)}, \label{naQj-naQ}
\end{align}
and $\<\Lambda Q, |y|^2 Q\> = - \|yQ\|_{L^2}^2$
we have
\begin{align} \label{ptUj-DUj-fUj}
     \lbb_j^2 {\rm Re}\langle i\partial_tU_{j}+\Delta U_{j}+|U_{j}|^{\frac{4}{d}}U_{j},\Lambda U_{j}\rangle
   = - \frac{1}{4}\|yQ\|_2^2(\lbb_j^2\dot{\g_j}+\g_j^2)
     + \calo(Mod  |\beta_j|).
\end{align}

Thus,
plugging \eqref{ptR-Ul}, \eqref{ptR-x2Uj} and \eqref{ptUj-DUj-fUj}  into \eqref{eq-Ul-Rl}
and rearranging the terms according to the orders of the renormalized viarable $\ve_j$
we get
\begin{align} \label{LHS-Modequa}
   & \lbb_j^2 \times ({\rm L.H.S.\ of}\ \eqref{eq-Ul-Rl})  \\
=&- \frac{1}{4}\|yQ\|_2^2(\lbb_j^2\dot{\g_j}+\g_j^2)
   + {\rm Re}\langle\Delta \varepsilon_{j}-\varepsilon_j+(1+\frac 2d)|Q_{j}|^{\frac{4}{d}}\varepsilon_{j}+\frac 2d|Q_{j}|^{\frac{4}{d}-2}Q_{j}^2\ol{{\varepsilon}_{j}},\Lambda Q_{j}\rangle  \nonumber \\
&- \g_j{\rm Im}\langle \Lambda \varepsilon_{j},\Lambda Q_{j}\rangle+2\b_j{\rm Im}\langle \Lambda \varepsilon_{j},\nabla Q_{j}\rangle
   + \calo((Mod+P^2 )\|R\|_{L^2} + Mod |\beta_j|).\nonumber
\end{align}

By \eqref{L+-L-}, straightforward computations show that,
if $\ve_j = \ve_{j,1} + i \ve_{j,2}$,
\begin{align}
& {\rm Re}\langle\Delta \varepsilon_{j}-\varepsilon_j+(1+\frac 2d)|Q_{j}|^{\frac{4}{d}}\varepsilon_{j}+\frac 2d|Q_{j}|^{\frac{4}{d}-2}Q_{j}^2\ol{{\varepsilon}_{j}},\Lambda Q_{j}\rangle   \nonumber \\
=& -\<L_+ \ve_{j,1}, \Lambda Q\>
  -\<L_+ \ve_{j,2}, (\beta_j\cdot y - \frac{\g_j}{4}|y|^2)\Lambda Q\>  \nonumber \\
 & -\<L_- \ve_{j,2}, (\beta_j\cdot y - \frac{\g_j}{2} |y|^2)Q\>
  + \calo(P^2 \|R\|_{L^2}).
\end{align}
Note that
\begin{align*}
   & L_+(\beta_j\cdot y - \frac{\g_j}{4}|y|^2) \Lambda Q
     = (\beta_j\cdot y - \frac{\g_j}{4}|y|^2) L_+ \Lambda Q
      + \g_j \Lambda^2 Q - 2 \beta_j \cdot \na \Lambda Q, \\
   & L_-(\beta_j\cdot y - \frac{\g_j}{2}|y|^2) Q
     = (\beta_j\cdot y - \frac{\g_j}{2}|y|^2)L_- Q
       + 2\g_j \Lambda Q - 2\beta_j \cdot \na Q.
\end{align*}
Taking into account the self-adjointness of $L_{\pm}$
and $L_+ \Lambda Q = -2 Q$, $L_-Q=0$
we get
\begin{align} \label{linear.1-Modequa}
& {\rm Re}\langle\Delta \varepsilon_{j}-\varepsilon_j+(1+\frac 2d)|Q_{j}|^{\frac{4}{d}}\varepsilon_{j}+\frac 2d|Q_{j}|^{\frac{4}{d}-2}Q_{j}^2\ol{{\varepsilon}_{j}},\Lambda Q_{j}\rangle   \nonumber \\
=& 2 \<\ve_{j,1}, Q\> + 2 \<\ve_{j,2}, (\beta_j\cdot y - \frac{\g_j}{4}|y|^2)Q\>  \nonumber \\
 &  + \g_j \<\Lambda \ve_{j,2}, \Lambda Q\>
   - 2 \beta_j\<\na \ve_{j,2}, \Lambda Q\>
   -2 \g_j \<\ve_{j,2}, \Lambda Q\>
   + 2 \beta_j \<\ve_{j,2}, \na Q\>.
\end{align}
Moreover,
we see that
\begin{align} \label{gj-Lvej-LQj}
    - \g_j {\rm Im} \<\Lambda \ve_j, \Lambda Q_j\>
    = -\g_j \<\Lambda \ve_{j,2}, \Lambda Q\> + \calo(P^2 \|R\|_{L^2}),
\end{align}
and by the almost orthogonality \eqref{Orth-almost},
\begin{align*}
   {\rm Im} \<\Lambda \ve_j, \na Q_j\>
    =  {\rm Im}  \<\na \ve_j, \Lambda  Q_j\>
      + {\rm Im} \<\ve_j, \na Q_j\>
   =  & {\rm Im}  \<\na \ve_j, \Lambda  Q_j\>
      + \calo(e^{-\frac{\delta}{T-t}}\|R\|_{L^2}),
\end{align*}
which yields that
\begin{align} \label{betaj-Lvej-naQj}
    2 \beta_j {\rm Im} \<\Lambda \ve_j, \na Q_j\>
   =& 2 \beta_j {\rm Im} \<\na \ve_j, \Lambda Q_j\> + \calo(e^{-\frac{\delta}{T-t}}\|R\|_{L^2}) \nonumber \\
   =& 2 \beta_j \<\na \ve_{j,2}, \Lambda Q\>
      + \calo((P^2 + e^{-\frac{\delta}{T-t}})\|R\|_{L^2}).
\end{align}

Thus,
we conclude from \eqref{linear.1-Modequa}, \eqref{gj-Lvej-LQj}, \eqref{betaj-Lvej-naQj}
and the almost orthogonality \eqref{Orth-almost}
that
\begin{align}   \label{linear-Modequa}
     &{\rm Re}\langle\Delta \varepsilon_{j}-\varepsilon_j+(1+\frac 2d)|Q_{j}|^{\frac{4}{d}}\varepsilon_{j}
         +\frac 2d|Q_{j}|^{\frac{4}{d}-2}Q_{j}^2\ol{{\varepsilon}_{j}},\Lambda Q_{j}\rangle \nonumber \\
    &  - \g_j{\rm Im}\langle \Lambda \varepsilon_{j},\Lambda Q_{j}\rangle+2\b_j{\rm Im}\langle \Lambda \varepsilon_{j},\nabla Q_{j}\rangle \nonumber \\
    =& 2 \<\ve_{j,1}, Q\> + 2 \<\ve_{j,2}, (\beta_j\cdot y - \frac{\g_j}{4}|y|^2)Q\>
        -2 \g_j \<\ve_{j,2}, \Lambda Q\>
   + 2 \beta_j \<\ve_{j,2}, \na Q\>  \nonumber \\
     & + \calo((P^2+e^{-\frac{\delta}{T-t}})\|R\|_{L^2})   \nonumber \\
     =&2{\rm Re}\<R_j, U_j\>
         + \calo((P^2+e^{-\frac{\delta}{T-t}})\|R\|_{L^2}).
\end{align}
This along with \eqref{LHS-Modequa} yields that
\begin{align} \label{Mod-LHS-bdd}
   \lbb_j^2 \times ({\rm L.H.S.\ of}\ \eqref{eq-Ul-Rl})
  =&- \frac{1}{4}\|yQ\|_2^2(\lbb_j^2\dot{\g_j}+\g_j^2)+ 2 {\rm Re}\langle R_j, U_j\rangle   \nonumber \\
  &+ \calo((Mod+P^2+e^{-\frac{\delta}{T-t}})\|R\|_{L^2} + Mod |\beta_j|).
\end{align}

Then, combining \eqref{Mod-RHS-bdd} and \eqref{Mod-LHS-bdd}
we obtain that for each $1\leq j\leq K$,
\begin{align}
|\lbb_j^2\dot{\g_j}+\g_j^2|
   \leq& C \big( Mod (P + \|R\|_{L^2}+e^{-\frac{\delta}{T-t}} )+|{\rm Re}\langle R_j, U_j\rangle|
          +(e^{-\frac{\delta}{T-t}}+P^{\nu_*+1})(1+ \|R\|_{L^2})  \nonumber \\
        &\quad + P^2 \|R\|_{L^2}+\|R\|_{L^2}^2
     + \|R\|_{H^1}^3 \big).
\end{align}

Similar arguments apply also to the remaining four modulation equations.
Actually, taking the inner products of equation (\ref{eq-U-R})
with $i(x-\a_j) U_j$, $i|x-\a_j|^2 U_j$, $\nabla {U_j}$, $ \varrho_{j}$,
respectively,
then taking the real parts
and using analogous arguments as above,
we can obtain the same bounds for
$|{\lbb_j \dot{\a}_j}-2\beta_j|$,
$|\lbb_j \dot{\lbb}_j + \g_j|$,
$|\lbb^2_j\dot{\beta}_j + \beta_j \g_j|$
and $|\lbb_j^2\dot{\theta}_j - 1-|\beta_j|^2|$,
respectively.
We then get
\begin{align}\label{Mod-bdd-l}
Mod_{j}(t)
\leq& C\big( Mod(P+\|R\|_{L^2}+e^{-\frac{\delta}{T-t}})
           +\sum_{j=1}^{K} |{\rm Re}\langle R_j, U_j\rangle|
           + (e^{-\frac{\delta}{T-t}}+P^{\nu_*+1}) (1+ \|R\|_{L^2})  \nonumber \\
    & \qquad  + P^2 \|R\|_{L^2}
           +\|R\|_{L^2}^2+  \|R\|_{H^1}^{3} \big).
\end{align}

Therefore,
taking $T$ possibly even smaller
such that
$$(1+C) (P+ \sup_{t^*\leq t\leq T_*}\|R(t)\|_{H^1} +e^{- \frac{\delta}{T}}) \leq \frac 12$$
and then summing over $j$
we obtain \eqref{Mod-bdd} and finish the proof.
\hfill $\square$

\section{Bootstrap estimates} \label{Sec-Bootstrap}

This section is mainly devoted to the bootstrap type estimates
of the remainder $R$ and the modulation parameters $\calp$,
which are the key ingredients in the construction of
multi-bubble blow-up solutions
in Section \ref{Sec-Exit} later.
The main result is formulated in Theorem \ref{Thm-u-Boot} below.

\begin{theorem}[Bootstrap] \label{Thm-u-Boot}
Let $\ve^*>0$ be sufficiently small,
$0<\zeta< \frac{1}{12}$.
For any $\ve\in (0,\ve^*]$,
let $T=T(M)$ be small enough,
satisfying \eqref{T-M-0},
and fix $T_*\in (0,T)$ .
Suppose that there exists $t^*\in(0,T_*)$ such that
$u$ admits the unique geometrical decomposition \eqref{u-dec} on $[t^*,T_*]$
and  the following estimates hold for  $\kappa:= \nu_* -3(\geq 2)$:

$(i)$ For the reminder term,
\begin{align} \label{R-Tt}
\|\nabla R(t)\|_{L^2}\leq (T-t)^{\kappa},\quad\|R(t)\|_{L^2}\leq (T-t)^{\kappa+1}.
\end{align}
$(ii)$ For the modulation parameters,  $1\leq j\leq K$,
\begin{align}
&|\la_{j}(t) - \omega_j (T-t) | + |\gamma_{j}(t)  - \omega_j^2 (T-t) |\leq (T-t)^{\kappa+1+\zeta},   \label{lbbn-Tt} \\
&|\al_{j}(t)-x_j|+|\beta_{j}(t)|\leq (T-t)^{\frac{\kappa}{2}+1+\zeta},  \label{anbn-Tt}\\
&|\theta_{j}(t) - (\omega_j^{-2}(T-t)^{-1} + \vartheta_j)| \leq (T-t)^{\kappa-1 +\zeta}. \label{thetan-Tt}
\end{align}

Then,
there exists $t_*\in [0,t^*)$ such that
\eqref{u-dec} holds on the larger interval $[t_*, T_*]$
and
the coefficients in estimates \eqref{R-Tt}-\eqref{thetan-Tt} can be refined to $1/2$,
i.e., for any $t\in [t_*, T_*]$, $1\leq j\leq K$,
\begin{align}
&\|\nabla R(t)\|_{L^2}\leq \frac{1}{2}(T-t)^{\kappa},\quad\|R(t)\|_{L^2}\leq \frac{1}{2}(T-t)^{\kappa+1},  \label{wn-Tt-boot-2} \\
&\lf|\la_{j}(t) -\omega_j (T-t) \rt| + \ \lf|\gamma_{j}(t) -\omega_j^2(T-t) \rt|\leq \frac{1}{2}(T-t)^{\kappa+1+ \zeta}, \label{lbbn-Tt-boot-2} \\
&|\al_{j}(t)-x_j|+|\beta_{j}(t)|\leq \frac{1}{2} (T-t)^{\frac{\kappa}{2}+1+ \zeta},  \label{anbn-Tt-boot-2} \\
&  |\theta_{j}(t) - (\omega_j^{-2}(T-t)^{-1} + \vartheta_j) |\leq \frac{1}{2} (T-t)^{\kappa-1 +\zeta}. \label{thetan-Tt-boot-2}
\end{align}
\end{theorem}

In order to prove Theorem \ref{Thm-u-Boot},
we may take $t_* \in [0,t^*)$, sufficiently close to $t^*$,
such that
$u$ still has the geometrical decomposition \eqref{u-dec}
on the larger interval  $[t_*, T_*]$
(this is possible because the Jacobian of transformation
is continuous in time).
Moreover,
by virtue of the local well-posedness theory and the $C^1$-regularity of modulation parameters,
taking $t_*$ possibly closer to $t^*$,
we have that for any $t\in [t_*,t_n]$,
\begin{align}
&\|\nabla R(t)\|_{L^2}\leq 2(T-t)^\kappa, \ \ \|R(t)\|_{L^2}\leq 2(T-t)^{\kappa+1}, \label{R-Tt2} \\
&\lf|\la_j(t)- \omega_j(T-t)\rt| + \lf|\gamma_j(t) - \omega_j^2(T-t)\rt|\leq 2(T-t)^{\kappa+1+\zeta}, \label{lbb-ga-Tt2}   \\
&|\al_j(t) - x_j| + |\beta_j(t)|\leq 2(T-t)^{\frac \kappa 2+1+\zeta}, \label{ab-Tt2}  \\
&   |\theta_j(t)- (\omega_j^{-2}(T-t)^{-1} + \vartheta_j) |\leq 2(T-t)^{\kappa-1+\zeta}. \label{theta-Tt2}
\end{align}
By \eqref{lbb-ga-Tt2} and \eqref{theta-Tt2},
we may also take $T$ sufficiently small such that \eqref{aj-xj} holds,
and thus Lemma \ref{Lem-inter-est} is applicable below.

\begin{remark} \label{Rem-lbbj-gj-P-Tt}
We infer from \eqref{lbb-ga-Tt2} that for $T$ small enough,
$\lbb_j, \g_j, P$ are comparable with $T-t$, i.e.,
\begin{align} \label{lbb-g-t}
     C_1 (T-t)\leq \lbb_j, \g_j, P \leq C_2 (T-t).
\end{align}
where
$C_1, C_2$ are positive constants
independent of $\ve$.
\end{remark}

By virtue of Proposition \ref{Prop-Mod-bdd}
and Proposition \ref{Prop-mass-local} below,
we have the refined estimate for the modulation parameters below.
\begin{lemma} \label{Lem-Mod-w-lbb}
Assume  estimates \eqref{R-Tt}-\eqref{thetan-Tt} to hold
with $T$ sufficiently small.
Then, there exists $C>0$ such that
\begin{align} \label{Mod-w-lbb}
Mod(t) \leq C (T-t)^{\kappa+3}, \ \ \forall t\in[t_*,T_*].
\end{align}
\end{lemma}
Moreover,
by equation \eqref{equa-u-RNLS},
the remainder $R$ satisfies the equation
\begin{align} \label{equa-R}
   i\partial_t R +\Delta R+(f(u)-f(U))+b \cdot \nabla R+c  R=-\eta,
\end{align}
where
\begin{align} \label{etan-Rn}
    \eta = i\partial_t U +\Delta U+f(U)+b \cdot \nabla U+c U.
\end{align}
The estimates of
$U_j$, $R$ and $\eta$ are contained in Lemmas \ref{Lem-P-ve-U} and \ref{Lem-eeta} below.
\begin{lemma}   \label{Lem-P-ve-U}
Assume   estimates \eqref{R-Tt}-\eqref{thetan-Tt}  to hold
and let $\ve_j$ be defined in \eqref{Rj-ej}.
Then, there exists $C>0$
such that  for all $t\in [t_*, T_*]$, $1\leq j\leq K$,
\begin{align*}
    &\|\varepsilon_j(t)\|_{L^2}=\|R_j(t)\|_{L^2} \leq C (T-t)^{\kappa+1}, \ \
     {\lambda_j^{-1}}\|\na \ve_j(t)\|_{L^2}=\|\na R_j(t)\|_{L^2} \leq C (T-t)^{\kappa},  \\
   &\|U_j(t)\|_{L^2}= \|Q\|_{L^2}, \ \ \|\nabla U_j(t)\|_{L^2} + \|\frac{x_j-\a_j(t)}{\lbb_j(t)} \cdot \na U_j(t)\|_{L^2} \leq  C (T-t)^{-1}.
\end{align*}
\end{lemma}
\begin{lemma}  \label{Lem-eeta}
Assume estimates \eqref{R-Tt}-\eqref{thetan-Tt}  to hold with $T$ sufficiently small.
Then, there exists a constant $C>0$
such that for  $t\in[t_*,T_*]$
and multi-index $\nu$ with $|\nu|\leq 2$,
\begin{align}
\|\partial^\nu \eta(t)\|_{L^2}
      \leq C(T-t)^{\kappa+1-|\nu|}.    \label{eta-L2}
\end{align}
\end{lemma}
The proof is postponed to the Appendix for simplicity.

The remainder of Section \ref{Sec-Bootstrap}
is devoted to the proof of Theorem \ref{Thm-u-Boot}.
We first derive the estimates of the localized mass and energy
in Subsection \ref{Subsec-Energy},
and then in Subsection \ref{Subsec-Mono}
we derive the key monotonicity property of the generalized energy,
involving a Morawetz type term and localized funcitons,
which actually constitutes the most technical part of this section.
The detailed proof of Theorem \ref{Thm-u-Boot} is then given in Subsection \ref{Subsec-Boot-proof}.
We shall assume estimates \eqref{R-Tt}-\eqref{thetan-Tt}  to hold on $[t_*,T_*] \subseteq [0,T)$
with $T$ small enough and satisfying \eqref{T-M-0}
throughout Subsections \ref{Subsec-Energy}-\ref{Subsec-Boot-proof}.

\subsection{Estimates of localized mass and energy} \label{Subsec-Energy}
\begin{proposition}  [Estimate of localized mass]    \label{Prop-mass-local}
There exists $C>0$
such that for any $t\in [t_*, T_*]$ and $1\leq j\leq K$,
\begin{align} \label{mass-local-esti}
    2{\rm Re}\int \overline{U_{j}}R_{j}dx+\int |R(t)|^2\Phi_jdx
  = \calo((T-t)^{2 \kappa +2}),
\end{align}
where $R_j := R\Phi_j$
with $\Phi_j$ the local functions defined in \eqref{phi-local}.
\end{proposition}

\begin{remark}
$(i)$. The estimate \eqref{mass-local-esti} allows us to control
the scalar product along the direction $Q$
when applying the localized coercivity in Corollary \ref{Cor-coer-f-local}.

$(ii)$. It should be mentioned that,
the proof of \eqref{mass-local-esti}
relies on the analysis of the localized mass $\int |u|^2 \Phi_j dx$,
instead of the usual whole mass $\|u\|_{L^2}^2$.
This is quite different from the single bubble case in \cite{R-S,SZ20}.
\end{remark}

{\bf Proof of Proposition \ref{Prop-mass-local}.}
Using \eqref{u-dec} and Lemma \ref{Lem-inter-est} we have that for some $\delta>0$,
\begin{align*}
\int |u |^2\Phi_jdx
=& \int |U|^2\Phi_jdx +\int |R |^2\Phi_jdx +2{\rm Re}\int \overline{U} R\Phi_j dx  \nonumber \\
=& \int |U |^2\Phi_jdx +\int |R |^2\Phi_jdx +2{\rm Re}\int \overline{U_j} R_j dx
   + \calo(e^{-\frac{\delta}{T-t}}\|R\|_{L^2}),
\end{align*}
which yields that
\begin{align} \label{uL2-phi}
|2{\rm Re}\int (\overline{U_{j}}R_{j})(t) dx+\int |R(t)|^2\Phi_j dx|
& \leq|\int |u(t)|^2\Phi_j dx-\int |u(T_*)|^2\Phi_j dx| \nonumber \\
&+|\int |u(T_*)|^2\Phi_j dx-\int |U(t)|^2\Phi_j dx|+C e^{-\frac{\delta}{T-t}} \|R\|_{L^2}.
\end{align}

For the first term on the right-hand side of \eqref{uL2-phi},
we use equation \eqref{equa-u-RNLS} to get
\begin{align} \label{du2-bc}
 \frac{d}{dt}\int |u|^2\Phi_jdx
 =& {\rm Im}\int (2\ol{u}\nabla u+b|u|^2)\cdot  \nabla\Phi_j dx  \nonumber  \\
\leq& \int_{ |x-x_l|\geq 4\sigma,1\leq l\leq K}2|\ol{u}\nabla u|+|b||u|^2 dx,
\end{align}
where $\sigma$ is given by \eqref{sep-xj-0}.
By \eqref{ab-Tt2},
we may take $t_*$ close to $T_*$ such that $|x_l(t) - \a_l(t)|\leq \sigma$
for any $t\in [t_*, T_*]$, $1\leq l\leq K$.
This along with \eqref{Q-decay} and \eqref{u-dec} yields that
\begin{align*}
     |\frac{d}{dt}\int |u|^2\Phi_jdx|
\leq& C\int_{|x-\alpha_{l}|\geq 3\sigma,1\leq l\leq K}
      |U+R||\na U + \na R| + |U+R|^2 dx  \\
\leq& C(\|R\|_{L^2}^2 + \|R\|_{L^2}\|\na R\|_{L^2}
         + \sum\limits_{l=1}^K (\int_{|y|\geq \frac{3\sigma}{\lbb_l}} |Q|^2 dy)^\frac 12 \|\na R\|_{L^2} \nonumber \\
    &\qquad      + \sum\limits_{l=1}^K
             \int_{|y|\geq \frac{3\sigma}{\lbb_l}}
             |Q|^2 + \lbb_l^{-2}|\na Q_l|^2 dy) \\
\leq& C(\|R\|_{L^2}^2 + \|R\|_{L^2}\|\na R\|_{L^2}  + e^{-\frac{\delta}{T-t}}),
\end{align*}
where $\delta >0$.
Hence, we obtain that for some $\delta>0$,
\begin{align}  \label{u-utn}
   &|\int |u(t)|^2\Phi_j dx-\int |u(T_*)|^2\Phi_j dx| \nonumber \\
\leq& C\int_{t}^{T_*} \|R(s)\|_{L^2}^2 + \|R(s)\|_{L^2}\|\na R(s)\|_{L^2}  ds
+C e^{-\frac{\delta}{T-t}}.
\end{align}

Regarding the second term on the right-hand side of \eqref{uL2-phi},
we apply Lemma \ref{Lem-inter-est} to
extract the main blow-up profile $U_j$
\begin{align*}
\int |U(t)|^2\Phi_jdx
=\int |U_{j}(t)|^2 \Phi_j dx+ \calo(e^{-\frac{\delta}{T-t}}).
\end{align*}
Since
\begin{align*}
   \int|U_j(t)|^2 \Phi_j dx
   = \int |Q|^2 dy
      + \int |Q(y)|^2 (\Phi_j(\lbb_j(t) y +\a_j(t))-1) dy,
\end{align*}
and for some $\delta >0$,
\begin{align*}
    \int |Q(y)|^2 (1- \Phi_j(\lbb_j(t) y +\a_j(t))) dy
    \leq \int_{|y|\geq \frac{3\sigma}{\lbb_j(t)}} Q^2(y) dy
   \leq C e^{-\frac{\delta}{T-t}},
\end{align*}
we infer that
\begin{align}\label{U-Q-esti}
   \int |U (t)|^2\Phi_jdx
   =\|Q\|_{L^2}^2 + \calo(e^{-\frac{\delta}{T-t}}).
\end{align}
Similarly, we have
\begin{align}  \label{utn-Q-esti}
\int |u(T_*)|^2\Phi_j dx
 =\int |\sum_{l=1}^KS_l(T_*)|^2\Phi_j dx
 =\|Q\|_2^2+ \calo(e^{-\frac{\delta}{T-t}}).
\end{align}
We infer from \eqref{U-Q-esti} and \eqref{utn-Q-esti} that
\begin{align} \label{uT-Ut}
   | \int |u(T_*)|^2\Phi_j dx-\int |U(t)|^2\Phi_j dx | \leq C e^{-\frac{\delta}{T-t}}.
\end{align}

Therefore,
plugging \eqref{u-utn} and \eqref{uT-Ut}
into \eqref{uL2-phi}
we obtain
\begin{align*}
  |2{\rm Re}\int (\overline{U_{j}}R_{j})(t) dx+\int |R(t)|^2\Phi_jdx|
  \leq C (\int_t^{T_*}  \|R\|_{L^2}^2 + \|R\|_{L^2}\|\na R\|_{L^2} ds + e^{-\frac{\delta}{T-t}}),
\end{align*}
which along with \eqref{R-Tt2} yields \eqref{mass-local-esti}
for $T$ small enough
and finishes the proof.
\hfill $\square$

Theorem \ref{Thm-energy} below contains the estimate
of the variation of energy.
Unlike the deterministic case,
the energy \eqref{energy} is no longer conserved
and the corresponding variation
plays an important role
in the derivation of the refined estimate of the modulation parameter $\beta$  later
(see Lemma \ref{Lem-betaj} below).

\begin{proposition} [Variation of the energy] \label{Thm-energy}
There exists $C>0$
such that  for any $t\in[t_*,T_*]$,
\begin{align} \label{esti-Eut-Eutn}
      |E(u(t))-E(u(T_*))| \leq C(T-t)^{\kappa+1}.
\end{align}
\end{proposition}

{\bf Proof.}
The proof is quite similar to that of \cite[Theorem 5.6]{SZ20},
based on the Gagliardo-Nirenberg inequality \eqref{G-N}
and the estimate \eqref{phik-Taylor} of the spatial functions of noises
under Assumption $(A1)$.
Actually,  as in \cite[(5.20)]{SZ20}, we have
\begin{align} \label{dtE}
\frac{d}{dt}E(u_n)
=&-2\sum\limits_{k=1}^N B_k {\rm Re}\int \nabla^2 \phi_k(\nabla u_n,\nabla \ol{u_n})dx
+\frac{1}{2}\sum\limits_{k=1}^N B_k \int \Delta^2 \phi_k|u_n|^2dx  \nonumber \\
& +\frac{2}{d+2}\sum\limits_{k=1}^N  B_k \int\Delta \phi_k|u_n|^{2+\frac{4}{d}}dx
 -\sum\limits_{j=1}^d {\rm Im}\int \na(\sum\limits_{k=1}^N \partial_j \phi_k B_k)^2 \cdot \na u_n \ol{u_n} dx.
\end{align}
This yields that
\begin{align*}
 |\frac{d}{dt} E(u_n)|
 \leq &  C \|R\|^2_{H^1}
        + C \sum\limits_{k=1}^K \sum\limits_{j=1}^d
            \big(\int (|\na^2 \phi_k| + |\Delta \phi_k| + |\partial_j\phi_k \na \partial_j \phi_k|)
               (|\na U|^2 + |U|^{2+\frac 4d}) dx  \nonumber \\
   & \qquad  \qquad   \qquad   \qquad  \ + \int (|\Delta^2 \phi_k| + |\partial_j \phi_k \na \partial_j \phi_k|)
               |U|^2dx\big),
\end{align*}
which, via the change of variables, can be further bounded by,
up to some constant,
\begin{align*}
  & \|R\|^2_{H^1}
       +  \sum\limits_{k,l=1}^K
          \int  (T-t)^{-2} \sum\limits_{|\nu|\leq 2} |\partial^\nu \phi_k|(\lbb_ly+\a_l)
               (|\na Q_l|^2 + |Q_l|^{2+\frac 4d})
               + \sum\limits_{|\nu|\leq 4} |\partial^\nu \phi_k| (\lbb_ly+\a_l)
               |Q|^2dy.
\end{align*}
Thus,
using \eqref{phik-Taylor} and \eqref{R-Tt2}
we obtain
\begin{align*}
     |\frac{d}{dt} E(u_n)|
     \leq C (T-t)^\kappa,
\end{align*}
which immediately yields \eqref{esti-Eut-Eutn},
thereby finishing the proof.
\hfill $\square$

\subsection{Monotonicity of generalized energy}  \label{Subsec-Mono}
This subsection is mainly devoted to the
monotonicity property of a new generalized energy,
which is the key ingredient in the proof of the bootstrap estimate \eqref{wn-Tt-boot-2} of the remainder.

It should be mentioned that,
unlike the single blow-up point case in \cite{R-S,SZ20},
the new generalized energy \eqref{def-I} below
includes also the localized functions in an appropriate way,
such that the different profiles can be decoupled completely
and the key monotonicity property is still preserved.

More precisely,
let $\chi(x)=\psi(|x|)$ be a smooth radial function on $\R^d$,
where $\psi$ satisfies
$\psi'(r) = r$ if $r\leq 1$,
$\psi'(r) = 2- e^{-r}$ if $r\geq2$,
and
\be\label{chi}
 |\frac{\psi^{'''}(r)}{\psi^{''}(r)} |\leq C,
\ \ \frac{\psi'(r)}{r}-\psi^{''}(r) \geq0.
\ee
Let $\chi_A(x) :=A^2\chi(\frac{x}{A})$, $A>0$,
$f(u):= |u|^{\frac 4d} u$,
and $F(u):= \frac{d}{2d+4} |u|^{2+\frac 4d}$.
We shall also use the notations $f'(U,R)\cdot R$ and $f''(U,R)\cdot R^2$
as in \eqref{f'R} and \eqref{f''R2}, respectively.

We define the generalized energy by
\begin{align} \label{def-I}
I(t) := &\frac{1}{2}\int |\nabla R|^2dx+\frac{1}{2}\sum_{j=1}^K\int\frac{1}{\lambda_{j}^2} |R|^2 \Phi_jdx
           -{\rm Re}\int F(u)-F(U)-f(U)\ol{R}dx \nonumber \\
&+\sum_{j=1}^K\frac{\gamma_{j}}{2\lambda_{j}}{\rm Im} \int (\nabla\chi_A) (\frac{x-\alpha_{j}}{\lambda_{j}})\cdot\nabla R \ol{R}\Phi_jdx.
\end{align}

The key monotonicity property of the generalized energy
is formulated in Theorem \ref{Thm-I-mono} below.

\begin{theorem}  \label{Thm-I-mono}
There exist $C_1, C_2(A), C_3>0$
such that for  any $t\in[t_*,T_*]$
\begin{align} \label{dIt-mono-case2}
\frac{d I}{dt}
\geq C_1
     \sum_{j=1}^{K}\frac{1}{\lambda_{j}}
     \int (|\nabla R_{j}|^2+\frac{1}{\lbb_{j}^2} |R_{j}|^2 )
     e^{-\frac{|x-\alpha_{j}|}{A\lbb_j}}dx
    -C_2(A)(T-t)^{2 \kappa}
    - C_3\ve^* (T-t)^{2\kappa-1}.
\end{align}

\end{theorem}

\begin{remark}
Theorem \ref{Thm-I-mono} yields that
the derivative of the generalized energy
is almost positive, up to some error terms,
and thus the generalized energy is almost monotone.
We also mention that,
the error term of order $(T-t)^{2\kappa-1}$
corresponds to the frequencies $\{\omega_j\}_{j=1}^K$,
and the small coefficient $\ve^*$
is important later in the derivation of
the bootstrap estimate \eqref{wn-Tt-boot-2}
of the remainder $R$,
and also in the iteration arguments in the proof of uniqueness.
\end{remark}

In order to prove Theorem \ref{Thm-I-mono},
we separate $I$ into two parts $I= I^{(1)}+ I^{(2)}$,
where
\begin{align}
I^{(1)}&:=\frac{1}{2} \int |\nabla R|^2dx+\frac{1}{2}\sum_{j=1}^K\int\frac{1}{\lambda_{j}^2} |R|^2 \Phi_j dx
           -{\rm Re}\int F(u)-F(U)-f(U)\ol{R}dx, \label{I1}\\
I^{(2)}&:=\sum_{j=1}^K\frac{\gamma_{j}}{2\lambda_{j}}{\rm Im}
            \int (\nabla\chi_A) (\frac{x-\alpha_{j}}{\lambda_j})\cdot\nabla R\ol{R}\Phi_jdx.   \label{I2}
\end{align}
Below we treat $I^{(1)}$ and  $I^{(2)}$ separately in Lemmas \ref{Lem-I1t} and \ref{Lem-I2t}.
Let us first show the estimate of $I^{(1)}$.

\begin{lemma}   \label{Lem-I1t}
Consider the situations as in Theorem \ref{Thm-I-mono}.
Then, for every $t\in[t_*,T_*]$ we have
that for some $C_1, C_2>0$,
\begin{align} \label{I1t-case1}
 \frac{d I^{(1)}}{dt}
\geq &\sum_{j=1}^{K}\frac{\gamma_j}{\lambda_j^4} \|R_j\|_{L^2}^2
     -\sum_{j=1}^{K}\frac{\gamma_j}{\lambda_j^2} {\rm Re}\int (1+\frac 2d)|U_j|^\frac 4d|R_j|^2
      +\frac 2 d   |U_j|^{\frac 4d -2} \ol{U_j}^2R_j^2  dx   \nonumber  \\
&-\sum_{j=1}^{K}\frac{\gamma_j}{\lambda_j} {\rm Re}
      \int(\frac{x-\alpha_j}{\lambda_j})\cdot\nabla \ol{U_j}
      \bigg\{ \frac 2d (1+\frac 2d) |U_j|^{\frac {4}{d}-2}U_j|R_j|^2    \nonumber \\
 & \qquad \qquad  + \frac 1d (1+ \frac 2d)|U_j|^{\frac {4}{d}-2}\ol{U_j}R_j^2
   + \frac 1d (\frac 2d -1) |U_j|^{\frac 4d -4} U_j^3 \overline{R_j}^2 \bigg\}dx \nonumber \\
 &  - C_1 (T-t)^{2 \kappa}
   - C_2 \ve^* (T-t)^{2\kappa-1}.
\end{align}
\end{lemma}

{\bf Proof.}
Using the identities
\begin{align*}
   \partial_t F(u)
= {\rm Re} \lf(f(u) \partial_t \ol{u}\rt),
\ \ \partial_t f(U)
    = \partial_z f(U) \partial_t U
      + \partial_{\ol z} f(U) \partial_t \ol U,
\end{align*}
and the expansion \eqref{g-gzz-expan} we have
\begin{align*}
   \frac{d I^{(1)}}{dt}
   =& {\rm Im} \<\Delta R + f(u) - f(U), i\partial_t R\>
      - \lbb_j^{-2} {\rm Im} \<R_j, i\partial_t R\>  \\
    &  - \dot \lbb_j \lbb_j^{-3} \int |R|^2 \Phi_j dx
     - {\rm Re} \<f''(U, R)\cdot R^2, \partial_t U\>.
\end{align*}
Then, in view of \eqref{equa-R}, we obtain
\begin{align} \label{equa-I1t}
 \frac{d I^{(1)}}{dt}
=& -\sum_{j=1}^K \dot{\lbb}_{j}\lbb_{j}^{-3} {\rm Im}\int|R|^2\Phi_jdx
  -\sum_{j=1}^K \lbb_{j}^{-2} {\rm Im} \<f^\prime(U)\cdot R, R_j\>  \nonumber\\
  &  -{\rm Re} \< f''(U,R)\cdot R^2,  \partial_t {U} \>
    - \sum_{j=1}^K \lbb_{j}^{-2}{\rm Im} \< R \nabla\Phi_j, \na R\> \nonumber\\
  &-\sum_{j=1}^K \lbb_{j}^{-2} {\rm Im} \< f''(U,R)\cdot R^2, R_j\>
-{\rm Im} \<\Delta R - \sum_{j=1}^K \lbb_{j}^{-2} R_j +f(u)-f(U), \eta\>  \nonumber\\
& - {\rm Im} \<\Delta R -\sum\limits_{j=1}^K \lbb_j^{-2} R_j +f(u)-f(U), b\cdot \nabla R+cR  \>
=:  \sum\limits_{j=1}^7 I^{(1)}_{t,j},
\end{align}
where $\eta$ is given by \eqref{etan-Rn}.

As we shall see below,
the main orders of $\frac{dI^{(1)}}{dt}$
are contributed by the first three terms $I^{(1)}_{t,1}$, $I^{(1)}_{t,2}$ and $I^{(1)}_{t,3}$,
the
fourth term will contributes the error of order $(T-t)^{2\kappa-1}$
for which we shall treat
{\rm Case (I)} and {\rm Case (II)} separately,
while the remaining three terms are of the
negligible order $(T-t)^{2 \kappa}$.

{\it $(i)$ Estimate of $I^{(1)}_{t,1}$.}
Since by \eqref{Mod-w-lbb},
$|\frac{\lbb_j\dot \lbb_j+\g_j}{\lbb_j^4}| \leq C \frac{Mod}{\lbb_j^4} \leq C (T-t)^{ \kappa-1}$,
it follows from \eqref{R-Tt} that
\begin{align} \label{I1t1-esti}
   I^{(1)}_{t,1}
  =& \sum_{j=1}^{K} (\frac{\g_j}{\lbb_j^4}\int |R|^2\Phi_jdx-\frac{\lbb_j \dot{\lbb}_j +\g_j}{\lbb_j^4}\int |R|^2\Phi_jdx) \nonumber \\
  =&  \sum_{j=1}^{K} \frac{\g_j}{\lbb_j^4}\int |R|^2\Phi_jdx+\calo((T-t)^{2 \kappa})  \nonumber \\
  \geq& \sum_{j=1}^{K} \frac{\g_j}{\lbb_j^4} \|R_j\|_{L^2}^2
        - C (T-t)^{2\kappa},
\end{align}
where we also used the inequality $\Phi_j \geq \Phi_j^2$ in the last step.

{\it $(ii)$ Estimates of $I^{(1)}_{t,2}$ and  $I^{(1)}_{t,3}$.}
We apply Lemma \ref{Lem-inter-est}
to decouple different blow-up profiles to obtain
\begin{align*}
   I^{(1)}_{t,2} +  I^{(1)}_{t,3}
   =  - \sum\limits_{j=1}^K
       \frac{1}{\lbb_j^2} {\rm Im} \int \frac 2d |U_j|^{\frac 4d-2} U_j^2 \ol{R_j}^2 dx
     - \sum\limits_{j=1}^K  {\rm Re}
        \int f'' (U_j, R_j)\cdot R_j^2 \partial_t \ol{U_j} dx
     + \calo(e^{-\frac{\delta}{T-t}}).
\end{align*}
Then,
for each $1\leq j\leq K$,
straightforward computations show that
(see also the proof of \cite[(5.46),(5.49)]{SZ20})
\begin{align} \label{I1t2t3-esti}
   I^{(1)}_{t,2} +  I^{(1)}_{t,3}
   =&- \sum\limits_{j=1}^K \frac{\gamma_j}{\lambda_j^2} {\rm Re}
      \int (1+\frac 2d)|U_j|^{\frac {4}{d}}|R_j|^2
     + \frac 2 d |U_j|^{\frac 4d-2} \ol{U_j}^2 R_j^2 \nonumber \\
 &- \sum\limits_{j=1}^K\frac{\gamma_j}{\lambda_j} {\rm Re}
    \int (\frac{x-\alpha_j}{\lambda_j})\cdot\nabla \ol{U_j}
     \bigg\{\frac 2d(1+\frac 2d) |U_j|^{\frac {4}{d}-2}U_j|R_j|^2 \nonumber\\
       &\ \ \  + \frac 1 d (1+\frac 2d) |U_j|^{\frac {4}{d}-2}\ol{U_j}R_j^2
        + \frac 1d (\frac 2d-1) |U_j|^{\frac 4d-4} U_j^3 \overline{R_j}^2\bigg\}dx
      +\calo ((T-t)^{2\kappa}).
\end{align}

{\it $(iii)$ Estimate of $I^{(1)}_{t,4}$.}
We consider {\rm Case (I)} and {\rm Case (II)} separately.
First,
in {\rm Case (I)},
since $\sum_{j=1}^{K} \na \Phi_j(x)=0$,
we see that
\begin{align*}
  |I^{(1)}_{t,4}|
  =& \sum\limits_{j=1}^K
     (\frac{1}{\lbb_j^2} - \frac{1}{\omega^2(T-t)^2})
     {\rm Im} \<R\na \Phi_j, \na R\> \\
  \leq&  \sum\limits_{j=1}^K
        \frac{|\lbb_j-\omega(T-t)||\lbb_j+\omega(T-t)|}{\lbb_j^2 \omega^2 (T-t)^2}
        \|R\na \Phi_j\|_{L^2} \|\na R\|_{L^2}\\
        \leq & \sum\limits_{j=1}^K
        \frac{|\lbb_j-\omega(T-t)||\lbb_j+\omega(T-t)|}{\lbb_j^2 \omega^2 (T-t)} (\frac{\|R\|^2_{L^2}}{(T-t)^2}+\|\nabla R\|_{L^2}^2) .
\end{align*}
Since
$|\omega-\omega_j|\leq \ve^*$ for any $1\leq j\leq K$,
using \eqref{lbb-ga-Tt2} and \eqref{lbb-g-t} we see that
\begin{align*}
        & \frac{|\lbb_j-\omega(T-t)||\lbb_j+\omega(T-t)|}{\lbb_j^2 \omega^2 (T-t)}  \\
    \leq&  \frac{(|\lbb_j-\omega_j(T-t)|+|(\omega_j-\omega)(T-t)|)|\lbb_j+\omega(T-t)|}{\lbb_j^2 \omega^2 (T-t)}  \\
    \leq& C ((T-t)^{\kappa-1+\zeta} + \ve^*(T-t)^{-1}).
\end{align*}
This along with \eqref{R-Tt2} yields that
\begin{align} \label{I1t4-esti-case1}
     |I^{(1)}_{t,4}|
     \leq C ((T-t)^{3\kappa-1+\zeta} + \ve^* (T-t)^{2\kappa-1})
     \leq C ((T-t)^{2\kappa} + \ve^* (T-t)^{2\kappa-1}).
\end{align}

In {\rm Case (II)},
we see that
\begin{align} \label{I1-t4-0}
    |I^{(1)}_{t,4}|
    \leq  \sum\limits_{j=1}^K
          \frac{1}{\lbb_j^2} \|\na \Phi_j\|_{L^\9} \|R\|_{L^2} \|\na R\|_{L^2}.
\end{align}
Since  $|\na \Phi_j|\leq C \sigma^{-1} \leq C \ve^*$ in {\rm Case (II)},
using \eqref{R-Tt2} and \eqref{lbb-g-t}
we have
\begin{align} \label{I1t4-esti-case2}
     |I^{(1)}_{t,4}|
     \leq C (T-t)^{-1}(\frac{\|R\|^2_{L^2}}{(T-t)^2}+\|\nabla R\|_{L^2}^2) \leq C \ve^* (T-t)^{2\kappa-1}.
\end{align}

{\it $(iv)$ Estimate of $I^{(1)}_{t,5}$.}
Since
\begin{align} \label{U-Tt}
   |U(t)| \leq C (T-t)^{-\frac d2},
\end{align}
using \eqref{f''vR2-bdd}, \eqref{G-N}, \eqref{fLp-fH1}
and  \eqref{R-Tt2}
we get
\begin{align} \label{I1t5-esti}
   | I^{(1)}_{t,5}|
   \leq& C(T-t)^{-2} (\int |U|^{\frac 4d -1} |R|^3 dx + \|R\|_{L^{\frac 4d +2}}^{\frac 4d +2})  \nonumber  \\
   \leq& C(T-t)^{-2} ((T-t)^{-2+\frac d2} \|R\|_{L^2}^{3-\frac d2} \|\na R\|_{L^2}^{\frac d2} + \|R\|_{H^1}^{\frac 4d+2})
   \leq  C(T-t)^{2\kappa}.
\end{align}

{\it $(v)$ Estimate of $I^{(1)}_{t,6}$.}
Regarding $I^{(1)}_{t,6}$,
since by \eqref{g-gz-expan} and \eqref{f'vR-bdd},
\begin{align*}
  |f(u)-f(U)|
  =|f'(U, R)\cdot R|
  \leq C(|U|^{\frac 4d} + |R|^{\frac 4d})|R|
  \leq C((T-t)^{-2} + |R|^{\frac 4d})|R|,
\end{align*}
using the integration by parts formula,
\eqref{fLp-fH1},
\eqref{R-Tt2}
and the estimate \eqref{eta-L2} of $\eta$,
we obtain
\begin{align} \label{I1t6-esti}
     |I^{(1)}_{t,6}|
     \leq& C (\|\na \eta\|_{L^2} \|\na R\|_{L^2}
              + (T-t)^{-2}\|R\|_{L^2} \|\eta\|_{L^2}  + \|R\|_{H^1}^{\frac 4d +1} \|\eta\|_{L^2} )  \nonumber \\
     \leq& C (T-t)^{2  \kappa}.
\end{align}

{\it $(vi)$ Estimate of $I^{(1)}_{t,7}$.}
The last term $I^{(1)}_{t,7}$ can be estimated similarly
as in the proof of \cite[Lemma 5.10]{SZ20}. Precisely,
using the explicit expressions \eqref{b}
of $b$,
$\sup_{0\leq t\leq T}|B_k|<\9$, a.s., $1\leq k\leq N$,
and integration by parts formula we have
\begin{align} \label{DR-bnaR.1}
       & |{\rm Im}\<\Delta R - \lbb_j^{-2} R_j + f(u) - f(U), b\cdot \na R\>|  \nonumber \\
   \leq&  C \sum\limits_{k=1}^N
          \big( | \int \na^2 \phi_k(\na \ol R, \na R)dx|
         + |\int \Delta \phi_k |\na R|^2 dx|
         + (T-t)^{-2} \|R\|_{L^2}^2  \nonumber \\
       &\qquad  + \int |R|^{\frac 4d +1} |\na \phi_k \cdot \na \ol R| dx
         + |\int (f(u) - f(U) - |R|^{\frac 4d} R) \na \phi_k \cdot \na \ol R dx | \big)    \nonumber \\
   \leq & C(T-t)^{2  \kappa}
         + C  \sum\limits_{k=1}^N  |\int (f(u)-f(U)-|R|^{\frac 4d}R)(\na \phi_k \cdot \na \ol R) dx|
\end{align}
where in the last step
we also used H\"older's inequality
and the inequality
\begin{align} \label{R-Lp}
   \|R\|_{L^p}^p
   \leq C \|R\|_{L^2}^{p+d-\frac 12 dp} \|\na R\|_{L^2}^{\frac 12 dp -d}
   \leq C (T-t)^{ \kappa p+p+d-\frac 12 dp}, \ \ \forall p\geq 2.
\end{align}
Moreover, by Lemma \ref{Lem-inter-est},
the last term in \eqref{DR-bnaR.1} is bounded by
\begin{align*}
    & C\sum\limits_{k=1}^{4/d} \sum\limits_{j=1}^K
         \int |R^k U_j^{1+\frac 4d -k} \na \phi_k \cdot \na \ol{R_j}| dx
         + C e^{-\frac{\delta}{T-t}} \\
   \leq& C \sum\limits_{j=1}^{4/d} \sum\limits_{k=1}^K
            \int \lbb_j^d |(R^k U_j ^{1+\frac 4d -k} \na \phi_k \cdot
                   \na \ol{R_j}) (\lbb_j y +\a_j)|  dy + C e^{-\frac{\delta}{T-t}}.
\end{align*}
which,
via \eqref{phik-Taylor} and \eqref{R-Lp},
can be further bounded by
\begin{align*}
   C \sum\limits_{j=1}^{4/d} \sum\limits_{k=1}^K
        (T-t)^{-\frac d2(1+\frac 4d-k)+\nu_*} \| R\|_{H^1} \|R\|^k_{L^{2k}}
   + C e^{-\frac{\delta}{T-t}}
         \leq C(T-t)^{2 \kappa}.
\end{align*}
Hence,
we obtain
\begin{align} \label{DR-bnaR}
   |{\rm Im} \<\Delta R- \lbb_j^{-2} R_j + f(u) - f(U), b\cdot \na R\>|
   \leq C(T-t)^{2 \kappa}.
\end{align}

Similarly,
since $|U(t)|\leq C (T-t)^{- \frac d2}$
and $\|c\|_{L^\9(t_*, T_*;L^\9)} <\9$,
using \eqref{R-Tt2} and \eqref{R-Lp} we get
\begin{align} \label{DR-cR}
      & |{\rm Im} \<\Delta R- \lbb_j^{-2} R_j + f(u) - f(U),c R\>|   \nonumber \\
   \leq& \|R\|^2_{H^1} + (T-t)^{-2}\|R\|_{L^2}^2
         +  \sum\limits_{k=1}^{1+\frac 4d} (T-t)^{-\frac d2(1+\frac 4d -k)} \|R\|_{L^{k+1}}^{k+1} \nonumber  \\
   \leq& C (T-t)^{2 \kappa}.
\end{align}

Thus, we conclude from \eqref{DR-bnaR} and \eqref{DR-cR} that
\begin{align} \label{I1t7-esti}
    |I^{(1)}_{t,7}|
  \leq  C(T-t)^{2 \kappa}.
\end{align}

Therefore,
plugging estimates \eqref{I1t1-esti}, \eqref{I1t2t3-esti}, \eqref{I1t4-esti-case1}, \eqref{I1t4-esti-case2},
\eqref{I1t5-esti}, \eqref{I1t6-esti} and \eqref{I1t7-esti}
into \eqref{equa-I1t}
we obtain \eqref{I1t-case1}
and finish the proof of Lemma \ref{Lem-I1t}.
\hfill $\square$

\begin{lemma} \label{Lem-I2t}
For all $t\in[t_*,T_*]$, we have
that for some $C(A)>0$,
\begin{align} \label{I2t}
 \frac{d I^{(2)}}{dt}
 \geq &  - \sum_{j=1}^{K} \frac{\gamma_j}{4\lambda_j^4}\int \Delta^2\chi_A(\frac{x-\alpha_j}{\lambda_j})|R_j|^2 dx
           + \sum_{j=1}^{K} \frac{\gamma_j}{\lambda_j^2} {\rm Re} \int \nabla^2\chi_A(\frac{x-\alpha_j}{\lambda_j})(\nabla R_j,\nabla \ol{R_j}) dx   \nonumber \\
&+ \sum_{j=1}^{K}\frac{\gamma_j}{\lambda_j}{\rm Re}\int
   \nabla\chi_A (\frac{x-\alpha_j}{\lambda_j} )\cdot\nabla \ol{U_j}
  \bigg\{\frac 2d(1+\frac 2d)|U_j|^{\frac {4}{d}-2}U_j|R_j|^2 \nonumber\\
   & \qquad + \frac 1d (1+\frac 2d) |U_j|^{\frac {4}{d}-2}\ol{U_j}R_j^2
   + \frac 1d (\frac 2d -1) |U_j|^{\frac 4d-4}U_j^3 \ol{R_j}^2 \bigg\}dx
     -C(A) (T-t)^{2 \kappa}.
\end{align}
\end{lemma}

\begin{remark}
The difficulty in the proof of Lemma \ref{Lem-I2t} lies in the analysis of
the interactions between the remainders,
of which the perturbation order is of only polynomial type.
This is different from the situation in Lemma \ref{Lem-I1t},
where the interactions involving $U_j$
are very weak,
because of the exponential decay of the ground state.
The point here is to gain additional decays
from the functions $\partial^\nu \chi_A$, where $|\nu|\geq 2$.
\end{remark}

{\bf Proof of Lemma \ref{Lem-I2t}.}
Straightforward computations show that
\begin{align}\label{equa-I2t}
   \frac{d I^{(2)}}{dt}
   =&-\sum_{j=1}^K \frac{\dot{\lambda}_{j}\gamma_{j}-\lambda_{j}\dot{\gamma}_{j}}{2\lambda_{j}^2}
           {\rm Im} \< \nabla\chi_A(\frac{x-\alpha_{j}}{\lambda_{j}})\cdot\nabla R, R_j\> \nonumber  \\
   &+\sum_{j=1}^K\frac{\gamma_{j}}{2\lambda_{j}}{\rm Im}
       \< \partial_t (\nabla\chi_A(\frac{x-\alpha_{j}}{\lambda_{j}}))\cdot\nabla R, {R}_j\>
    +\sum_{j=1}^K \frac{\gamma_{j}}{2\lambda_{j}^2}   {\rm Im}
       \< \Delta\chi_A(\frac{x-\alpha_{j}}{\lambda_{j}})R_j, \pa_t {R} \>  \nonumber \\
   &+   \sum_{j=1}^K \frac{\gamma_{j}}{2\lambda_{j}}    {\rm Im}
       \< \nabla\chi_A(\frac{x-\alpha_{j}}{\lambda_{j}})\cdot ( \nabla R_j + \na R \Phi_j), \partial_t  R \> \nonumber \\
    =:&  \sum\limits_{j=1}^K
         (I_{t,j1}^{(2)} + I_{t,j2}^{(2)} + I_{t,j3}^{(2)} + I_{t,j4}^{(2)}).
\end{align}
We shall estimate $I_{t,jk}^{(2)}$, $1\leq k\leq 4$, separately.
The main contributions come from the last two terms
$I_{t,j3}^{(2)}$ and $I_{t,j4}^{(2)}$,
which requires a delicate analysis of the interactions between remainders.

$(i)$ {\it Estimate of $ I^{(2)}_{t,j1}$ and $I^{(2)}_{t,j2}$.}
Since $\sup_{y} |\na^2 \chi_A(y)(1+|y|)| \leq C(A)$,
by Lemmas \ref{Lem-Mod-w-lbb} and \ref{Lem-P-ve-U},
\begin{align}
|\partial_t(\nabla\chi_A(\frac{x-\alpha_j}{\lambda_j}))|
=& |\nabla^2\chi_A(\frac{x-\alpha_j}{\lambda_j})\cdot
\big((\frac{x-\a_j}{\lbb_j})\cdot \frac{\dot{\lbb}_j \lbb_j +\g_j}{\lbb_j^2} \nonumber \\
  &\ \ -(\frac{x-\a_j}{\lbb_j})\cdot \frac{\g_j}{\lbb_j^2}
  + \frac{\lbb_j \dot{\a}_j - 2\beta_j}{\lbb_j^2}
  + \frac{2\beta_j}{\lbb_j^2} \big) |   \nonumber  \\
\leq&  C(A) {\lambda_j^{-2}}( Mod_j+P_j)
\leq C(A)(T-t)^{-1}.
\end{align}
Taking into account
$|\frac{\dot{\lambda}_j\gamma_j-\lambda_j\dot{\gamma}_j}{2\lambda_j^2}|
\leq C\frac{Mod_j(t)}{\lambda_j^3} \leq C(T-t)^{ \kappa}$
and \eqref{R-Tt2}
we obtain
\begin{align}  \label{esti-I2t1t2}
 |I^{(2)}_{t,j1} + I^{(2)}_{t,j2}|
 \leq& C(A) ((T-t)^{ \kappa} \|\na R\|_{L^2} \|R\|_{L^2}
             + (T-t)^{-1}  \|\na R\|_{L^2} \|R\|_{L^2} ) \nonumber \\
 \leq& C(A) (T-t)^{2 \kappa}.
\end{align}

$(ii)$ {\it Estimate of $ I^{(2)}_{t,j3}$.}
We claim that
\begin{align} \label{I2t3.1-esti}
 I^{(2)}_{t,j3}
 =& -\frac{\gamma_j}{4\lambda_j^4}{\rm Re}\int \Delta^2\chi_A (\frac{x-\alpha_j}{\lambda_j} )|R_j|^2 dx+\frac{\gamma_j}{2\lambda_j^2}{\rm Re}\int \Delta\chi_A(\frac{x-\alpha_j}{\lambda_j})|\nabla R_j|^2 dx  \nonumber \\
&-\frac{\gamma_{j}}{2\lambda_{j}^2}
   {\rm Re} \< \Delta\chi_A(\frac{x-\alpha_{j}}{\lambda_{j}})R_j, f'(U_j)\cdot R_j \>
   + \calo((T-t)^{2 \kappa}).
\end{align}

In order to prove \eqref{I2t3.1-esti},
we infer from \eqref{g-gzz-expan} and equation \eqref{equa-R} that
\begin{align} \label{Deltachi-esti}
 I^{(2)}_{t,j3} =& - \frac{\gamma_{j}}{2\lambda_{j}^2}
    {\rm Re}  \< \Delta\chi_A(\frac{x-\alpha_{j}}{\lambda_{j}})R_j,
       \Delta R+f'(U)\cdot R +f''(U,R)\cdot R^2+(b\cdot \na +c) R+\eta\>.
\end{align}
The main contributions come from the terms involving $\Delta R$ and $f'(U)\cdot R$.

First,
since for any $j\not =l$,
$|x-\a_j|\geq 4 \sigma$ on the support of $R_l$,
taking $t_*$ close to $T$
we may let $|x_j-\a_j|\leq \sigma$,
and thus $|x - \a_j|\geq 3\sigma$
on the support $R_l$.
By the integration by parts formula,
for $1\leq j\neq l\leq K$,
\begin{align} \label{Deltachi.1-esti}
&|\frac{\gamma_{j}}{2\lambda_{j}^2}{\rm Re}\int  \Delta\chi_A(\frac{x-\alpha_{j}}{\lambda_{j}})R_j\Delta\overline{ R_l}dx|   \nonumber \\
\leq&  |\frac{\gamma_{j}}{2\lambda_{j}^2}{\rm Re}\int_{|x-\al_j|\geq 3\sigma}  \Delta\chi_A(\frac{x-\alpha_{j}}{\lambda_{j}})\nabla R_j\cdot\nabla\overline{ R_l}dx| \nonumber \\
   & + |\frac{\gamma_{j}}{2\lambda_{j}^3}{\rm Re}\int_{|x-\al_j|\geq 3\sigma}
      \nabla\Delta\chi_A(\frac{x-\alpha_{j}}{\lambda_{j}})\cdot\nabla\overline{ R_l}R_jdx| \nonumber \\
=:& K_1 + K_2.
\end{align}

The key observation here is that,
because of the decay of  $\pa^\nu \chi_A$, $|\nu|\geq 2$,
the different remainders $R_j$ and $R_l$
have weak interactions of order
$(T-t)^{2\kappa}$,
which is important for the bootstrap estimate of the remainder.

To be precise,
since
$$\Delta \chi(y) = \psi''(|y|) + (d-1)\psi'(|y|)|y|^{-1} \leq C |y|^{-1},  \ \
if\ |y|\geq 2,$$
we have
\begin{align} \label{K1-esti}
   K_1
   \leq C \frac{\g_j}{ \lbb_j^2}
          (\frac{\lbb_j A}{3\sigma}) \| R\|^2_{H^1}
   \leq CA (T-t)^{2 \kappa}.
\end{align}
Similarly,
since
$$\partial_j \Delta \chi(y) = \psi'''(|y|) \frac{y_j}{|y|}
 + (d-1) (\psi''(|y|) \frac{y_j}{|y|^2} - \psi'(|y|) \frac{y_j}{|y|^3})
 \leq C |y|^{-2}, \ \
if\ |y|\geq 2,$$
we get
\begin{align} \label{K2-esti}
    K_2 \leq C \frac{\g_j}{ \lbb_j^3} (\frac{\lbb_j A}{3\sigma})^2 \|\na R\|_{L^2} \|R_j\|_{L^2}
    \leq C  (T-t)^{2 \kappa+1}.
\end{align}
Hence,
plugging \eqref{K1-esti} and \eqref{K2-esti} into \eqref{Deltachi.1-esti}
we conclude that
the interactions between different remainders $R_j$ and $R_l$
have the negligible order $(T-t)^{2 \kappa}$, i.e.,
\begin{align*}
   |\frac{\gamma_{j}}{2\lambda_{j}^2}{\rm Re}\int  \Delta\chi_A(\frac{x-\alpha_{j}}{\lambda_{j}})R_j\Delta\overline{ R_l}dx|
   \leq C (T-t)^{2 \kappa}.
\end{align*}
This along with the integration by parts formula yields that
\begin{align} \label{Dchi-Rj-DR}
    -\frac{\g_j}{2\lbb_j^2} {\rm Re} \< \Delta\chi_A(\frac{x-\alpha_{j}}{\lambda_{j}})R_j, \Delta R\>
 =& - \frac{\g_j}{4\lbb_j^4} {\rm Re} \int \Delta^2 \chi_A(\frac{x - \a_j}{\lbb_j}) |R_j|^2 dx  \nonumber  \\
  & + \frac{\g_j}{2\lbb_j^2} {\rm Re} \int \Delta \chi_A(\frac{x - \a_j}{\lbb_j}) |\na R_j|^2 dx
   + \calo( (T-t)^{2 \kappa}).
\end{align}

We also apply Lemma \ref{Lem-inter-est} to obtain
\begin{align} \label{Dchi-Rj-f'UR}
   {\rm Re} \< \Delta\chi_A(\frac{x-\alpha_{j}}{\lambda_{j}})R_j, f'(U)\cdot R \>
   = {\rm Re} \< \Delta\chi_A(\frac{x-\alpha_{j}}{\lambda_{j}})R_j, f'(U_j)\cdot R_j\> + \calo(e^{-\frac{\delta}{T-t}}).
\end{align}

Moreover,
since by \eqref{f''vR2-bdd} and \eqref{U-Tt},
\begin{align*}
  |f''(U,R)\cdot R^2|
  \leq C (|U|^{\frac 4d-1} + |R|^{\frac 4d-1}) |R|^2
  \leq C ((T-t)^{-2+\frac d2} + |R|^{\frac 4d-1})|R|^2,
\end{align*}
using \eqref{fghLp-H1},
\eqref{R-Tt2} and \eqref{R-Lp}
we have
\begin{align} \label{Dchi-Rj-f''UR2}
      & |\frac{\g_j}{2\lbb_j^2} {\rm Re} \< \Delta\chi_A(\frac{x-\alpha_{j}}{\lambda_{j}})R_j, f''(U,R)\cdot R^2\>|  \nonumber \\
\leq& C(A)  \int (T-t)^{-1} ((T-t)^{-2+\frac d2} + |R|^{\frac 4d -1})|R|^3 dx \nonumber \\
\leq& C(A) ((T-t)^{-3+ \frac d2}\|R\|_{L^3}^3
           + (T-t)^{-1}\|R\|_{H^1}^{\frac 4d+2}  )
\leq  C(A)(T-t)^{2  \kappa}.
\end{align}

Furthermore,
using H\"older's inequality, \eqref{R-Tt2} and \eqref{eta-L2} we have
\begin{align} \label{Dchi-Rj-bcR-eta}
       & |\frac{\g_j}{2\lbb_j^2} {\rm Re} \< \Delta\chi_A(\frac{x-\alpha_{j}}{\lambda_{j}})R_j, (b\cdot \na +c)R +\eta \>|   \nonumber \\
  \leq & C(A)(T-t)^{-1}(\|R\|_{L^2}\|\nabla R\|_{L^2}+\|R\|_{L^2}^2)
         + C(A) (T-t)^{-1} \|\eta\|_{L^2} \|R\|_{L^2}   \nonumber \\
  \leq& C(A) (T-t)^{2 \kappa}.
\end{align}

Hence,
plugging \eqref{Dchi-Rj-DR}, \eqref{Dchi-Rj-f'UR},
\eqref{Dchi-Rj-f''UR2} and \eqref{Dchi-Rj-bcR-eta} into \eqref{Deltachi-esti}
we obtain \eqref{I2t3.1-esti}, as claimed.

$(i)$ {\it Estimate of $ I^{(2)}_{t,j4}$.}
We claim that
\begin{align}   \label{I2t3.2-esti}
 I^{(2)}_{t,j4} =&\frac{\gamma_j}{\lambda_j^2}{\rm Re}\int \nabla^2\chi_A(\frac{x-\alpha_j}{\lambda_j})(\nabla R_j,\nabla \ol{R_j}) dx
-\frac{\gamma_j}{2\lambda_j^2}{\rm Re}\int \Delta\chi_A(\frac{x-\alpha_j}{\lambda_j})|\nabla R_j|^2 dx \nonumber  \\
&-\frac{\gamma_{j}}{\lambda_{j}}\< \nabla\chi_A(\frac{x-\alpha_{j}}{\lambda_{j}})\cdot
         \nabla R_j,f^\prime(U_j)\cdot R_j\>
         + \calo((T-t)^{2 \kappa}).
\end{align}

For this purpose,
we infer from equation \eqref{equa-R} and Lemma \ref{Lem-inter-est}  that
\begin{align} \label{I2t3.2-esti.0}
I^{(2)}_{t,j4} =& -\frac{\gamma_{j}}{2\lambda_{j}} {\rm Re}
     \< \nabla\chi_A(\frac{x-\alpha_{j}}{\lambda_{j}})\cdot
         (\nabla R_j + \nabla R\Phi_j), \nonumber \\
 & \qquad \qquad \ \ \Delta R+f^\prime(U_j)\cdot R_j
+f''(U,R)\cdot R^2+(b\cdot \na +c) R+\eta\>
   + \calo(e^{-\frac{\delta}{T-t}}).
\end{align}

We first show that
\begin{align} \label{nachi-naR-DR}
   & - {\rm Re}\< \nabla\chi_A(\frac{x-\alpha_{j}}{\lambda_{j}})\cdot
        (\nabla R_j + \nabla R \Phi_j ),  \Delta {R} \>    \\
 =&  {\rm Re} \int \frac{2}{\lbb_j} \na^2 \chi_A (\frac{x-\alpha_{j}}{\lambda_{j}}) (\na R_j, \na \ol{R_j})
     - \frac{1}{\lbb_j} \Delta \chi_A (\frac{x-\alpha_{j}}{\lambda_{j}}) |\na R_j|^2 dx
     + \calo((T-t)^{2 \kappa}). \nonumber
\end{align}
In order to prove \eqref{nachi-naR-DR},
using integration by parts formula we see that
\begin{align} \label{nachi-naR-DR.0}
     & - {\rm Re} \< \nabla\chi_A(\frac{x-\alpha_{j}}{\lambda_{j}})\cdot\nabla R_j,  \Delta {R} \> \nonumber \\
   =& {\rm Re} \int \frac{1}{\lbb_j} \na^2 \chi_A(\frac{x-\alpha_{j}}{\lambda_{j}}) (\na R_j, \na \ol{R})
      - \frac{1}{\lbb_j} \Delta \chi_A(\frac{x-\alpha_{j}}{\lambda_{j}}) \na R_j \cdot \na \ol{R}  dx  \nonumber  \\
    &  - \sum\limits_{1\leq k,l\leq d} {\rm Re} \int \partial_k \chi_A(\frac{x-\alpha_{j}}{\lambda_{j}}) \partial_l R_j \partial_{kl} \ol{R} dx.
\end{align}
Then,  as in the proof of \eqref{Dchi-Rj-DR},
using the decay of $\pa_\nu \chi_A$ with $|\nu|=2$
we obtain that
the interactions between different remainders have negligible contributions
and thus
\begin{align} \label{nachi-naR-DR.1}
      &- {\rm Re} \< \nabla\chi_A(\frac{x-\alpha_{j}}{\lambda_{j}})\cdot\nabla R_j , \Delta {R}  \>    \nonumber \\
    =&  {\rm Re} \int \frac{1}{\lbb_j} \na^2 \chi_A (\frac{x-\alpha_{j}}{\lambda_{j}})(\na R_j, \na \ol{R_j})
      -  \frac{1}{\lbb_j} \Delta \chi_A  (\frac{x-\alpha_{j}}{\lambda_{j}})|\na {R_j}|^2 dx  \nonumber \\
     &  - \sum\limits_{1\leq k,l\leq d} {\rm Re} \int \partial_k \chi_A(\frac{x-\alpha_{j}}{\lambda_{j}}) \partial_l R_j \partial_{kl} \ol{R} dx
      + \calo((T-t)^{2\kappa}).
\end{align}
Similarly, we have
\begin{align} \label{nachi-naR-DR.2}
    & - {\rm Re} \<  \na \chi_A(\frac{x-\alpha_{j}}{\lambda_{j}}) \cdot \nabla R \Phi_j, \Delta{R} \> \nonumber  \\
   =&{\rm Re} \int \frac{1}{\lbb_j} \na^2 \chi_A(\frac{x-\alpha_{j}}{\lambda_{j}}) (\na R, \na \ol{R}) \Phi_j
     + \sum\limits_{1\leq k,l\leq d}  \pa_k \chi_A(\frac{x-\alpha_{j}}{\lambda_{j}}) \pa_{kl} R \pa_l \ol{R} \Phi_j dx\nonumber \\
    & + {\rm Re} \int (\na \chi_A(\frac{x-\alpha_{j}}{\lambda_{j}}) \cdot \na R) (\na \ol{R} \cdot \na \Phi_j)dx    \nonumber   \\
   =&   {\rm Re} \int \frac{1}{\lbb_j} \na^2 \chi_A(\frac{x-\alpha_{j}}{\lambda_{j}})(\na R_j, \na \ol{R_j})
     + \sum\limits_{1\leq k,l\leq d}   \pa_k \chi_A(\frac{x-\alpha_{j}}{\lambda_{j}}) \pa_{kl} \ol{R} \pa_l {R} \Phi_j  dx  \nonumber \\
     & + \calo((T-t)^{2 \kappa}).
\end{align}
Moreover,
for the two terms involving $\partial_{kl}\ol R$ in \eqref{nachi-naR-DR.1} and \eqref{nachi-naR-DR.2},
we see that the cancellation appears
and the integration by parts formula and \eqref{R-Tt} give
\begin{align}  \label{nachi-naR-DR.3}
  & {\rm Re} \int \pa_k \chi_A(\frac{x-\alpha_{j}}{\lambda_{j}}) \pa_{kl} \ol{R} \pa_l {R} \Phi_j
   -   \partial_k \chi_A(\frac{x-\alpha_{j}}{\lambda_{j}}) \partial_l R_j \partial_{kl} \ol{R} dx  \nonumber \\
  =& - {\rm Re} \int \pa_k \chi_A(\frac{x-\alpha_{j}}{\lambda_{j}}) R \partial_l \Phi_j \pa_{kl} \ol{R}  d x    \nonumber \\
  =&  {\rm Re} \int \frac{1}{\lbb_j} \Delta \chi_A(\frac{x-\alpha_{j}}{\lambda_{j}}) \pa_l R \ol{R} \pa_l \Phi_j
     +  \pa_k \chi_A(\frac{x-\alpha_{j}}{\lambda_{j}}) \pa_l \Phi_j \pa_l R \pa_k \ol{R} dx \nonumber \\
   &  + {\rm Re} \int  \pa_k \chi_A(\frac{x-\alpha_{j}}{\lambda_{j}}) \pa_{kl} \Phi_j  \pa_l R \ol{R}  dx
  = \calo((T-t)^{2 \kappa}).
\end{align}
Thus,
plugging \eqref{nachi-naR-DR.1}-\eqref{nachi-naR-DR.3}
into \eqref{nachi-naR-DR.0}
we obtain \eqref{nachi-naR-DR}, as claimed.

We also apply Lemma \ref{Lem-inter-est}
to decouple different profiles between $\{U_j\}$ and $\{R_j\}$ to obtain that,
similarly to \eqref{Dchi-Rj-f'UR},
\begin{align*}
     &\frac{\gamma_{j}}{2\lambda_{j}} {\rm Re}
     \< \nabla\chi_A(\frac{x-\alpha_{j}}{\lambda_{j}})\cdot
         (\nabla R_j + \nabla R\Phi_j), f'(U)\cdot R\>  \nonumber \\
    =& \frac{\gamma_{j}}{2\lambda_{j}} {\rm Re}
     \< \nabla\chi_A(\frac{x-\alpha_{j}}{\lambda_{j}})\cdot
         (\nabla R_j + \nabla R_j \Phi_j), f'(U_j)\cdot R_j\>
         + \calo(e^{-\frac{\delta}{T-t}}).
\end{align*}
Since $|x-\a_j|\geq 3\sigma$ on the support of $1-\Phi_j$,
we have the exponential decay of $U_j$ that
$|U_j(x)| \leq C \lbb_j^{-\frac d 2} e^{-\delta \frac{3\sigma}{\lbb_j}}$,
and thus
\begin{align*}
   {\rm Re}  \< \nabla\chi_A(\frac{x-\alpha_{j}}{\lambda_{j}})
       \cdot \nabla R \Phi_j), f'(U_j)\cdot R_j\>
   =  {\rm Re}  \< \nabla\chi_A(\frac{x-\alpha_{j}}{\lambda_{j}}) \cdot \nabla R, f'(U_j)\cdot R_j\>
      + \mathcal{O}(e^{-\frac{\delta}{T-t}}).
\end{align*}
Hence, we obtain
\begin{align} \label{nachi-naR-f'UR}
     &\frac{\gamma_{j}}{2\lambda_{j}} {\rm Re}
     \< \nabla\chi_A(\frac{x-\alpha_{j}}{\lambda_{j}})\cdot
         (\nabla R_j + \nabla R\Phi_j), f'(U)\cdot R\>  \nonumber \\
    =& \frac{\gamma_{j}}{\lambda_{j}} {\rm Re}
     \< \nabla\chi_A(\frac{x-\alpha_{j}}{\lambda_{j}})\cdot
         \nabla R_j, f'(U_j)\cdot R_j\>
         + \calo(e^{-\frac{\delta}{T-t}}).
\end{align}

Moreover,
using H\"older's inequality,
\eqref{R-Tt2}
and \eqref{eta-L2} we easily get
\begin{align} \label{nachi-naR-f''UR2}
   |\frac{\gamma_{j}}{2\lambda_{j}} {\rm Re}
     \< \nabla\chi_A(\frac{x-\alpha_{j}}{\lambda_{j}})\cdot
         (\nabla R_j + \nabla R\Phi_j),
         f''(U)\cdot R^2 + (b\cdot \na +c)R+\eta\>|
   \leq C(T-t)^{2 \kappa}.
\end{align}

Hence, we conclude from \eqref{nachi-naR-DR}, \eqref{nachi-naR-f'UR}
and \eqref{nachi-naR-f''UR2} that \eqref{I2t3.2-esti} holds.

Now,
putting the estimates \eqref{esti-I2t1t2},
\eqref{I2t3.1-esti} and \eqref{I2t3.2-esti} altogether
we obtain
\begin{align*}
  \frac{dI^{(2)}}{dt}
 =&- \sum_{j=1}^K \frac{\gamma_j}{4\lambda_j^4}\int \Delta^2\chi_A (\frac{x-\alpha_j}{\lambda_j})|R_j|^2 dx
  +\sum_{j=1}^K\frac{\gamma_j}{\lambda_j^2}{\rm Re}\int \nabla^2\chi_A(\frac{x-\alpha_j}{\lambda_j})(\nabla R_j,\nabla \ol{R_j}) dx  \nonumber \\
& -\sum_{j=1}^K{\rm Re} \<\frac{\gamma_j}{2\lambda_j^2}\Delta\chi_A(\frac{x-\alpha_j}{\lambda_j}) R_j
                + \frac{\gamma_j}{\lambda_j}\nabla\chi_A(\frac{x-\alpha_j}{\lambda_j})\cdot  \nabla R_j,
                  f'(U_j) \cdot R_j \>
  + \calo((T-t)^{2 \kappa}) .
\end{align*}
Because the profiles are decoupled completely,
we can treat each profile individually by using similar computations as in the proof of \cite[Lemma 5.11]{SZ20}
and thus obtain that
the third term on the right-hand-side above is equal to
\begin{align} \label{esti-I2.31}
  & \sum_{j=1}^{K}\frac{\gamma_j}{\lambda_j}{\rm Re}\int \nabla\chi_A (\frac{x-\alpha_j}{\lambda_j} )\cdot\nabla \ol{U_j}
   \bigg\{\frac 2d(1+\frac 2d) |U_j|^{\frac {4}{d}-2}U_j|R_j|^2   \nonumber \\
  & \qquad \qquad \qquad   +\frac 1d(1+\frac 2d) |U_j|^{\frac {4}{d}-2}\ol{U_j} R_j^2
   +\frac 1d(\frac 2d-1) |U_j|^{\frac 4d-4}U_j^3\overline{R_j}^2\bigg\} dx,
\end{align}
which immediately yields \eqref{I2t},
thereby finishing the proof of Lemma \ref{Lem-I2t}.
\hfill $\square$

We are now ready to prove Theorem \ref{Thm-I-mono}.

{\bf Proof of Theorem \ref{Thm-I-mono}.}
At this stage,
the blow-up profiles are decoupled in \eqref{I1t-case1} and \eqref{I2t},
up to the acceptance order $\calo((T-t)^{2 \kappa})$,
and thus we are able to treat each profile separately
by using   similar arguments as in the proof of \cite[Theorem 5.8]{SZ20}.
For the reader's convenience,
we sketch the proof below.

Combining \eqref{I1t-case1} and \eqref{I2t} altogether
and then using
the renormalized variable $\ve_j$ in \eqref{Rj-ej}
we obtain that for all $t\in[t_*,T_*]$,
\begin{align} \label{dI}
\frac{dI}{dt}
   \geq&
   \sum\limits_{j=1}^K
   \frac{\gamma_j}{\lambda_j^4}
   \big(\int \nabla^2\chi (\frac{y}{A} )(\nabla \varepsilon_j,\nabla \ol{\varepsilon_j}) dy
          + \int |\varepsilon_j|^2dy
          -\int (1+\frac 4d)Q^\frac{4}{d}\varepsilon_{j,1}^2+Q^\frac{4}{d}\varepsilon_{j,2}^2  dy \nonumber \\
&  \qquad \qquad -\frac{1}{4A^2}\int \Delta^2\chi (\frac{y}{A} )|\varepsilon_j|^2 dy\big) \nonumber \\
&+ \sum\limits_{j=1}^K  \frac 2d  \frac{\gamma_j}{\lambda_j^4}\int (A\nabla\chi (\frac{y}{A} )-y ) \cdot \nabla Q
Q^{\frac 4d -1} ((1+\frac 4d)\varepsilon_{j,1}^2+ \varepsilon_{j,2}^2)dy \nonumber \\
& -C \ve^* (T-t)^{2\kappa-1} -C(A)(T-t)^{2 \kappa},
\end{align}
where $\ve_{j,1}$ and $\ve_{j,2}$ denote the real and imaginary parts of $\ve_j$, respectively.

Then,
since
$$\int \na^2 \chi(\frac yA)(\na \ve_j, \na \ol{\ve_j}) dy
\geq \int \psi''(|\frac yA|) |\na \ve_j|^2 dy,$$
applying Corollary \ref{Cor-coer-f-local} with $\phi(x) :=\psi''(|x|)$,
Lemma \ref{Lem-almost-orth} and Proposition \ref{Prop-mass-local},
and using the estimate
$Scal(\ve_j) \leq C (T-t)^{2\kappa +4}$
we obtain for some $\wt C>0$
\begin{align} \label{dIdt-esti}
\frac{dI}{dt}
\geq& \wt C \sum\limits_{j=1}^K \frac{\gamma_j}{\lambda_j^4} \int \psi'' (|\frac{y}{A}| )(|\varepsilon_j|^2+|\nabla \varepsilon_j|^2) dy
    -  \sum\limits_{j=1}^K \frac{1}{4A^2} \frac{\g_j}{\lbb_j^4} \int \Delta^2\chi (\frac{y}{A} )|\varepsilon_j|^2 dy  \nonumber \\
&+  \sum\limits_{j=1}^K\frac 2d \frac{\gamma_j}{\lambda_j^4}\int (A\nabla\chi (\frac{y}{A} )-y ) \cdot \nabla Q
Q^{\frac 4d -1} ((1+\frac 4d) \varepsilon_{j,1}^2+ \varepsilon_{j,2}^2)dy   \nonumber \\
&  -C \ve^* (T-t)^{2\kappa-1}   -C(A)(T-t)^{2 \kappa}.
\end{align}
Taking into account that
for $A$ large enough,
$$\frac{1}{4A^2} |\Delta^2\chi (\frac{y}{A} )| \leq \frac {1}{4} \wt C \psi'' (|\frac{y}{A}| ), $$
and
$$\frac 2d(2+\frac 4d) |A\nabla\chi (\frac{y}{A} )-y ||\na QQ^{\frac 4d -1}| \leq \frac {1}{4} \wt C \psi'' (|\frac{y}{A}| ), $$
we arrive at
\begin{align*}
  \frac{dI}{dt}
\geq& \frac{1}{2} \wt C \sum\limits_{j=1}^K
        \frac{\gamma_j}{\lambda_j^4} \int \psi'' (|\frac{y}{A}| )(|\varepsilon_j|^2+|\nabla \varepsilon_j|^2) dy
     -C \ve^* (T-t)^{2\kappa-1}  - C(A) (T-t)^{2 \kappa}.
\end{align*}

Therefore,
as $\psi''(r) \geq \delta e^{-r}$ for some $\delta>0$,
we obtain   \eqref{dIt-mono-case2}
and finish the proof.  \hfill $\square$

\subsection{Proof of bootstrap estimates}  \label{Subsec-Boot-proof}

In this subsection
we prove the crucial bootstrap estimates in Theorem \ref{Thm-u-Boot}.
To begin with,
we first obtain the refined estimate for the modulation parameter $\beta$.

\begin{lemma}[Refined estimate for $\beta$] \label{Lem-betaj}
There exists $C>0$ such that
for all  $t\in[t_*,T_*]$,
\begin{align}  \label{b-g-ve-lbb}
\sum_{j=1}^{K}|\beta_{j}(t)|^2
\leq C\sum_{j=1}^{K}|\omega_j^2\lambda_{j}^2(t)-\g_{j}^2(t)|+ C(T-t)^{\kappa+3}.
\end{align}
\end{lemma}

\begin{remark}
Unlike single
bubble case in \cite{SZ20},
the proof of Lemma \ref{Lem-betaj}
requires the localized mass in Proposition \ref{Prop-mass-local}
and also a delicate treatment of the localized function $\Phi_j$
and the radial function $\phi_A$ in Corollary \ref{Cor-coer-f-local},
in order to derive the coercivity of energy.
\end{remark}

{\bf Proof of Lemma \ref{Lem-betaj}.}
Using the expansion \eqref{F-Taylor*} of $F(u) = \frac{d}{2d+4} |u|^{2+\frac 4d}$
we have
\begin{align}
E(u) =& \frac{1}{2}\int|\nabla U|^2dx - \frac{d}{2d+4} \int |U|^{2+\frac 4d} dx
   -  {\rm Re}\int ({\Delta U+|U|^{\frac 4d}U}) \overline{R} dx  \nonumber \\
&+\frac{1}{2} {\rm Re} \int |\nabla R|^2-(1+\frac2d)|U|^{\frac 4d} |R|^2
   - \frac 2d |U|^{\frac 4d-2}U^2\overline{R}^2 dx+ o(\|R\|^2_{H^1}).
\end{align}
Note that,
by Lemma \ref{Lem-inter-est} and the explicit expression \eqref{Uj-Qj-Q} of $U_j$,
\begin{align}
 &\frac{1}{2}\int|\nabla U|^2dx - \frac{d}{2d+4} \int |U|^{2+\frac 4d}dx \nonumber \\
=&\sum_{j=1}^K(\frac{1}{2}\int|\nabla U_{j}|^2dx - \frac{d}{2d+4} \int |U_{j}|^{2+\frac 4d} dx)
   + \calo(e^{-\frac{\delta}{T-t}})  \nonumber \\
=&\sum_{j=1}^K( \frac{|\beta_{j}|^2}{2\lambda_{j}^2} \|Q\|^2_{L^2}
+ \frac{\gamma_{j}^2}{8\lambda_{j}^2} \|yQ\|_{L^2}^2)
+  \calo(e^{-\frac{\delta}{T-t}}) ,
\end{align}
and
\begin{align}
  {\rm Re}\int (\Delta U+|U|^{\frac 4d}U)\overline{R} dx
= \sum_{j=1}^K{\rm Re}\int (\Delta U_{j}+|U_{j}|^{\frac 4d}U_{j})\overline{R_j} dx
+ \calo(e^{-\frac{\delta}{T-t}}\|R\|_{L^2}).
\end{align}
Taking into account Proposition \ref{Prop-mass-local}
and rearranging the terms according to the orders of $R$
we obtain that for $T$ small enough,
\begin{align} \label{Eu-esti1}
E(u)
=& E(u)+\sum_{j=1}^{K}\frac{{1}}{\lbb_{j}^2} {\rm Re} \int \ol{U_j} R_j
       + \frac12 |R|^2\Phi_jdx
       + \calo((T-t)^{2 \kappa})   \nonumber  \\
=& \sum_{j=1}^K( \frac{|\beta_{j}|^2}{2\lambda_{j}^2} \|Q\|_{L^2}^2
+ \frac{\gamma_{j}^2}{8\lambda_{j}^2} \|yQ\|_{L^2}^2)
 -  \sum_{j=1}^{K} {\rm Re} \int ({\Delta U_{j}-\frac{1}{\lbb_{j}^2} U_{j}+|U_{j}|^{\frac 4d}U_{j}}) \overline{R_{j}} dx  \nonumber  \\
&+\frac{1}{2}{\rm Re} \int |\nabla R|^2+\sum_{j=1}^{K}\frac{1}{\lbb_{j}^2}|R|^2  \Phi_j
     -(1+\frac2d)|U|^{\frac 4d} |R|^2
     - \frac 2d |U|^{\frac 4d-2}U^2\overline{R}^2 dx  \nonumber \\
&  + \calo((T-t)^{2 \kappa}).
\end{align}

On one hand,
using the identity \eqref{DQj-Qj-fQj}
and the change of variables
we get
\begin{align}
  &  {\rm Re} \int ({\Delta U_{j}-\frac{1}{\lbb_{j}^2} U_{j}+|U_{j}|^{\frac 4d}U_{j}}) \overline{R_{j}} dx  \nonumber \\
= &  \frac{1}{\lambda_{j}^2} {\rm Im} \int(\gamma_{j}\Lambda Q_{j} -2\beta_{j}\cdot\nabla Q_{j}) \ol{\varepsilon_j} dx
    + \frac{1}{\lbb_{j}^2} {\rm Re} \int |\beta_{j} - \frac {\gamma_{j}}{ 2} y|^2 Q_{j} \ol{\varepsilon_j} dx,
\end{align}
which along with the almost orthogonality in Lemma \ref{Lem-almost-orth}
and \eqref{ab-Tt2} yields that
\begin{align}  \label{Eu-esti2}
  |\sum_{j=1}^{K} {\rm Re} \int ({\Delta U_{j}-\frac{1}{\lbb_{j}^2} U_{j}+|U_{j}|^{\frac 4d}U_{j}}) \overline{R_{j}} dx|
  \leq C \|R\|_{L^2}
  \leq C(T-t)^{2 \kappa+1}.
\end{align}

On the other hand,
we claim that
there exist $\wt c, C>0$ such that
for the quadratic terms of $R$ on the right-hand side of \eqref{Eu-esti1},

\begin{align} \label{R-quad-decoup}
E_2(u) :=& {\rm Re} \int |\nabla R|^2+\sum_{j=1}^{K}\frac{1}{\lbb_{j}^2}|R|^2  \Phi_j
     -(1+\frac2d)|U|^{\frac 4d} |R|^2
     - \frac 2d |U|^{\frac 4d-2}U^2\overline{R}^2 dx  \nonumber \\
   \geq& \wt c ( \int |\nabla R|^2+\sum_{j=1}^{K}\frac{1}{\lbb_{j}^2}|R|^2 \Phi_j dx) \nonumber \\
       & -C ( (T-t)^{-1} \|R\|_{L^2}^2
         +  e^{-\frac{\delta}{T-t}} \|R\|_{H^1}^2
         + (T-t)^{2\kappa+2}).
\end{align}
(Note that,
this does not follow directly from Corollary \ref{Cor-coer-f-local},
because the localized function $\Phi_j$ does not satisfy the conditions there.)

For this purpose,
using the partition of unity
and Lemma \ref{Lem-inter-est}
we have
\begin{align*}
 E_2(u)
 =& \sum_{j=1}^{K}{\rm Re}\int (|\nabla R|^2
          +\frac{1}{\lbb_{j}^2}|R|^2)\Phi_j
          -(1+\frac2d)|U_{j}|^{\frac 4d}|R|^2
          -\frac 2d |U_{j}|^{\frac 4d-2}U_{j}^2\overline{R}^2 dx \nonumber \\
   & +\calo( e^{-\frac{\delta}{T-t}}\|R\|_{L^2}^2).
\end{align*}
In order to obtain the coercivity of the energy,
we use $\phi_{A,j}(x):=\phi_A(\frac{x-\alpha_j}{\lambda_{j}})$ with $\phi_A(x)$ as in Corollary \ref{Cor-coer-f-local}
and the renormalized variable $\wt \ve_j$
defined by
$$R(t,x)=\lbb_{j}^{-\frac d2} \wt \ve_{j} (t,\frac{x-\a_j}{\lbb_j}) e^{i \theta_j},$$
to reformulate $E_2(u)$ as follows
\begin{align} \label{R-quad-decoup-2}
 E_2(u)=
  & \sum\limits_{j=1}^K
  {\rm Re}\int (|\nabla R|^2+\frac{1}{\lbb_{j}^2} |R|^2)\phi_{A,j}
    -(1+\frac2d)|U_{j}|^{\frac 4d} |R|^2
    -\frac 2d |U_{j}|^{\frac 4d-2}U_{j}^2\overline{R}^2 dx \nonumber \\
& +\sum\limits_{j=1}^K \int (|\nabla R|^2+\frac{1}{\lbb_{j}^2}|R|^2)(\Phi_j-\phi_{A,j})dx
   +\calo(e^{-\frac{\delta}{T-t}}\|R\|_{L^2}^2) \nonumber \\
  =& \sum\limits_{j=1}^K  \frac{1}{\lbb_{j}^2}
     {\rm Re}\int (|\nabla \wt \ve_{j}|^2+|\wt \ve_{j}|^2)\phi_{A}
     -(1+\frac2d) Q^{\frac 4d}|\wt \ve_{j}|^2
     - \frac 2d  Q^{\frac 4d-2}\overline{Q_{j}}^2\wt \ve_{j}^2 dy  \nonumber  \\
&+\sum\limits_{j=1}^K \frac{1}{\lbb_{j}^2} \int (|\nabla \wt \ve_{j}|^2+|\wt \ve_{j}|^2)(\Phi_j(\lbb_{j}y+\alpha_{j})
-\phi_{A}(y))dy
  +\calo(e^{-\frac{\delta}{T-t}}\|R\|_{L^2}^2) \nonumber \\
=&: \sum\limits_{j=1}^K E_{21,j} + \sum\limits_{j=1}^K  E_{22,j}
+ \calo(e^{-\frac{\delta}{T-t}}\|R\|_{L^2}^2).
\end{align}

Since $Q_j= Q+ \calo(P\<y\>^2 Q)$,
the  localized coercivity in Corollary \ref{Cor-coer-f-local},
the estimate $Scal(\wt \ve_j) \leq C(T-t)^{2\kappa +4}$
and \eqref{R-Tt} yield immediately that
there exists $\wt c_j>0$ such that
\begin{align} \label{t1}
   E_{21,j}
   \geq& \frac{1}{\lbb^2_j}{\rm Re} \int (|\na \wt \ve_j|^2 + |\wt \ve_j|^2) \phi_A
      -(1+\frac 2d) Q^\frac 4d |\wt \ve_j|^2
      - \frac 2d Q^ \frac 4d \wt \ve_j^2 dy
      - C (T-t)^{-1} \|R\|_{L^2}^2 \nonumber \\
   \geq&  \wt c_j\int (|\nabla R|^2+\frac{1}{\lbb_{j}^2}|R|^2)\phi_{A, j} dx
          - C \frac{1}{\lbb^2_j} Scal(\wt \ve_j) - C (T-t)^{-1}\|R\|_{L^2}^2  \nonumber \\
   \geq&  \wt c_j\int (|\nabla R|^2+\frac{1}{\lbb_{j}^2}|R|^2)\phi_{A, j} dx
          - C(T-t)^{2\kappa+2} - C (T-t)^{-1}\|R\|_{L^2}^2.
\end{align}

Moreover,
set $\wt c:=\min\{\half, \wt c_j, 1\leq j\leq K\}>0$.
Since $\Phi_j(\lbb_{j}y+\alpha_{j}) - \phi_{A}(y)) \geq 0$
if $|y\cdot  {\bf v_1}| \leq \frac{3\sigma}{\lbb_j}$,
we get
\begin{align*}
 E_{22,j}
 \geq&\frac{\wt c}{\lbb_{j}^2} \int_{|y\cdot  {\bf v_1}|\leq \frac{3\sigma}{\lbb_{j}}} (|\nabla \wt \ve_j|^2+|\wt \ve_j|^2)
     (\Phi_j(\lbb_{j}y+\alpha_{j}) -\phi_{A}(y))dy\\
&+\frac{1}{\lbb_{j}^2}\int_{|y\cdot  {\bf v_1}|\geq \frac{3\sigma}{\lbb_{j}}} (|\nabla \wt \ve_j|^2
+|\wt \ve_j|^2)(\Phi_j(\lbb_{j}y+\alpha_{j}) -\phi_{A}(y))dy.
\end{align*}
By the positivity of $\Phi_j$ and
the exponential decay of $\phi_A$,
the second term on the right-hand side above
is bounded from below by
\begin{align*}
  &\frac{\wt c}{\lbb_{j}^2} \int_{|y\cdot  {\bf v_1}|\geq \frac{3\sigma}{\lbb_{j}}}
    (|\nabla \wt \ve_j|^2+|\wt \ve_j|^2)(\Phi_j(\lbb_{j}y+\alpha_{j}) -\phi_{A}(y))dy
 -\frac{1-\wt c}{\lbb_{j}^2} \int_{|y\cdot  {\bf v_1}|\geq \frac{3\sigma}{\lbb_{j}}}
   (|\nabla \wt \ve_j|^2+|\wt \ve_j|^2) \phi_{A}dy\\
\geq&\frac{\wt c}{\lbb_{j}^2}\int_{|y\cdot  {\bf v_1}|\geq \frac{3\sigma}{\lbb_{j}}}
   (|\nabla \wt \ve_j|^2+|\wt \ve_j|^2)(\Phi_j(\lbb_{j}y+\alpha_{j}) -\phi_{A}(y))dy
 -\frac{1-\wt c}{\lbb_{j}^2}  e^{-\frac{3\sigma}{A\lbb_{j}}}  \|\wt \ve_j\|_{H^1}.
\end{align*}
This yields that for $t$ close to $T$,
\begin{align} \label{t2}
 E_{22,j}
 \geq&\frac{\wt c}{\lbb_{j}^2}\int (|\nabla \wt \ve_j|^2+|\wt \ve_j|^2)
    (\Phi_j(\lbb_{j}y+\alpha_{j}) -\phi_{A}(y))dy
    -\frac{1-\wt c}{\lbb_{j}^2} e^{-\frac{3\sigma}{A\lbb_{j}}} \|\wt \ve_j\|_{H^1}  \nonumber \\
\geq& {\wt c} \int (|\nabla R|^2+\frac{1}{\lbb_{j}^2}|R|^2)(\Phi_j-\phi_{A,j})dx
    - C e^{-\frac{\delta}{T-t}} \|R\|_{H^1}^2.
\end{align}
Then,
plugging (\ref{t1}) and (\ref{t2}) into \eqref{R-quad-decoup-2}
we obtain \eqref{R-quad-decoup},
as claimed.

Therefore, taking into account $u(T_*) = S_T(T_*)$ we have
\begin{align*}
E(u(T_*))
=&\frac{1}{2}\int|\sum_{j=1}^K\nabla S_j(T_*)|^2 dx
   - \frac{d}{2d+4} \int |\sum_{j=1}^K S_j(T_*)|^{2+\frac 4d} dx \\
=&\sum_{j=1}^K\frac{\omega_j^2}{8}\|yQ\|_{L^2}^2
  + \calo(e^{-\frac{\delta}{T-t}}),
\end{align*}
and then plugging \eqref{Eu-esti2} and \eqref{R-quad-decoup} into \eqref{Eu-esti1}
we arrive at
\begin{align*}
\sum_{j=1}^{K}\frac{1}{2\lbb_{j}^2}|\beta_{j}(t)|^2 \|Q\|_{L^2}^2
\leq&
\sum_{j=1}^{K}\frac{1}{8\lbb_{j}^2} \|yQ\|_{L^2}^2 |\omega_j^2\lambda_{j}^2(t) - \g_{j}^2(t)| \\
&+|E(u)(t)-E(u)(T_*)|
+ \calo((T-t)^{2 \kappa}),
\end{align*}
which, via Theorem \ref{Thm-energy} and \eqref{lbb-g-t},
yields \eqref{b-g-ve-lbb}
and  finishes the proof.
\hfill $\square$

We are now ready to prove the bootstrap estimates in Theorem \ref{Thm-u-Boot}.

{\bf Proof of Theorem \ref{Thm-u-Boot}}.
{\it (i) Estimate of $R$.}
On one hand,
by \eqref{F-Taylor*},
\begin{align*}
    I
  = & \frac{1}{2}\int |\nabla R|^2
         +\frac{1}{2}\sum_{j=1}^K {\rm Re} \int\frac{1}{\lambda_{j}^2}|R|^2\Phi_j
         - (1+\frac 2d) |U|^{\frac 4d} |R|^2
         - \frac{2}{d}|U|^{\frac 4d -2} U^2 \ol{R}^2 dx \\
  & +\calo( (T-t)^{-2+\frac d2}\sum\limits_{k=3}^{2+\frac 4d} \|R\|_{H^1}^k)
         + \calo(\|R\|_{L^2} \|\na R\|_{L^2}),
\end{align*}
which, via \eqref{lbb-g-t} and \eqref{R-quad-decoup}, yields that for some $0<\wt c <1$,
\begin{align*}
   I \geq& \wt c ( |\nabla R|_{L^2}^2
           +  \frac{1}{(T-t)^2}\|R\|_{L^2}^2 ) \nonumber \\
         &- C(e^{-\frac{\delta}{T-t}} \|R\|_{H^1}^2
             + (T-t)^{-1}\|R\|_{L^2}^2
             + (T-t)^{\frac d2}\|R\|_{H^1}^2
             + \|R\|_{L^2} \|\na R\|_{L^2}).
\end{align*}
Then, taking $T$ small enough
we get that
\begin{align} \label{I-lowbdd}
   I(t) &\geq  \frac {1}{2} \wt c
           (\|\nabla R\|_{L^2}^2
         +  \frac{1}{(T-t)^2}\|R\|_{L^2}^2).
\end{align}

On the other hand,
Theorem \ref{Thm-I-mono} yields that for any $t\in[t_*,T_*]$,
\begin{align}  \label{dIdt-lowbdd}
\frac{dI}{dt}
\geq   - C(A)(T-t)^{2\kappa} - C \ve^*(T-t)^{2\kappa-1}.
\end{align}

Thus, combining \eqref{I-lowbdd} and \eqref{dIdt-lowbdd}
and using the fundamental theorem of calculus
we obtain that for any $t\in[t_*,T_*]$,
\begin{align*}
    \frac{1}{2}\wt c (\|\nabla R(t)\|_{L^2}^2+ \frac{1}{(T-t)^2} \|R(t)\|_{L^2}^2)
 \leq  I(T_*) + \int_{t}^{T_*} C \ve^*(T-r)^{2\kappa-1} + C(A)(T-r)^{2 \kappa} dr.
\end{align*}
Taking into account $I(T_*)=0$ we obtain
\begin{align} \label{naR-R-esti-case1}
\|\nabla R(t)\|_{L^2}^2+\frac{1}{(T-t)^2}\|R(t)\|_{L^2}^2
\leq \frac{C\ve^*}{\kappa  \wt c} (T-t)^{2 \kappa}
     + \frac{2C(A)}{(2 \kappa+1) \wt c} (T-t)^{2 \kappa+1},
\end{align}
which yields \eqref{wn-Tt-boot-2} immediately,
as long as $\ve^*$ and $T$ are sufficiently small
such that
$$\frac{C}{ \kappa \wt c} \ve^* + \frac{2C(A)}{(2 \kappa+1)\wt c} T \leq \frac 18.$$

{\it $(ii)$ Estimates of $\lambda_j$ and $\g_j$.}
Since $(\frac {\gamma_{j}}{ \lambda_{j}})(T_*)=\omega_j$
and by \eqref{Mod-w-lbb},
\begin{align}
  |\frac{d}{dt}(\frac{\gamma_{j}}{\lambda_{j}} )|
=\frac{|\lambda_{j}^2\dot{\gamma}_{j}-\lambda_{j}\dot{\lambda}_{j}\gamma_{j}|}{\lambda_{j}^3}
\leq 2 \frac{Mod}{\lambda_{j}^3}
\leq C (T-t)^{ \kappa},
\end{align}
we infer that for  $T$ small enough such that $C T^{4\zeta} \leq \frac 12$,
\begin{align}  \label{gamlbb-1}
 | (\frac{\gamma_j}{\lambda_j} )(t)-\omega_j(t) |
\leq\int_{t}^{T}  | \frac{d}{dr} (\frac{\gamma_{j}}{\lambda_{j}} ) |dr
\leq C  (T-t)^{ \kappa+1}
\leq \frac{1}{2}(T-t)^{ \kappa+6  \zeta}.
\end{align}
This along with \eqref{Mod-w-lbb} yields that
\begin{align*}
|\frac{d}{dt}(\lambda_{j} -\omega_j (T-t))|
=|\dot{\lambda}_{j}+\frac{\gamma_{j}}{\lambda_{j}}+\omega_j-\frac{\gamma_{j}}{\lambda_{j}}|
\leq\frac{Mod}{\lambda_{j}}+\frac{1}{2} (T-t)^{ \kappa+6  \zeta}
\leq C (T-t)^{ \kappa+6  \zeta},
\end{align*}
which implies that for $T$ possibly  smaller
such that $C T^{ \zeta} \leq \frac{1}{2}$,
\begin{align} \label{lbb-Tt*}
|\lambda_{j} -\omega_j(T-t)|
\leq\int_{t}^{T}  |\frac{d}{dr}(\lambda_{j}-\omega_j(T-r)) |dr
\leq \frac{1}{2} (T-t)^{ \kappa+1+ 5 \zeta},
\end{align}
thereby yielding the estimate of $\lbb_j$ in \eqref{lbbn-Tt-boot-2}.

Similarly,
by \eqref{Mod-w-lbb} and \eqref{gamlbb-1},
\begin{align*}
   |\frac{d}{dt}(\g_j - \omega_j^2 (T-t))|
   =  |\dot\g_{j} + \frac{\g_j^2}{\lbb_j^2} + \omega_j^2 - \frac{\g_j^2}{\lbb_j^2}|
   \leq  \frac{Mod}{\lbb_j^2} + C|\omega_j - \frac{\g_j}{\lbb_j}|
   \leq C (T-t)^{ \kappa+ 6 \zeta},
\end{align*}
which along with $\g_j(T_*) = \omega^2_j (T-T_*)$ yields that
\begin{align*}
   |\g_j(t) - \omega_j^2 (T-t)|
   \leq \int_t^{T_*}  |\frac{d}{dr} (\g_j(r) - \omega_j^2 (T-r)) | dr
   \leq C (T-t)^{ \kappa+1+6 \zeta}.
\end{align*}
Hence, for $T$ very small such that $C T^\zeta \leq \frac 12$ we obtain
\begin{align} \label{g-Tt*}
|\gamma_{j} - \omega_j^2(T-t)|\leq \frac{1}{2} (T-t)^{ \kappa+1+ 5\zeta},
\end{align}
which yields the estimate of $\g_j$ in \eqref{lbbn-Tt-boot-2}.

{\it $(iii)$ Estimates of $\beta_j$ and $\alpha_j$.}
We use the refined estimate of $\beta_j$ in Lemma \ref{Lem-betaj} to get
\begin{align}
|\beta_{j}|^2
\leq C \sum\limits_{j=1}^K
      |\omega_j^2\lambda_{j}^2-\gamma_{j}^2|
       + C(T-t)^{ \kappa+3}
\leq C \sum\limits_{j=1}^K  (\lbb_{j}^2 |\omega_j-\frac {\g_{j}}{ \lbb_{j}}| + \lbb_{j}^{\kappa+3}),
\end{align}
which along with \eqref{lbb-g-t} and \eqref{gamlbb-1}
yields that  for $T$ small enough,
\begin{align} \label{b-Tt*}
|\beta_j|
\leq  C \sum\limits_{j=1}^K  (\lbb_{j} |\omega_j-\frac {\g_{j}}{ \lbb_{j}}|^\frac 12 + \lbb_{j}^{ \frac \kappa 2+\frac 32})
\leq C (T-t)^{\frac{ \kappa}{2}+1+ 3\zeta}
\leq\frac{1}{2}(T-t)^{\frac{ \kappa}{2}+1+2\zeta}.
\end{align}

Moreover, since  $\alpha_j(T_*)=x_j$
and by \eqref{Mod-w-lbb} and \eqref{b-Tt*},
\begin{align}
|\dot{\alpha}_{j}|=
|\frac{\lambda_j\dot{\alpha}_{j}-2\beta_{j}}{\lambda_{j}}+\frac{2\beta_{j}}{\lambda_{j}}|
\leq\frac{Mod}{\lambda_{j}}+\frac{2|\beta_{j}|}{\lambda_{j}}
\leq C (T-t)^{\frac{ \kappa}{2}+ 2\zeta},
\end{align}
we infer that for sufficiently small $T$,
\begin{align}
|\alpha_{j}(t)-x_j|\leq\int_{t}^{T_*}|\dot{\alpha}_{j}(r)|dr
\leq C |T-t|^{\frac{\kappa}{2}+1+2\zeta}
\leq \frac{1}{2} (T-t)^{\frac{ \kappa}{2}+1+\zeta},
\end{align}
which yields \eqref{anbn-Tt-boot-2}.

{\it $(iv)$ Estimate of $\theta_j$.}
By \eqref{lbb-g-t}, \eqref{lbb-Tt*} and \eqref{b-Tt*},
\begin{align}  \label{esti-thetan.0*}
   |\frac{d}{dt}(\theta_{j} - {\omega_j^{-2}(T-t)^{-1}}+ \vartheta_j)|
=& |\frac{\lbb_{j}^2 \dot{\theta}_{j} -1 -|\beta_{j}|^2}{\lbb_{j}^2}
  +\frac{|\beta_{j}|^2}{\lbb_{j}^2}
  + \frac{1}{\lbb_{j}^2}
  - \frac{1}{\omega_j^2(T-t)^2}| \nonumber \\
\leq& \frac{Mod}{\lbb_j^2} + \frac{|\beta_j|^2}{\lbb_j^2} +
     \frac{|\lbb_{j}-\omega_j(T-t)||\lbb_{j}+\omega_j(T-t)|}{\omega_j^2\lbb_{j}^2(T-t)^2 } \nonumber \\
   \leq& C (T-t)^{ \kappa-2+5\zeta},
\end{align}
which yields that for $t$  sufficiently small,
\begin{align}
 |\theta_{j}- ({\omega_j^{-2}(T-t)^{-1}} + \vartheta_j)|
 \leq& \int_{t}^{T} | \frac{d}{dr}(\theta -  {\omega_j^{-2}(T-r)^{-1}} +\vartheta_j) |dr \nonumber \\
 \leq& C (T-t)^{ \kappa-1+5\zeta}
\leq \frac{1}{2} (T-t)^{ \kappa-1+ 4\zeta},
\end{align}
thereby yielding \eqref{thetan-Tt-boot-2}.
The proof of Theorem \ref{Thm-u-Boot} is complete.
\hfill $\square$

\section{Existence of multi-bubble solutions} \label{Sec-Exit}

In this section,
we shall fix $\ve^*>0$ to be sufficiently small,
and for any $0<\ve\leq \ve^*$,
take $\tau^*$ to be very small such that
for a large universal constant $C$,
\begin{align*}
    C (1+\max\limits_{1\leq j\leq K} |x_j|) \tau_*^\frac 14 \leq \frac 12.
\end{align*}
For any $T\in (0,\tau^*)$,
take any increasing sequence $\{t_n\}$ converging to $T$ and
consider the approximating solutions $u_n$
satisfying the equation
\be    \label{equa-u-t}
\left\{ \begin{aligned}
 &i\partial_t u_n+\Delta u_n+|u_n|^{\frac 4d}u_n+ (b \cdot \nabla +c) u_n =0,   \\
 &u_n(t_n)=\sum_{j=1}^{K}S_j(t_n),
\end{aligned}\right.
\ee
where the coefficients $b,c$ are given by \eqref{b} and \eqref{c} respectively,
and  $S_j$
are the pseudo-conformal blow-up solutions defined in \eqref{Sj-blowup},
$1\leq j\leq K$.
We also note that for each $n\geq 1$,
$t_n$ plays the same role as $T_*$ in the previous sections.

We first have the uniform estimates of approximating solutions in Theorem \ref{Thm-u-Unibdd} below.

\begin{theorem}[Uniform estimates]\label{Thm-u-Unibdd}
For $n$ large enough,
$u_n$ admits the unique
geometrical decomposition $u_n=\omega_n+R_n$ on $[0, t_n]$ as in \eqref{u-dec},
and
estimates \eqref{R-Tt}-\eqref{thetan-Tt} hold on $[0,t_n]$.
Moreover, we have
\begin{align}\label{u-sigma}
  \sup\limits_{t\in[0,t_n]} \|xu_n(t)\|_{L^2}\leq C,
\end{align}
and for any $t\in [0,t_n]$,
\begin{align}\label{R-Sigma}
\|R_n(t)\|_{\Sigma}\leq C(T-t)^{\kappa}.
\end{align}
\end{theorem}

{\bf Proof.}
The proof of \eqref{R-Tt}-\eqref{thetan-Tt} is quite similar to that of \cite[Theorem 5.1]{SZ20},
mainly based on the bootstrap estimate in Theorem \ref{Thm-u-Boot}
and the abstract bootstrap principle (see, e.g., \cite[Proposition 1.21]{T06}).
For simplicity,
the details are omitted here.

Below we prove estimates \eqref{u-sigma} and \eqref{R-Sigma}.
Let  $\varphi(x)\in C^1(\R^d,\R)$ be a radial cutoff function such that
$\varphi(x)=0$ for $|x|\leq r$,
and $\varphi(x)=(|x|-r)^2$ for $|x|>r$,
where $r=2\max_{1\leq j\leq K}\{|x_j|,1\}$.
Note that, $|\nabla \varphi|\leq C\varphi^{\half}$.

Using integration by parts formula
we have for some constant $C>0$ independent of $n$,
\begin{align}   \label{1}
 &|\frac{d}{dt}\int |u_n|^2\varphi dx|
 =|{\rm Im}\int (2\ol{u_n} \nabla u_n+b|u_n|^2)\cdot  \nabla\varphi dx| \nonumber \\
\leq& C\int_{ |x-x_j|\geq 1,1\leq j\leq K}(|u_n||\nabla u_n|+|u_n|^2)\varphi^{\half} dx  \nonumber  \\
\leq& C\big((\int_{ |x-x_j|\geq 1,1\leq j\leq K}|\nabla u_n|^2dx)^{\half}+(\int_{ |x-x_j|\geq 1,1\leq j\leq K}|u_n|^2dx)^{\half}\big)
   (\int|u_n|^2\varphi dx)^\half
\end{align}
By \eqref{u-dec} and the exponential decay of the ground state,
\begin{align}
|\int_{ |x-x_j|\geq 1,1\leq j\leq K} |u_n(t)|^2+ |\nabla u_n(t)|^2dx|
\leq C(\|R_n(t)\|_{H^1}^2+e^{-\frac{\delta}{T-t}}).
\end{align}
Thus, taking into account the uniform estimate \eqref{R-Tt}  we get
\begin{align} \label{2}
       |\frac{d}{dt}\int |u_n(t)|^2\varphi dx|
 \leq& C(\|R_n(t)\|_{H^1} +e^{-\frac{\delta}{2(T-t)}}) (\int|u_n(t)|^2\varphi dx)^\frac 12 \nonumber \\
 \leq& C(T-t)^\kappa (\int|u_n(t)|^2\varphi dx)^\half.
\end{align}

Moreover,
using the boundary condition $
u_n(t_n)=\sum_{j=1}^{K}S_j(t_n)$ and $ R_n(t_n)=0$
we  have
\begin{align}
|\int |u_n(t_n)|^2\varphi dx|\leq Ce^{-\frac{\delta}{T-t_n}},
\end{align}

Thus,
integrating from $t$ to $t_n$ we get for $t\in[0,t_n]$,
\begin{align}
\int|u_n(t)|^2\varphi dx\leq C(T-t)^{2\kappa+2}+ C e^{-\frac{\delta}{T-t}} \leq C(T-t)^{2\kappa+2},
\end{align}
which yields that
\begin{align}
   \int|R_n(t)|^2\varphi dx
   \leq& C(\int|U_n(t)|^2\varphi dx+\int|u_n(t)|^2\varphi dx)   \nonumber \\
   \leq& C(T-t)^{2\kappa+2}+ C e^{-\frac{\delta}{T-t}}
   \leq C (T-t)^{2\kappa+2}.
\end{align}
But   $\varphi(x)\geq \half|x|^2$ for $|x|$ large enough.
Hence, we infer that for $t\in[0,t_n]$
\begin{align}
\int|xu_n(t)|^2 dx\leq C(\int|u_n(t)|^2\varphi dx+\int|u_n(t)|^2 dx)\leq C,
\end{align}
and
\begin{align}
\int|xR_n(t)|^2 dx\leq C(\int|R_n(t)|^2\varphi dx+\int|R_n(t)|^2 dx)\leq C(T-t)^{2\kappa+2}.
\end{align}

Therefore,
taking into account \eqref{R-Tt} we obtain \eqref{u-sigma} and \eqref{R-Sigma}
and finish the proof.
\hfill $\square$

{\bf Proof of Theorem \ref{Thm-Multi-Blowup-RNLS}.}
Let $\tau^*, \ve^*$ be as in Theorem \ref{Thm-u-Unibdd}
and let $T\in (0,\tau^*]$, $\ve\in (0,\ve^*]$ be fixed below.
By virtue of Theorem \ref{Thm-u-Unibdd},
we have the geometrical decomposition
\begin{align} \label{dec-un-UnRn}
u_n =U_n +R_n ,\ \ \forall t\in[0,t_n],
\end{align}
with
\begin{align}  \label{STtn}
u_n(t_n) = S_T(t_n):=\sum_{j=1}^{K}S_j(t_n),\  \ R_n(t_n)=0,
\end{align}
where $U_n = \sum_{j=1}^K U_{n,j}$ is as in \eqref{U-Uj}
with the modulation parameters $\calp_{n,j}$ satisfying
\be\ba
    \calp_{n,j}(t_n): &= \lf(\la_{n,j}(t_{n}),\alpha_{n,j}(t_n),\beta_{n,j}(t_n),\gamma_{n,j}(t_n),\theta_{n,j}(t_n)\rt)\\
    &=(\omega_j(T-t_n),x_j,0,\omega_j^2(T-t_n), {\omega_j^{-2}(T-t_n)^{-1}}+\vartheta_j).
\ea\ee
and $S_j$ are the pseudo-conformal blow-up solutions given by \eqref{Sj-blowup}, $1\leq j\leq K$.
Moreover,
the uniform estimates \eqref{R-Tt}-\eqref{thetan-Tt} hold on $[0,t_n]$.

In particular,
$\{u_n(0)\}$ are uniformly bounded in $\Sigma$,
and thus $u_n(0)$ converges weakly to some $u_0 \in \Sigma$.

We claim that $u_n(0)$ indeed converges strongly in $L^2$, i.e.,
\begin{align} \label{unt0-u0-L2}
    u_n(0) \to u_0,\ \  in\ L^2,\ as\ n \to \9.
\end{align}
This follows immediately from the uniform integrability of $\{u_n(0)\}$,
that is, by \eqref{u-sigma},
\begin{align}  \label{un-uninteg-L2}
 \sup\limits_{n\geq 1} \|u_n(0)\|_{L^2(|x|>A)}
  \leq \frac{1}{A} \sup\limits_{n\geq 1}\|x u_n(0)\|_{L^2(|x|>A)}
  \leq \frac{C}{A} \to 0, \ \ as\ A\to \9.
\end{align}

Thus,
by virtue of \eqref{unt0-u0-L2}
and
the $L^2$ local well-posedness theory (see, e.g. \cite{BRZ14})
we obtain a unique $L^2$-solution $u$ to \eqref{equa-u-t} on $[0, T)$
satisfying that $u(0)=u_0$,
and
\begin{align} \label{un-u-0-L2}
\lim_{n\rightarrow \infty}\|u_n-u\|_{C([0,t];L^2)}=0,\ \ t\in [0, T).
\end{align}
Moreover,
since $u_0\in H^1$,
using the $H^1$ local well-posedness result (see, e.g., \cite{BRZ16})
we also have $u\in C([0,t]; H^1)$ for any $0<t<T$.
Such solution is indeed the desirable blow-up solution that explodes at the given $K$ points $\{x_j\}_{j=1}^K$.

As a matter of fact,
let
\begin{align} \label{P0}
   (\lbb_{0,j}, \a_{0,j}, \beta_{0,j}, \g_{0,j}, \theta_{0,j})
           : =  (\omega_j(T-t), x_j, 0, \omega_j^2(T-t), \omega_j^{-2} (T-t)^{-1} +\vartheta_j),
\end{align}
and $Q_{0,j} (t,y):= Q(y)e^{i(\beta_{0,j}(t)\cdot y - \frac{\g_{0,j}(t)}{4}|y|^2)}$.
Since
\begin{align} \label{Unj-Sj}
      U_{n,j} - S_j
   =& (\lbb_{n,j}^{-\frac d2} - \lbb_{0,j}^{-\frac d2}) Q_{n,j}(t, \frac{x-\a_{n,j}}{\lbb_{n,j}}) e^{i\theta_{n,j}}  \nonumber \\
    & + \lbb_{0,j}^{-\frac d2} (Q_{n,j}(t, \frac{x-\a_{n,j}}{\lbb_{n,j}}) - Q_{n,j}(t, \frac{x-\a_{0,j}}{\lbb_{0,j}}))  e^{i\theta_{n,j}} \nonumber \\
    & +  \lbb_{0,j}^{-\frac d2} (Q_{n,j}(t, \frac{x-\a_{0,j}}{\lbb_{0,j}}) - Q_{0,j}(t, \frac{x-\a_{0,j}}{\lbb_{0,j}}))  e^{i\theta_{n,j}} \nonumber \\
    &  + \lbb_{0,j}^{-\frac d2} Q_{0,j}(t, \frac{x-\a_{0,j}}{\lbb_{0,j}}) (e^{i\theta_{n,j}} - e^{i\theta_{0,j}}) ,
\end{align}
using the change of variables
and \eqref{lbb-g-t} we infer that,
if $M:= 1+ \max_{1\leq j\leq K}|x_j|$,
\begin{align*}
   \|x(U_n - S_T)\|_{L^2}
   \leq &  \lbb_{0,j}^{-\frac d2} \||x|(Q_{n,j}(t, \frac{x-\a_{n,j}}{\lbb_{n,j}}) - Q_{n,j}(t, \frac{x-\a_{0,j}}{\lbb_{0,j}}))\|_{L^2} \\
    &  + \lbb_{0,j}^{-\frac d2} \||x|(Q_{n,j}(t, \frac{x-\a_{0,j}}{\lbb_{0,j}}) - Q_{0,j}(t, \frac{x-\a_{0,j}}{\lbb_{0,j}}))\|_{L^2} \\
     & + C M ( |\frac{\lbb_{n,j}^\frac d2 - \lbb_{0,j}^\frac d2}{\lbb_{0,j}^\frac d2}|
                + |\theta_{n,j}-\theta_{0,j}|) \| \<y\> Q\|_{L^2}.
\end{align*}
Then, taking into account
the uniform estimates \eqref{lbbn-Tt}-\eqref{thetan-Tt}
of the modulation parameters and the well localized property of $Q$
we get that for $T$ small enough such that
$MT\leq 1$,
\begin{align*}
   \|x(U_n - S_T)\|_{L^2}
        \leq& C M \sum\limits_{j=1}^K
        \big( | \frac{\lbb_{0,j}}{\lbb_{n,j}} -1|  + |\frac{\a_{n,j}-\a_{0,j}}{\lbb_{n,j}}|
        + |\beta_{n,j} -\beta_{0,j}|   \\
       & \qquad  \ \ +   | \g_{n,j} - \g_{0,j}|
         + |\frac{\lbb_{n,j}^\frac d2 - \lbb_{0,j}^\frac d2}{\lbb_{0,j}^\frac d2}| + |\theta_{n,j} - \theta_{0,j}| \big) \\
   \leq&  CM (T-t)^{\frac \kappa 2 + \zeta}
   \leq  C (T-t)^{\frac \kappa 2 -1 + \zeta},
\end{align*}
where $C$ is a universal constant,
independent of $n, t, M, T$.

Similarly, we have
(see also the proof of \cite[Theorem 2.12]{SZ20})
\begin{align*}
   \|U_n - S_T\|_{L^2}
   \leq& C \sum\limits_{j=1}^K
        \big( | \frac{\lbb_{0,j}}{\lbb_{n,j}} -1|
        + |\frac{\a_{n,j}-\a_{0,j}}{\lbb_{n,j}}|
        + |\beta_{n,j} -\beta_{0,j}|
        +  +    | \g_{n,j} - \g_{0,j}|  \\
       & \qquad  \ \
         + |\frac{\lbb_{n,j}^\frac d2 - \lbb_{0,j}^\frac d2}{\lbb_{0,j}^\frac d2}|
         + |\theta_{n,j} - \theta_{0,j}|  \big)
   \leq C (T-t)^{\frac \kappa 2 + \zeta},
\end{align*}
and
\begin{align*}
    \|\na U_n - \na S_T\|_{L^2}
   \leq& C \sum\limits_{j=1}^K
        \big(  \frac{1}{\lbb_{0,j}} | \frac{\lbb_{0,j}}{\lbb_{n,j}} -1|
             + |\frac{\a_{n,j}-\a_{0,j}}{\lbb_{0,j}\lbb_{n,j}}|
        + | \frac{\beta_{n,j} -\beta_{0,j}}{\lbb_{0,j}} |
         +  |\frac{\g_{n,j} - \g_{0,j}}{\lbb_{0,j}} |   \\
       & \qquad  \ \
         + |\frac{\lbb_{n,j}^{1+\frac d2} - \lbb_{0,j}^{1+\frac d2}}{\lbb_{n,j} \lbb_{0,j}^{1+\frac d2}}|
           + |\frac{\theta_{n,j} - \theta_{0,j}}{\lbb_{0,j}} |   \big)
   \leq C (T-t)^{\frac \kappa 2 -1 + \zeta}.
\end{align*}
Hence,
we conclude that
\begin{align}
   \|U_n(t) - S_T(t)\|_{\Sigma} \leq C (T-t)^{ \frac{\kappa}{2}-1+\zeta}.
\end{align}
This along with \eqref{R-Tt} and \eqref{dec-un-UnRn} yields that
\begin{align}
   \|u_n(t) - S_T(t)\|_{\Sigma}
   \leq \|U_n(t) - S_T(t)\|_{\Sigma} + \|R_n(t)\|_{\Sigma}
   \leq C(T-t)^{\frac{\kappa}{2}-1+\zeta},
\end{align}
where $C$ is independent of $n$.

Hence,
in view of \eqref{un-u-0-L2},
we infer that for some subsequence (still denoted by $\{n\}$),
\begin{align*}
   u_n(t) -S_T(t) \rightharpoonup u(t) - S_T(t),\ \ weakly\ in\ \Sigma,\ as\ n\to \9,
\end{align*}
which yields that
\begin{align*}
    \|u(t) - S_T(t)\|_{\Sigma}
    \leq \liminf\limits_{n\to \9}  \|u_n(t) - S_T(t)\|_{\Sigma}
    \leq C (T-t)^{\frac{\kappa}{2}-1+\zeta}.
\end{align*}

Therefore, the proof of Theorem \ref{Thm-Multi-Blowup-RNLS} is complete.
\hfill $\square$

\section{Uniqueness of multi-bubble solutions}  \label{Sec-Unique}

\subsection{Geometrical decomposition} \label{Subsec-Decomp-Uniq}

In this subsection we obtain the geometrical decomposition and
uniform estimates for the blow-up solution  constructed in the proof of Theorem \ref{Thm-Multi-Blowup-RNLS}
which, actually, are inherited from those of the approximating solutions.

Let us start with the boundedness of the remainders in the more regular space  $H^\frac 32$,
which will be used in Theorem \ref{Thm-wtI-mono} later
to derive the key monotonicity formula of the
generalized functional defined on the difference.

Below we use the same notations $u,u_n, U_n,R_n$ and $t_n$ as in Section \ref{Sec-Exit}.
We also
set $M:= \max_{1\leq j\leq K}|x_j|+1$
and keep using the notation $\kappa := \nu_*-3$.

\begin{proposition} \label{Prop-Rn-H23-bdd}
Assume $(A0)$ and $(A1)$ with $\nu_*\geq 5$.
Then, for $T$ small enough such that
\begin{align} \label{T-M}
    C (1+\max_{1\leq j\leq K}|x_j|) T^\frac 14 \leq \frac 12,
\end{align}
where $C$ is a large universal constant, independent of $\ve, T, n$,
we have
\begin{align} \label{Rn-H23-bdd}
\|R_n(t)\|_{{H}^\frac{3}{2}}\leq (T-t)^{\kappa-2}, \ \  t\in[0,t_n).
\end{align}
\end{proposition}

{\bf Proof.} We rewrite the equation \eqref{equa-R} of $R_n$ as follows
\begin{align}   \label{equa-Rn}
  & i\partial_t R_n +\Delta R_n+b \cdot \nabla R_n+c  R_n
    =  -\eta_n -f(R_n) -(f(u_n)-f(U_n)-f(R_n)),
\end{align}
with $ R_n(t_n)=0$ and $\eta_n$  as in \eqref{etan-Rn},
where $U$ is replaced by $U_n$ given by \eqref{dec-un-UnRn}.
Then, applying the operator $\<\na\>^{\frac 32}$ to both sides of \eqref{equa-Rn}
we obtain
\begin{align}   \label{equa-na32Rn}
  & i\partial_t (\<\na\>^\frac 32 R_n) +\Delta (\<\na\>^\frac 32 R_n) +(b \cdot \nabla +c) (\<\na\>^\frac 32 R_n)  \nonumber \\
   =& [b\cdot \na + c, \<\na\>^{\frac 32}] R_n
     -\<\na\>^\frac 32 \eta_n
    -\<\na\>^\frac 32f(R_n)
     -\<\na\>^\frac 32(f(u_n)-f(U_n)-f(R_n)),
\end{align}
where $[b\cdot \na + c, \<\na\>^{\frac 32}]$ is the commutator
$(b\cdot \na + c)\<\na\>^{\frac 32} - \<\na\>^{\frac 32} (b\cdot \na + c)$.
We regard \eqref{equa-na32Rn} as the equation for the unknown $\<\na\>^\frac 32 R_n$
and apply the Strichartz estimates and local smoothing estimates
(see \cite[Theorem 2.13]{Z17})
to get
\begin{align}\label{Rn-Est-32}
      \|R_n\|_{C([t,t_n]; {H}^\frac{3}{2})}
  \leq& C(\|[b\cdot \na +c,\<\na\>^\frac 32] R_n\|_{L^2(t,t_n; H^{-\frac 12}_1)}
           + \|\<\nabla\>^{\frac{3}{2}}\eta_n\|_{L^{\frac{4+2d}{4+d}}(t,t_n;L^{\frac{4+2d}{4+d}})} \nonumber \\
      &  +\|\<\nabla\>^{\frac{3}{2}}(f(R_n))\|_{L^{\frac{4+2d}{4+d}}(t,t_n;L^{\frac{4+2d}{4+d}})}  \nonumber \\
      &   +\|\<\na\>^{\frac 32}(f(u_n)-f(U_n)-f(R_n))\|_{L^2(t,t_n;H_1^{-\frac 12})})
        =: \sum\limits_{j=1}^4 R_j.
\end{align}
Below we estimate each term $R_j$, $1\leq j \leq 4$, separately.

{\it $(i)$ Estimate of $R_1$.}
Since by Assumption $(A0)$,
$|\pa^\nu b| + |\pa^\nu c| \leq C \<x\>^{-2}$ for any multi-index $\nu$,
using the calculus of pseudo-differential operators (see, e.g., \cite{Z17})
and \eqref{R-Tt} we get
\begin{align} \label{na-bnac-commu}
   R_1 \leq C \|R_n\|_{L^2(t,t_n; H^{1}_{-1})}
   \leq  C (T-t)^\frac 12 \|R_n\|_{C([t,t_n];H^1)}
   \leq C (T-t)^{\kappa+\frac 12}.
\end{align}

{\it $(ii)$ Estimate of $R_2$.}
Using \eqref{equa-Ut} and \eqref{etan-Rn}
we  have
the pointwise estimate of $\eta_n$ that
for any multi-index $\nu$ with $|\nu| \leq 2$
and for $y:=\frac{x-\a_{n,j}}{\lbb_{n,j}}$,
\begin{align*}
  |\partial_x^\nu \eta_n(x)|
  \leq& \sum\limits_{j=1}^K \lbb_{n,j}^{-4-\frac d2}
        |\partial_y^\nu ((\lbb_{n,j} \wt b\cdot \na + \lbb_{n,j}^2 \wt c)Q_{n,j})|
        + \sum\limits_{j=1}^K \lbb_{n,j}^{-4-\frac d2}
         Mod_n \<y\>^2 |\sum\limits_{|\wt \nu|\leq |\nu|+1} \partial_y^{\wt \nu} Q_{n,j}|   \nonumber \\
  =:&\sum\limits_{j=1}^K \eta_{1,j} + \eta_{2,j},
\end{align*}
where $\wt b$ and $\wt c$ are as in \eqref{bc-wtbwtc}.
Note that, for  $p:= \frac{4+2d}{4+d}$, by \eqref{phik-Taylor},
\begin{align*}
   \|\eta_{1,j}\|_{L^p}
   \leq& \sum\limits_{j=1}^K \lbb_{n,j}^{-\frac d2 -4} \lbb_{n,j}^{\frac dp}
         \| \pa_y^\nu (( \lbb_{n,j}\wt b\cdot \na +  \lbb^2_{n,j} \wt c) Q_{n,j}) \|_{L^p}  \\
   \leq& \sum\limits_{j=1}^K \lbb_{n,j}^{d(\frac 1p - \frac 12) -4} \lbb_{n,j}^{\nu_*+1}
   \leq (T-t)^{\kappa+\frac{d}{2+d}}.
\end{align*}
Moreover, by \eqref{Mod-w-lbb},
\begin{align*}
   \| \eta_{2,j}\|_{L^p}
   \leq C(T-t)^{d(\frac 1p- \frac 12)-4}Mod_n
   \leq C(T-t)^{\kappa-\frac{2}{2+d}}.
\end{align*}
Hence,
we obtain that for any multi-index $\nu$ with $|\nu|\leq 2$,
\begin{align*}
   \|\pa_x^\nu \eta_{n}\|_{L^p}
   \leq C(T-t)^{\kappa - \frac{2}{2+d}}.
\end{align*}
This yields that
\begin{align}\label{Rn-Est-32-1}
 R_2
 \leq \| \eta_n\|_{L^{\frac{4+2d}{4+d}}([t,t_n];W^{2,\frac{4+2d}{4+d}})}
\leq C(T-t)^{\frac{4+d}{4+2d} +\kappa- \frac{2}{2+d}}
\leq C(T-t)^{\kappa}.
\end{align}

{\it $(iii)$ Estimate of $R_3$.}
Using the product rule in Lemma \ref{Lem-Pro-rule},
\eqref{fLp-fH1}
and \eqref{R-Tt}
we get
\begin{align}
      \|\<\nabla\>^{\frac{3}{2}}(f(R_n))\|_{L^{\frac{4+2d}{4+d}}}
\leq& C\|\<\nabla\>^{\frac{3}{2}}R_n\|_{L^2}\|R_n\|_{L^{\frac{8+4d}{d}}}^{\frac{4}{d}}   \nonumber \\
\leq& C\|R_n\|_{H^{\frac 32}}\|R_n\|_{H^1}^{\frac{4}{d}}
\leq C(T-t)^{\frac 4d \kappa}  \|R_n\|_{H^{\frac 32}},
\end{align}
which yields that
\begin{align} \label{Rn-Est-32-2}
  R_3 \leq C(T-t)^{\frac{4}{d} \kappa + \frac{4+d}{4+2d}}\|R_n\|_{C([t,t_n]; {H}^\frac{3}{2})}.
\end{align}

{\it $(iv)$ Estimate of $R_4$.}
We first see that
\begin{align}  \label{x-fun-fUn-fRn}
 R_4\leq & C \|\langle x\rangle(f(u_n)-f(U_n)-f(R_n))\|_{L^2(t,t_n;H^{1})} \nonumber \\
\leq&  C \sum\limits_{j=1}^{4/d}
       \big( \|\<x\>|U_n|^{1+\frac 4d- j} |R_n|^j\|_{L^2(t,t_n; L^2)}
       + \|\<x\>|\na U_n| |U_n|^{\frac 4d- j} |R_n|^j\|_{L^2(t,t_n; L^2)}  \nonumber \\
     & \qquad  + \|\<x\>|U_n|^{1+\frac 4d- j} |\na R_n|  |R_n|^{j-1}\|_{L^2(t,t_n; L^2)}
       \big) \nonumber \\
     =:& \sum\limits_{j=1}^{4/d} (R_{4,j1} + R_{4,j2} + R_{4,j3}).
\end{align}
Since
$|\a_j| \leq |x_j-\a_j| +|x_j| \leq M$, $1\leq j\leq K$,
and $\sup_{y\in \bbr^d} \<y\>Q(y) <\9$,
we infer that
\begin{align}  \label{xUn-M}
\|\<x\>U_n(t)\|_{L^\9}
\leq \sum\limits_{j=1}^K \sup_{y\in \bbr^d} \<\lbb_j y + \a_j\> \lbb_j^{-\frac d2} Q(y)
\leq C M (T-t)^{-\frac d2},
\end{align}
which along with \eqref{R-Lp} yields that
\begin{align} \label{H1-esti}
   R_{4,j1}
   \leq C M(T-t)^{-\frac d2(1+\frac 4d -j)+\frac 12} \|R_n\|^j_{C([t,t_n]; L^{2j})}
   \leq C M (T-t)^{\kappa-\frac 12}.
\end{align}

Moreover,
since for any $1<p<\9$,
\begin{align} \label{xUn-M-Lp}
   \|\<x\> |\na U_n| |U_n|^{\frac 4d -j} \|_{L^p}
   \leq C M (T-t)^{\frac d2(j-1) -3 + \frac dp} ,
\end{align}
using H\"older's inequality and \eqref{R-Lp} we get that
\begin{align}  \label{H2-esti}
  R_{4,j2}
   \leq& C \|\<x\>|\na U_n| |U_n|^{\frac 4d -j} \|_{L^2(t,t_n; L^p)}
       \|R_n\|^j_{C([t,t_n]; L^{qj})} \nonumber \\
   \leq& C M  (T-t)^{\frac d2(j-1) - \frac 52 + \frac dp} (T-t)^{j(\kappa+1-\frac d2)+\frac dq}
    \leq C M (T-t)^{\kappa- \frac 32},
\end{align}
where $p,q$ are any positive numbers such that $ \frac 1p + \frac 1q =\frac 12$.

It remains to treat the last term $R_{4,j3}$.
In the case where $j=1$ we have
\begin{align*}
  R_{4,j3} \leq C M (T-t)^{-\frac 32} \|\na R_n\|_{C([t,t_n]; L^2)}
   \leq CM (T-t)^{\kappa-\frac 32}.
\end{align*}
In the case where $d=1$ and $2\leq j\leq \frac 4d$,
using Sobolev's embedding $H^1(\bbr) \hookrightarrow L^\9(\bbr)$ we get
\begin{align*}
 R_{4,j3} \leq C M (T-t)^{-\frac d2(1+\frac 4d -j)+\frac 12} \|R_n\|^j_{C([t,t_n];H^1)}
   \leq C M (T-t)^{2\kappa-1}.
\end{align*}
Furthermore,
in the case where $d=2$ and $2\leq j\leq \frac 4d$,
using the Sobolev embedding $H^{\frac 12, 2}(\bbr^2) \hookrightarrow L^{\frac{2d}{d-1}}(\bbr^2)$ instead
we get
\begin{align*}
  R_{4,j3} \leq& C M (T-t)^{-\frac d2(1+\frac 4d-j)+\frac 12}
     \|R_n^{j-1}\|_{C([t,t_n];L^{2d})}  \| R_n\|_{C([t,t_n];W^{1,\frac{2d}{d-1}})} \\
    \leq& C M (T-t)^{-\frac d2(1+\frac 4d -j)+\frac 12} \|R_n\|^{j-1}_{C([t,t_n];H^1)} \|R_n\|_{C([t,t_n];H^\frac 32)}   \\
    \leq& C M  (T-t)^{\kappa-\frac 32}  \|R_n\|_{C([t,t_n];H^\frac 32)}.
\end{align*}
Hence we conclude that
\begin{align} \label{H4-esti}
    R_{4,j3} \leq C M  ((T-t)^{\kappa-\frac 32} \|R_n\|_{C([t,t_n];H^\frac 32)} + (T-t)^{\kappa-\frac 32}).
\end{align}
Thus, plugging \eqref{H1-esti}-\eqref{H4-esti} into \eqref{x-fun-fUn-fRn}
we obtain
\begin{align} \label{Rn-Est-32-3}
  R_4 \leq  C M ((T-t)^{\kappa-\frac 32} + (T-t)^{\kappa-\frac 32} \|R_n\|_{C([t,t_n]; {H}^\frac{3}{2})}).
\end{align}

Therefore, plugging the estimates
\eqref{na-bnac-commu}, \eqref{Rn-Est-32-1}, \eqref{Rn-Est-32-2} and \eqref{Rn-Est-32-3}
into (\ref{Rn-Est-32})
we obtain
\begin{align*}
   \|R_n\|_{C([t,t_n]; {H}^\frac{3}{2})}
   \leq C M ((T-t)^{\kappa- \frac 32} + (T-t)^{\kappa-\frac 32} \|R_n\|_{C([t,t_n]; {H}^\frac{3}{2})}).
\end{align*}
Hence,
taking $T$ very small such that \eqref{T-M} holds,
we obtain \eqref{Rn-H23-bdd}
and finish the proof.
\hfill $\square$

Proposition \ref{Prop-unST-H23} below shows that,
in the case where $\kappa \geq 3$ (i.e., $\nu_*\geq 6$),
one may enhance the approximation \eqref{u-Sj-H1} of $u$ and $u_n$  in the space $\dot H^\frac 32$.
\begin{proposition} \label{Prop-unST-H23}
Consider the situation as in Proposition \ref{Prop-Rn-H23-bdd}
but with $\nu_*\geq 6$.
Then, we have
\begin{align} \label{un-ST-H23}
    \|u_n(t) - S_T(t) \|_{\dot H^\frac{3}{2}} \leq C (T-t)^{\frac 12 (\kappa-3)+ \zeta},
\end{align}
where $\zeta\in (0,\frac{1}{12})$.
In particular,
for the blow-up solution $u$  constructed in Theorem \ref{Thm-Multi-Blowup-RNLS},
\begin{align} \label{u-ST-H23}
    \|u(t) - S_T(t) \|_{\dot H^\frac{3}{2}} \leq C (T-t)^{\frac 12 (\kappa-3)+\zeta},
\end{align}
and we also have the strong $H^1$ convergence that
for any $t\in (0,T)$,
\begin{align} \label{un-u-0-H1}
    \|u_n - u\|_{C([0,t]; H^1)} \to 0,\ \ as\ n\to \9.
\end{align}
\end{proposition}

{\bf Proof.}
Since
\begin{align*}
   \|u_n(t) - \sum\limits_{j=1}^K S_j(t)\|_{\dot H^\frac{3}{2}}
   \leq \|R_n(t)\|_{\dot H^\frac{3}{2}}
         + \sum\limits_{j=1}^K \|U_{n,j}(t) - S_j(t)\|_{\dot H^\frac{3}{2}},
\end{align*}
in view of Proposition \ref{Prop-Rn-H23-bdd},
we only need to prove that for each $1\leq j\leq K$,
\begin{align} \label{Uj-Sj-t}
    \|U_{n,j}(t) - S_j(t) \|_{\dot H^\frac{3}{2}} \leq C (T-t)^{\frac 12 (\kappa-3)+\zeta}.
\end{align}

For this purpose,
we
let $(\lbb_{0,j}, \a_{0,j}, \beta_{0,j}, \g_{0,j}, \theta_{0,j})$
and
$Q_{0,j}(t,y)$
be as in the proof of  Theorem \ref{Thm-Multi-Blowup-RNLS}.
Then, as in \eqref{Unj-Sj},
we decompose
\begin{align*}
   U_{n,j} - S_j
   =& (\lbb_{n,j}^{-\frac d2} - \lbb_{0,j}^{-\frac d2}) Q_{n,j}(t, \frac{x-\a_{n,j}}{\lbb_{n,j}}) e^{i\theta_{n,j}}
     + \lbb_{0,j}^{-\frac d2} (Q_{n,j}(t, \frac{x-\a_{n,j}}{\lbb_{n,j}}) - Q_{n,j}(t, \frac{x-\a_{0,j}}{\lbb_{0,j}}))  e^{i\theta_{n,j}} \\
    & +  \lbb_{0,j}^{-\frac d2} (Q_{n,j}(t, \frac{x-\a_{0,j}}{\lbb_{0,j}}) - Q_{0,j}(t, \frac{x-\a_{0,j}}{\lbb_{0,j}}))  e^{i\theta_{n,j}}
      + \lbb_{0,j}^{-\frac d2} Q_{0,j}(t, \frac{x-\a_{0,j}}{\lbb_{0,j}}) (e^{i\theta_{n,j}} - e^{i\theta_{0,j}}) \\
   =:& T_1 + T_2 + T_3 +T_4.
\end{align*}
Note that
\begin{align*}
   \|T_1\|_{\dot{H}^\frac 32}
   =  (\lbb_{n,j}^{-\frac d2} - \lbb_{0,j}^{-\frac d2}) \lbb_{n,j}^{\frac d2 -\frac 32}  \|Q\|_{\dot{H}^\frac 32}
   \leq C (T-t)^{ \kappa- \frac 32 + \zeta}.
\end{align*}
Moreover, for $T_2$,
we have
\begin{align*}
   \|T_2\|_{\dot{H}^\frac 32}
   = \lbb_{0,j}^{-\frac d2}
      \| \xi^\frac 32 (\lbb_{n,j}^d e^{-i\a_{n,j} \xi} \wh{Q_{n,j}}(\lbb_{n,j}\xi) - \lbb_{0,j}^d e^{-i\a_{0,j} \xi} \wh{Q_{n,j}}(\lbb_{0,j}\xi))\|_{L^2}
   \leq C(T-t)^{\frac 12(\kappa-3)+ \zeta  }.
\end{align*}
Regarding $T_3$ and $T_4$ we have the bounds
\begin{align*}
   \|T_3\|_{\dot{H}^\frac 32}
   \leq C \lbb_{0,j}^{-\frac 32} (|\beta_{n,j} - \beta_{0,j}| + |\g_{n,j} - \g_{0,j}|)
   \leq C (T-t)^{\frac 12(\kappa-1) + \zeta},
\end{align*}
and
\begin{align*}
   \|T_4\|_{\dot{H}^\frac 32}
   = \lbb_{0,j}^{-\frac 32} |\theta_{n,j} - \theta_{0,j}| \|Q\|_{\dot{H}^\frac 32}
   \leq C (T-t)^{\kappa - \frac 52+\zeta}.
\end{align*}
Thus, putting the estimates above altogether we obtain \eqref{Uj-Sj-t},
and thus \eqref{un-ST-H23} follows.

In view of \eqref{un-u-0-L2} and \eqref{un-ST-H23},
we also infer that
\begin{align} \label{unS-uS-weak}
u_n(t)-S_T(t)\rightharpoonup u(t)-S_T(t),\ \ weakly\ in\ \dot H^\frac{3}{2},\ as\ n\rightarrow\infty,
\end{align}
which yields \eqref{u-ST-H23} immediately.
Moreover, the strong convergence \eqref{un-u-0-H1} in $H^1$
follows from the strong $L^2$ convergence \eqref{un-u-0-L2},
the uniform $H^{\frac 32}$ boundedness \eqref{Rn-H23-bdd} and
standard interpolation arguments.
Therefore, the proof is complete.
\hfill $\square$

Below we show that the constructed blow-up solution $u$
indeed admits the geometrical decomposition as in Proposition \ref{Prop-dec-un}
on the maximal existing time interval $[0,T)$.

For each $0<t<T$ fixed,
Lemma \ref{Lem-Mod-w-lbb}
yields that the derivatives of modulation parameters $\dot \calp_n$ are uniformly bounded
on $[0,t]$,
and thus $\calp_n$ are equicontinuous on $[0,t]$, $n\geq 1$.
Then, by  the Arzel\`{a}-Ascoli Theorem,
$\calp_n$ converges uniformly on $[0,t]$ up to some subsequence
(which may depend on $t$).
But, using the diagonal arguments
one may extract a universal subsequence (still denoted by $\{n\}$)
such that
for some $\calp:= (\calp_1, \cdots, \calp_K)$,
where
$\mathcal{P}_j:=(\lbb_j,\alpha_j,\beta_j,\gamma_j,\theta_j) \in C([0,t];R^{2d+3})$,
$1\leq j\leq K$,
and for every $t\in (0,T)$,
one has
\begin{align} \label{Pn-P}
   \mathcal{P}_n \rightarrow\mathcal{P}\ \ in\ C([0,t];R^{(2d+3)K}).
\end{align}

Then, taking into account the uniform estimates \eqref{R-Tt}-\eqref{thetan-Tt}
we obtain that
for each $t\in[0,T)$ and for $1\leq j\leq K$,
\begin{align}
&\lf|\la_{j}(t) -\omega_j(T-t) \rt| + \ \lf|\gamma_{j}(t)  -\omega_j^2(T-t) \rt|\leq (T-t)^{\kappa+1+\zeta},   \label{lbb-Tt-Uniq} \\
&|\al_{j}(t)-x_j|+|\beta_{j}(t)|\leq (T-t)^{\frac{\kappa}{2}+1+\zeta},  \label{anb-Tt-Uniq}\\
&|\theta_{j}(t) - (\omega_j^{-2}(T-t)^{-1} +\vartheta_j)|\leq (T-t)^{\kappa-1+\zeta}. \label{theta-Tt-Uniq}
\end{align}
In particular,
as in Remark \ref{Rem-lbbj-gj-P-Tt},
$\lbb_j, \g_j, P$ are comparable to $T-t$, i.e.,
there exist $C_1, C_2>0$, independent of $\ve, T$,
such that
\begin{align} \label{lbbj-gj-P-Tt-Uniq}
   C_1(T-t) \leq \lbb_j, \g_j, P \leq C_2(T-t),\ \ \forall 0\leq t<T.
\end{align}

Let
\begin{align} \label{U-Uj-Uniq}
U(t,x) : =\sum_{j=1}^{K}\lbb_{j}^{-\frac d2} Q_{j} (t,\frac{x-\a_{j}}{\lbb_{j}}) e^{i\theta_{j}}\
(:=\sum_{j=1}^{K}U_{j}(t,x)),
\end{align}
with
\begin{align} \label{Qj-Q-Uniq}
   Q_{j}(t,y) := Q(y) e^{i(\beta_{j}(t) \cdot y - \frac 14 \g_{j}(t) {|y|^2})},
\end{align}
define $\varrho_j$ as in \eqref{rhon},
and define the remainder $R$ by
\begin{align} \label{R-Uniq}
   R(t,x):= u(t,x)-U(t,x), \ \ 0\leq t<T, x\in \bbr^d.
\end{align}
Then, for each $0<t<T$,
using the explicit expression \eqref{U-Uj-Uniq}
of $U$
and the convergence \eqref{Pn-P} of mudulation parameters,
we infer that
$\<x\>^2U_n \to  \<x\>^2U$ in $C([0,t];H^1)$,
$U_n \to U$ in $C([0,t];H^\frac 32)$,
$\Lambda U_n \to \Lambda U$ in  $C([0,t];L^2)$,
and $\varrho_{n,j}\to \varrho_j$ in $C([0,t];L^2)$.
Then,
in view of \eqref{un-u-0-L2} and \eqref{unS-uS-weak},
we obtain that $R_n \to R$ in $C([0,t];L^2)$
and $R_n(t) \rightharpoonup R(t)$ weakly in $H^\frac 32$,
which by interpolation yields that
$R_n(t) \to R(t)$ in $H^1$.

Hence, in view of \eqref{R-Tt} and Proposition \ref{Prop-Rn-H23-bdd},
we get that for
$T=T(M)$ small enough satisfying \eqref{T-M}
and for any $t\in [0,T)$,
\begin{align}   \label{R-Tt-Uniq}
\|R(t)\|_{L^2}\leq (T-t)^{\kappa+1},\ \
 \|R(t)\|_{H^1}\leq (T-t)^{\kappa}, \ \
 \|  R(t)\|_{H^{\frac 32}}\leq (T-t)^{\kappa-2}.
\end{align}
Furthermore,
the following orthogonality conditions hold on $[0,T)$  for each $1\leq j\leq K$:
\be\ba\label{ortho-cond-Rn-wn2}
&{\rm Re}\int (x-\a_{j}) U_{j}\ol{R}dx=0,\ \
{\rm Re} \int |x-\a_{j}|^2 U_{j} \ol{R}dx=0,\\
&{\rm Im}\int \nabla U_{j} \ol{R}dx=0,\ \
{\rm Im}\int \Lambda U_{j} \ol{R}dx=0,\ \
{\rm Im}\int \varrho_{j} \ol{R}dx=0.
\ea\ee

We also see that the modulation parameters in $\mathcal{P}$ are $C^1$ functions.
Actually,
similar arguments as in the proof of Proposition \ref{Prop-Mod-bdd}
show that the modulation equations $Mod_n$
can be expressed in terms of the inner products of polynomials of $\pa^\nu U_n$ and $R_n$,
where $|\nu|\leq 2$.
Hence,
the convergence of $\pa^\nu U_n$ and $R_n$
also yields that of $\dot \calp_n$,
which in turn yields the desirable $C^1$-regularity of
the modulation parameter $\calp$.

Thus,
similarly to Lemma \ref{Lem-Mod-w-lbb},
we have
\begin{lemma} \label{Lem-Mod-Uniq}
There exists $C>0$ such that
\begin{align} \label{Mod-bdd-Uniq}
Mod(t)\leq
C(T-t)^{\kappa+3}, \ \   t\in [0,T).
\end{align}
\end{lemma}

\subsection{Energy estimate of the difference}  \label{Subsec-Energy-Diff-Uniq}

Below we assume Assumption $(A1)$ with $\nu_*\geq 11$.
Let $v$ be any blow-up solution to \eqref{equa-u-RNLS}
satisfying \eqref{u-Sj-Tt3-Uniq-RNLS}, i.e.,
\begin{align} \label{v-Sj-t}
\|v(t)-\sum_{j=1}^KS_j(t)\|_{H^1}\leq C(T-t)^{3+\zeta},\ \ 0<t<T.
\end{align}
Set
\begin{align} \label{w-wj}
   w:=v-u =\sum_{j=1}^{K}w_j,\ \ with\ \ w_j :=w\Phi_j,
\end{align}
where the localization functions $\Phi_j$ are given by
\eqref{phi-local}, $1\leq j\leq K$.
Since $\kappa : =\nu_*-3 \geq 8$,
it follows from \eqref{u-Sj-Tt3-Uniq-RNLS} and \eqref{v-Sj-t} that
\begin{align} \label{w-Tt-Uniq}
\|w(t)\|_{H^1}\leq C(T-t)^{3+\zeta}.
\end{align}

For $1\leq j\leq K$,
define the renormalized variable $\epsilon_j$ by
\begin{align} \label{wj-vej-def}
w_j(t,x) :=\lbb_{j}(t)^{-\frac d2} \epsilon_{j} (t,\frac{x-\a_{j}(t)}{\lbb_{j}(t)}) e^{i\theta_{j}(t)}.
\end{align}

Set
\begin{align} \label{D-def}
    D(t):= \ \|\nabla w(t)\|_{L^2}^2+\sum_{j=1}^{K}\frac{\| w_j(t)\|_{L^2}^2}{\lbb_j^2(t)}.
\end{align}
We note that
\begin{align} \label{w-Dt}
    \|w_j(t)\|_{L^2} \leq C(T-t)\sqrt{D(t)},\ \
    \|w_j(t)\|_{H^1} \leq C \sqrt {D(t)}, \ \ 1\leq j\leq K.
\end{align}

Moreover,
as in \eqref{Scal-def},
we set
$Scal_j(t):= Scal(\epsilon_j)$,
i.e.,  for $\epsilon_{j,1} := {\rm Re} \epsilon_j$,
$\epsilon_{j,2} := {\rm Im} \epsilon_j$,
\begin{align} \label{Scaj-def}
Scal_j(t):= \<\epsilon_{j,1},Q\>^2+\<\epsilon_{j,1},yQ\>^2+\<\epsilon_{j,1},|y|^2Q\>^2
           +\<\epsilon_{j,2},\nabla Q\>^2+\<\epsilon_{j,2},\Lambda Q\>^2+\<\epsilon_{j,2},\rho\>^2,
\end{align}
which actually measures the deviations of the remainder term $\epsilon_j$
with respect to the six instable directions of the linearized operator $L$.

By equation \eqref{equa-u-RNLS},
$w$ satisfies the equation
\begin{align}  \label{equa-w}
   & i\partial_t w +\Delta w+f(u+w)-f(u)+b \cdot \nabla w+c  w=0, \ \ t\in (0,T),
\end{align}
and
$\lim_{t\to T} \|w(t)\|_{H^1} =0$, due to \eqref{w-Tt-Uniq}.

The strategy to prove that $w \equiv 0$
is to show that $D \equiv 0$.
As in the proof of the bootstrap estimate
of remainder in Section \ref{Sec-Bootstrap} above,
such result will be derived from the estimate of a generalized energy.
This leads to the definition of $\wt I$ below
\begin{align} \label{def-wtI}
 \wt {I} := &\frac{1}{2}\int |\nabla w|^2dx
            +\frac{1}{2}\sum_{j=1}^K \frac{1}{\lambda_{j}^2} \int  |w|^2 \Phi_jdx
           -{\rm Re}\int F(u+w)-F(u)-f(u)\ol{w}dx \nonumber \\
&+\sum_{j=1}^K\frac{\gamma_{j}}{2\lambda_{j}}{\rm Im} \int (\nabla\chi_A) (\frac{x-\alpha_{j}}{\lambda_{j}})\cdot\nabla w \ol{w}\Phi_jdx.
\end{align}
Note that,
$u$ and $w$ play similar roles as  $U$ and $R$ in \eqref{def-I}, respectively.

Lemma \ref{Lem-wtI-Nt} relates the generalized energy $\wt I(t)$
and the two quantities $D(t)$ and $Scal_j(t)$.

\begin{lemma}    \label{Lem-wtI-Nt}
For $t\in[0,T)$,
there exist $C_1, C_2, C_3>0$ such that
\begin{align}  \label{wtI-Nt}
C_1 D(t)- C_2\sum_{j=1}^{K}\frac{Scal_j(t)}{\lbb_j^2}
\leq  \wt {I}(t)
\leq C_3 D(t).
\end{align}
\end{lemma}

{\bf Proof.}
We first show that,
in the formulation of the generalized energy $\wt I$,
one may replace the blow-up solution $u$ with the main blow-up profile $U$
given by \eqref{U-Uj-Uniq},
at the cost of the error $\calo((T-t)D(t))$, i.e.,
\begin{align} \label{Fuw-Fu-fuw}
  & {\rm Re}\int F(u+w)-F(u)-f(u)\ol{w}dx \nonumber \\
= & {\rm Re}\int F(U+w)-F(U)-f(U)\ol{w}dx
  +\calo((T-t)D(t)).
\end{align}

In order to prove \eqref{Fuw-Fu-fuw},
using Taylor's expansion \eqref{g-gzz-expan}
we see that
\begin{align} \label{Fuw-FUw}
  & {\rm Re}\int F(u+w)-F(u)-f(u)\ol{w}dx  \\
=&{\rm Re}\int F(U+w)-F(U)-f(U)\ol{w}dx
  +\calo(\int (|U|^{\frac 4d-1} + |w|^{\frac 4d-1} + |R|^{\frac 4d-1} ) |R||w|^2dx).  \nonumber
\end{align}
By \eqref{fLp-fH1}, \eqref{U-Tt} and \eqref{R-Tt-Uniq},
\begin{align} \label{U-Rw2-esti}
    \int |U|^{\frac 4d-1} |R||w|^2 dx
    \leq C (T-t)^{-\frac d2(\frac 4d-1)} \|R\|_{L^2} \|w\|^2_{H^1}
    \leq C (T-t)  \|w\|^2_{H^1}.
\end{align}
Moreover,
by \eqref{fLp-fH1}, \eqref{R-Tt-Uniq} and \eqref{w-Tt-Uniq},
\begin{align} \label{w-R-Rw2}
   \int(|w|^{\frac 4d-1} + |R|^{\frac 4d-1} ) |R||w|^2dx
   \leq C (\|R\|_{H^1} \|w\|_{H^1}^{\frac 4d+1}
          + \|R\|^{\frac 4d}_{H^1} \|w\|_{H^1}^{2})
   \leq C(T-t) \|w\|^2_{H^1}.
\end{align}
Hence, plugging \eqref{U-Rw2-esti} and \eqref{w-R-Rw2} into \eqref{Fuw-FUw}
and using \eqref{w-Dt}
we obtain \eqref{Fuw-Fu-fuw}, as claimed.

Next  we analyze the right-hand-side of \eqref{Fuw-Fu-fuw}.
Note that,
by \eqref{g-gzz-expan},
\begin{align} \label{FUw-expan2}
  {\rm Re} (F(U+w) - F(U) - f(U)\ol w)
  =&   \frac 12 (1+\frac 2d) |U|^{\frac 4d} |w|^2
     + \frac{1}{d} |U|^{\frac 4d -2} {\rm Re} (U^2 \ol{w}^2)  \nonumber  \\
   &  + \calo((|U|^{\frac 4d-1} + |w|^{\frac 4d-1})|w|^3).
\end{align}
Note that
\begin{align}  \label{U4d-w3}
   \int (|U|^{\frac 4d-1} + |w|^{\frac 4d-1})|w|^3 dx
    \leq& C(T-t)^{-2} (\|w\|^3_{H^1} + \|w\|_{H^1}^{2+\frac 4d})
   \leq C (T-t) D(t).
\end{align}

Moreover,
a direct application of H\"older's inequality
also shows that
the last term on the right-hand-side of \eqref{def-wtI} is bounded by
\begin{align} \label{wtI-last-bdd}
   C \|w\|_{L^2} \|\na w\|_{L^2}
   \leq C(T-t) D(t).
\end{align}

Thus, we conclude from \eqref{Fuw-Fu-fuw}, \eqref{FUw-expan2}, \eqref{U4d-w3} and \eqref{wtI-last-bdd}  that
\begin{align} \label{wtI-expan2nd}
   \wt{I}
  =&  \frac{1}{2} {\rm Re}
       \int |\nabla w|^2
         + \sum_{j=1}^K  \frac{1}{\lambda_{j}^2}|w|^2\Phi_j
       - (1+\frac 2d) |U|^{\frac 4d} |w|^2
         - \frac{2}{d}|U|^{\frac 4d -2} U^2 \ol{w}^2 dx \nonumber \\
   &   + \calo((T-t)D(t)).
\end{align}
Arguing as in the proof of \eqref{R-quad-decoup}
we have that for some $C_1, C_2>0$,
\begin{align*}
  \wt I
  \geq& C_1 (\int |\na w|^2 + \sum\limits_{j=1}^K \lbb_j^{-2} |w|^2 \Phi_j dx)   \\
      & - C_2( (T-t)D+ (T-t)^{-1}\|w\|_{L^2}^2 + \sum\limits_{j=1}^K \lbb_j^{-2} Scal_j + e^{-\frac{\delta}{T-t}}\|w\|_{H^1}^2).
\end{align*}
Then, taking $T$ small enough
such that
$C_2(2(T-t)+e^{-\frac{\delta}{T-t}}) \leq \frac 12 C_1$, $\forall t\in [0,T]$,
and taking into account \eqref{D-def}
and $\Phi_j \geq \Phi_j^2$, we get
\begin{align*}
  \wt I
  \geq& C_1 (\int |\na w|^2 + \sum\limits_{j=1}^K \lbb_j^{-2} |w|^2 \Phi_j dx)
       - \frac 12 C_1 D(t)
       - C_2  \sum_{j=1}^{K} \lbb_j^{-2}(t) Scal_j(t) \nonumber \\
  \geq& \frac 12 C_1 D(t)
         - C_2 \sum_{j=1}^{K} \lbb_j^{-2}(t) Scal_j(t).
\end{align*}
But, since $|U(t)|\leq C(T-t)^{-\frac d2}$,
using H\"older's inequality,
\eqref{D-def} and \eqref{wtI-expan2nd}
we also have
\begin{align*}
 \wt{I}(t)\leq C D(t).
\end{align*}
Therefore, combining two estimates together
we obtain \eqref{wtI-Nt} and finish the proof.
\hfill $\square$

Similarly to Theorem \ref{Thm-I-mono},
we have the monotonicity property of $\frac{d \wt {I}}{dt}$.

\begin{theorem}  \label{Thm-wtI-mono}
There exist $C_1, C_2,C_3>0$ such that
for any $t\in [0,T)$,
\begin{align}  \label{dIt-mono-new}
\frac{d \wt {I}}{dt}
\geq& \sum_{j=1}^{K}\frac{C_1}{\lambda_{j}}
      \int (|\nabla w_{j}|^2+\frac{1}{\lbb_{j}^2} |w_{j}|^2 )
      e^{-\frac{|x-\alpha_{j}|}{A\lbb_j}}dx   \nonumber \\
    & -C(D(t) +\sum_{j=1}^{K}\frac{Scal_j(t)}{\lbb_j^3(t)})
      - C_3 \ve^* \frac{D(t)}{T-t}.
\end{align}
\end{theorem}

{\bf Proof. }
Similarly to \eqref{equa-I1t} and \eqref{equa-I2t},
using equation \eqref{equa-w} we compute
\begin{align} \label{equa-wtI}
 \frac{d \wt{I}}{dt}
=&-\sum_{j=1}^K\frac{\dot{\lbb}_{j}}{\lbb_{j}^3} {\rm Im}\int|w|^2 \Phi_jdx
-\sum_{j=1}^K\frac{1}{\lbb_{j}^2}{\rm Im} \< f^\prime(u)\cdot w, {w}_j\>
   -{\rm Re} \<f''(u,w)\cdot w^2, \pa_t u\> \nonumber\\
& - \sum_{j=1}^K\frac{1}{\lbb_{j}^2}{\rm Im} \< w\nabla \Phi_j, { \nabla w}\>
 -\sum_{j=1}^K\frac{1}{\lbb_{j}^2} {\rm Im} \< f''(u,w)\cdot w^2, {w}_j\> \nonumber\\
& - \sum\limits_{j=1}^K {\rm Im}
      \<\Delta w - \frac{1}{\lbb_j^2} w_j +f(u+w)-f(u), {b\cdot \nabla w+cw} \>   \nonumber\\
&-\sum_{j=1}^K\frac{\dot{\lambda}_{j}\gamma_{j}-\lambda_{j}\dot{\gamma}_{j}}{2\lambda_{j}^2}
        {\rm Im}  \< \nabla\chi_A(\frac{x-\alpha_{j}}{\lambda_{j}})\cdot\nabla w, {w}_j\>   \nonumber \\
& +\sum_{j=1}^K\frac{\gamma_{j}}{2\lambda_{j}}{\rm Im}
        \< \partial_t (\nabla\chi_A(\frac{x-\alpha_{j}}{\lambda_{j}}))\cdot\nabla w, {w}_j\> \nonumber \\
&+\sum_{j=1}^K {\rm Im}  \< \frac{\gamma_{j}}{2\lambda_{j}^2} \Delta\chi_A(\frac{x-\alpha_{j}}{\lambda_{j}})w_j
+\frac{\gamma_{j}}{2\lambda_{j}}  \nabla\chi_A(\frac{x-\alpha_{j}}{\lambda_{j}})\cdot ( \nabla w_j + \na w \Phi_j),  \pa_t {w} \> \nonumber \\
=:&\sum_{j=1}^{9}\wt{I}_{t,j}.
\end{align}
Below we replace the each appearance of $u$ by the blow-up profile $U$ in
the terms $\wt{I}_{t,2}$, $\wt{I}_{t,3}$, $\wt{I}_{t,5}$, $\wt{I}_{t,6}$ and $\wt{I}_{t,9}$.

{\it $(i)$ Estimate of  $\wt I_{t,2}$}.
Since by \eqref{f'R},
\begin{align} \label{f'u-f'U}
   |f'(u)\cdot w - f'(U)\cdot w| \leq C(|U|^{\frac 4d-1} + |R|^{\frac 4d-1}  + |w|^{\frac 4d-1}) |R||w|,
\end{align}
using Lemma \ref{Lem-inter-est} we get
\begin{align*}
|\wt{I}_{t,2}+\sum_{j=1}^K\frac{1}{\lbb_{j}^2}{\rm Im}
       \< f^\prime(U)\cdot w, {w}_j\>|
\leq& C \sum_{j=1}^K\frac{1}{\lbb_{j}^2}
        \int (|U|^{\frac{4}{d}-1}+ |R|^{\frac 4d-1} +|w|^{\frac 4d-1})|R||w|| {w}_j|dx.
\end{align*}
By H\"older's inequality, \eqref{G-N},
\eqref{U-Tt}, \eqref{R-Lp}  and \eqref{R-Tt-Uniq},
we have that for $\frac {1}{p} + \frac {1}{q} =\frac 12$,
\begin{align*}
    \frac{1}{\lbb_{j}^2} \int  |U|^{\frac{4}{d}-1} |R||w|| {w_j}|dx
    \leq  (T-t)^{-4+\frac d2} \|R\|_{L^{p}}  \|w\|_{L^{q}} \|w_j\|_{L^2}
    \leq  (T-t)^{\kappa-3+\frac{d}{p}} \|w\|_{L^2}^{2-\frac{d}{p}} \|\na w\|_{L^2}^{\frac{d}{p}}.
\end{align*}
Then, by Young's inequality
$ab\leq \frac{a^{\wt p}}{\wt p} + \frac{b^{\wt p'}}{\wt p'}$
with $\wt p'= \frac{2p}{d}$, $\frac{1}{\wt p} + \frac{1}{\wt p'} =1$ and $0<\frac{d}{p}<1$,
\begin{align*}
   \frac{1}{\lbb_{j}^2} \int  |U|^{\frac{4}{d}-1} |R||w|| {w_j}|dx
   \leq (T-t)^{(\kappa-3)\wt p+d} \|w_j\|_{L^2}^2 + \|\na w_j\|_{L^2}^2
   \leq C D(t).
\end{align*}
Moreover, by \eqref{fLp-fH1} and \eqref{R-Tt-Uniq},
\begin{align*}
   \frac{1}{\lbb_j^2}\int |R|^{\frac 4d} |w_j| |w| dx
   \leq C(T-t)^{-2}\|R\|_{H^1}^{\frac 4d} \|w\|_{H^1}^2
   \leq C \|w\|_{H^1}^2 \leq CD(t),
\end{align*}
and
\begin{align*}
   \frac{1}{\lbb_j^2}\int |w|^{\frac 4d} |R||w_j| dx
   \leq C (T-t)^{-2} \|R\|_{L^2} \|w\|_{H^1}^{1+\frac 4d}
   \leq C \|w\|_{H^1}^2 \leq C D(t).
\end{align*}
Hence,
we conclude that
\begin{align} \label{wtIt3-U-Nt}
    \wt I_{t,2}
    = - \sum\limits_{j=1}^K \frac{1}{\lbb_j^2}
       {\rm Im} \< f'(U)\cdot w,  w_j\>
     + \calo(D(t)).
\end{align}

{\it $(ii)$ Estimate of $\wt{I}_{t,3}$}.
Since by \eqref{R-Uniq},
\begin{align*}
   {\rm Re} \< f''(u,w)\cdot w^2,  \partial_t {u} \>
   = {\rm Re} \< f''(u,w)\cdot w^2,  \partial_t {U} \>
     + {\rm Re} \< f''(u,w)\cdot w^2,  \partial_t {R} \>,
\end{align*}
we shall treat the two terms on the right-hand side separately below.

First, using \eqref{f''R2} we see that
\begin{align} \label{f''uw-f''Uw}
|f''(u,w)\cdot w^2 - f''(U,w)\cdot w^2| \leq C(|U|^{\frac 4d-2} + |R|^{\frac 4d-2}  + |w|^{\frac 4d-2}) |R||w|^2,
\end{align}
which yields that
\begin{align}
&|{\rm Re} \< f''(u,w)\cdot w^2,  \partial_t {U} \>
   -{\rm Re} \< f''(U,w)\cdot w^2, \partial_t{U}\>| \nonumber \\
\leq& C\sum_{j=1}^K\|\partial_tU_j\|_{L^\infty}
    \int (|U|^{\frac{4}{d}-2}+ |R|^{\frac{4}{d}-2}+ |w|^{\frac{4}{d}-2}) |R||w|^2dx.
\end{align}
Note that,
by \eqref{equa-Ut},
$\|\partial_t U_j\|_{L^\9} \leq C (T-t)^{-2-\frac d2}$,
and by Gagliardo-Nirenberg's inequality \eqref{G-N},
\begin{align*}
   \|\pa_t U_j\|_{L^\9}  \int |U|^{\frac{4}{d}-2} |R||w|^2 dx
   \leq C(T-t)^{-4+\frac d2} \|R\|_{L^2} \|w\|_{L^4}^2
   \leq C(T-t)^{\kappa+ \frac d2-3} \|w\|_{L^2}^{2-\frac d2} \|\na w\|_{L^2}^\frac d2,
\end{align*}
which,
via Young's inequality $ab\leq \frac{a^p}{p} + \frac{b^q}{q}$
with $p=\frac{4}{4-d}$ and $q=\frac 4d$,
can be bounded by
\begin{align*}
   C(T-t)^{(\kappa+\frac d2 - 3) \frac{4}{d-2}}
   \|w\|_{L^2}^2 + \|\na w\|_{L^2}^2
   \leq C(T-t)^{-2} \|w\|_{L^2}^2 + \|\na w\|_{L^2}^2
   \leq CD(t).
\end{align*}
Moreover, by H\"older's inequality,
\eqref{fghLp-H1}, \eqref{R-Lp},
\eqref{R-Tt-Uniq} and \eqref{w-Tt-Uniq},
\begin{align*}
   \|\partial_tU_j\|_{L^\infty}
    \int  ( |R|^{\frac{4}{d}-2}+ |w|^{\frac{4}{d}-2}) |R||w|^2dx
    \leq& C(T-t)^{-2-\frac d2} (\|R\|_{L^{\frac 8d-2}}^{\frac 4d-1} \|w\|_{H^1}^2 + \|R\|_{L^2}\|w\|_{H^1}^{\frac 4d})  \\
    \leq& C\|w\|^2_{H^1}
    \leq C D(t).
\end{align*}
Hence, we obtain
\begin{align} \label{f''uww2-ptU}
     {\rm Re} \< f''(u,w)\cdot w^2,  \partial_t {U} \>
   ={\rm Re} \< f''(U,w)\cdot w^2, \partial_t{U}\>
     + \calo(D(t)).
\end{align}

Next we show that
\begin{align} \label{f''uww2-ptR}
     {\rm Re} \< f''(u,w)\cdot w^2,  \partial_t {R} \>
   = \calo(D(t)).
\end{align}
For this purpose,
using equation (\ref{equa-R}) we get
\begin{align*}
    |{\rm Re} \< f''(u,w)\cdot w^2,   \partial_t {R} \> |
= |{\rm Im} \< f''(u,w)\cdot w^2 , \Delta R+f(u)-f(U)+b\cdot \nabla R+cR+\eta \>|.
\end{align*}

Note that, by \eqref{R-Tt-Uniq},
\begin{align} \label{f''uw2DR}
  |{\rm Im} \< f''(u,w)\cdot w^2, \Delta R \>|
 \leq& C\|R\|_{\dot{H}^{\frac{3}{2}}}\|f''(u,w)\cdot w^2\|_{\dot{H}^{\frac{1}{2}}}  \nonumber \\
 \leq& C (T-t)^{\kappa-2} \|f''(u,w)\cdot w^2\|_{\dot{H}^{\frac{1}{2}}} .
\end{align}
Using  \eqref{fghH12-H1}, \eqref{R-Tt-Uniq} and $\|U(t)\|_{H^1} \leq C(T-t)^{-1}$
we get
\begin{align} \label{f''uw2-H12}
   \|f''(u,w)\cdot w^2  \|_{\dot{H}^{\frac{1}{2}}}
   \leq& C \sum\limits_{j=2}^{1+\frac 4d}
          \|u\|_{H^1}^{1+\frac 4d -j} \|w\|_{H^1}^j
   \leq C \sum\limits_{j=2}^{1+\frac 4d}
          (\|U\|_{H^1}^{1+\frac 4d -j} + \|R\|_{H^1}^{1+\frac 4d -j}) \|w\|_{H^1}^j \nonumber \\
   \leq&  C  \sum\limits_{k=2}^{1+\frac 4d} ((T-t)^{-(1+\frac 4d -j)} + (T-t)^{\kappa(1+\frac 4d -j)}) (T-t)^{(3+\zeta)(j-2)} \|w\|_{H^1}^2  \nonumber   \\
   \leq& C (T-t)^{1-\frac 4d} \|w\|_{H^1}^2,
\end{align}
which along with \eqref{f''uw2DR}
and $\kappa \geq 8$ yields  that
\begin{align} \label{f''uww2-DR}
     |{\rm Im} \< f''(u,w)\cdot w^2,  \Delta {R} \>|
    \leq C \|w\|_{H^1}^2
   \leq C D(t).
\end{align}

Moreover,
since
\begin{align} \label{fu-fU-bdd}
   |f(u)-f(U)|\leq C(|U|^{\frac 4d}+|R|^{\frac 4d})|R|,
\end{align}
and
\begin{align} \label{f''uw-w2-bdd}
   |f''(u,w)\cdot w^2|
   \leq C (|U|^{\frac 4d-1} + |R|^{\frac 4d -1} + |w|^{\frac 4d -1}) |w|^2,
\end{align}
by \eqref{fghLp-H1} and \eqref{R-Tt-Uniq} we get
\begin{align} \label{f''uww2-fufUbcR}
  |{\rm Im} \< f''(u,w)\cdot w^2,  f(u)-f(U)+b\cdot \nabla R+cR \>|
  \leq C \|w\|_{H^1}^2 \leq C D(t).
\end{align}

Furthermore,
using \eqref{fghLp-H1},
\eqref{eta-L2} and \eqref{f''uw-w2-bdd} again
we also have
\begin{align} \label{f''uww2-esta}
  |{\rm Im} \< f''(u,w)\cdot w^2, \eta \>|
  \leq& C((T-t)^{-2+\frac d2} + \|R\|_{H^1}^{\frac 4d-1} + \|w\|_{H^1}^{\frac 4d -1}) \|\eta\|_{L^2} \|w\|_{H^1}^2  \nonumber \\
  \leq& C \|w\|_{H^1}
   \leq C D(t).
\end{align}

Thus, combining  estimates \eqref{f''uww2-DR}, \eqref{f''uww2-fufUbcR} and \eqref{f''uww2-esta}
we obtain \eqref{f''uww2-ptR}, as claimed.

Therefore, we infer from \eqref{f''uww2-ptU} and \eqref{f''uww2-esta} that
\begin{align}
   \wt I_{t,3}
    = {\rm Re} \<f''(U,w)\cdot w^2, \pa_t U\> + \calo(D(t)).
\end{align}

{\it $(iii)$ Estimate of $\wt I_{t,5}$}.
Using \eqref{f''uw-f''Uw}
we have
\begin{align*}
|\wt{I}_{t,5}+\sum_{j=1}^K\frac{1}{\lbb_{j}^2}{\rm Im} \< f''(U,w)&\cdot w^2, {w}_j\>|
\leq C (T-t)^{-2} \int(|U|^{\frac 4d-2} + |R|^{\frac 4d-2} +  |w|^{\frac 4d-2}) |R||w|^2 |w_j| dx.
\end{align*}
By H\"older's inequality, \eqref{fLp-fH1}, \eqref{U-Tt}, \eqref{R-Tt-Uniq} and \eqref{w-Tt-Uniq},
\begin{align*}
   (T-t)^{-2} \int  |U|^{\frac{4}{d}-2}   |R||w|^2|w_j|dx
   \leq& (T-t)^{-2-\frac d2(\frac 4d-2)} \|R\|_{L^2} \|w\|^3_{H^1}  \nonumber \\
   \leq& C (T-t)^{d+\kappa+\zeta} \|w\|_{H^1}^2
   \leq C D(t).
\end{align*}
Similarly,
using \eqref{fLp-fH1}, \eqref{R-Tt-Uniq} and \eqref{w-Tt-Uniq}
we also have
\begin{align*}
  (T-t)^{-2} \int  (|R|^{\frac{4}{d}-2} +|w|^{\frac{4}{d}-2})  |R||w|^2|w_j|dx
  \leq& C (T-t)^{-2} (\|R\|_{H^1}^{\frac 4d-1} \|w\|_{H^1}^3 +\|R\|_{L^2} \|w\|_{H^1}^{\frac 4d+1})  \\
  \leq& C D(t).
\end{align*}
Hence, we obtain
\begin{align} \label{wtI5-U-N}
   \wt I_{t,5}
   = - \sum\limits_{j=1}^N \frac{1}{\lbb_j^2}
     {\rm Im}\< f''(U,w)\cdot w^2, w_j\>
     + \calo(D(t)).
\end{align}

{\it $(iv)$ Estimate of $\wt I_{t,6}$}.
Since by \eqref{g-gz-expan},
\begin{align} \label{fuwfu-fUwfU}
   |f(u+w)-f(u) -(f(U+w)-f(U))|
\leq C (|U|^{\frac 4d-1} + |w|^{\frac 4d-1} + |R|^{\frac 4d-1})|R||w|,
\end{align}
taking into account \eqref{fLp-fH1} we infer that
\begin{align*}
   & {\rm Im} \< f(u+w)-f(u) -(f(U+w)-f(U)),{(b\cdot \na +c)w} \>  \\
   \leq& C \int (|U|^{\frac 4d-1} + |w|^{\frac 4d-1} + |R|^{\frac 4d-1}) |R||w| |(b\cdot \na +c)w| dx \\
   \leq& C [(T-t)^{-\frac d2(\frac 4d-1)} + \|w\|_{H^1}^{\frac 4d-1} + \|R\|_{H^1}^{\frac 4d -1}] \|R\|_{H^1} \|w\|_{H^1} \|(b\cdot \na+c)w\|_{L^2} \\
   \leq& C \|w\|_{H^1}^2
   \leq C D(t).
\end{align*}
This yields that
\begin{align} \label{wtIt6-U-Nt}
   \wt I_{t,6}
   = - \sum\limits_{j=1}^K {\rm Im}
     \<\Delta w - \frac{1}{\lbb_j^2} w_j +f(U+w)-f(U), {b\cdot \nabla w+cw}  \>
    + \calo(D(t)).
\end{align}

{\it $(v)$ Estimate of $\wt I_{t,9}$}.
The arguments are similar to those in
the previous cases $(i)$ and $(iv)$.
Actually,
using equation \eqref{equa-w} we infer that
\begin{align*}
  \wt I_{t,9}
  =& {\rm Im} \big\< \frac{\gamma_{j}}{2\lambda_{j}^2} \Delta\chi_A(\frac{x-\alpha_{j}}{\lambda_{j}})w_j
+\frac{\gamma_{j}}{2\lambda_{j}}  \nabla\chi_A(\frac{x-\alpha_{j}}{\lambda_{j}})\cdot ( \nabla w_j + \na w \Phi_j),  \\
   &\qquad  i\Delta w + i(f(u+w)-f(u)) + i (b\cdot \na +c)w \big\> .
\end{align*}
Then, in view of \eqref{fuwfu-fUwfU} we see that
\begin{align*}
  & |\< \frac{\gamma_{j}}{2\lambda_{j}^2} \Delta\chi_A(\frac{x-\alpha_{j}}{\lambda_{j}})w_j
      +\frac{\gamma_{j}}{2\lambda_{j}}  \nabla\chi_A(\frac{x-\alpha_{j}}{\lambda_{j}})\cdot ( \nabla w_j + \na w \Phi_j),
      i(f(u+w)-f(u)) \> \\
  & \ - \< \frac{\gamma_{j}}{2\lambda_{j}^2} \Delta\chi_A(\frac{x-\alpha_{j}}{\lambda_{j}})w_j
         +\frac{\gamma_{j}}{2\lambda_{j}}  \nabla\chi_A(\frac{x-\alpha_{j}}{\lambda_{j}})\cdot ( \nabla w_j + \na w \Phi_j),
         i(f(U+w)-f(U)) \>| \\
  \leq&C(A) \int (T-t)^{-1} |w_j| (|U|^{\frac 4d-1} + |w|^{\frac 4d-1} + |R|^{\frac 4d-1}) |R||w| dx \\
      & +  C(A)  \int |\na w_j + \na w \Phi_j| (|U|^{\frac 4d-1} + |w|^{\frac 4d-1} + |R|^{\frac 4d-1}) |R||w| dx.
\end{align*}
Arguing as in the proof of \eqref{wtIt3-U-Nt}  and \eqref{wtIt6-U-Nt}, respectively,
we can also bound the two integrations on the right-hand side above by $C D(t)$.
This yields that
\begin{align} \label{wtIt9-U-Nt}
   \wt I_{t,9}
   =& {\rm Im} \big\< \frac{\gamma_{j}}{2\lambda_{j}^2} \Delta\chi_A(\frac{x-\alpha_{j}}{\lambda_{j}})w_j
+\frac{\gamma_{j}}{2\lambda_{j}}  \nabla\chi_A(\frac{x-\alpha_{j}}{\lambda_{j}})\cdot ( \nabla w_j + \na w \Phi_j), \nonumber  \\
   &\qquad  i\Delta w + i(f(U+w)-f(U)) + i (b\cdot \na +c)w \big\>
     + \mathcal{O}(D(t)).
\end{align}

At this stage,
we have replaced $u$ by the blow-up profile $U$ in \eqref{equa-wtI}.
Then, arguing as in the proof of Theorem \ref{Thm-I-mono}
with $R$ replaced by $\omega$
and using \eqref{w-Tt-Uniq}
we obtain \eqref{dIt-mono-new}
and thus finish the proof.
\hfill $\square$

As a consequence of Lemma \ref{Lem-wtI-Nt} and  Theorem \ref{Thm-wtI-mono}
we have

\begin{theorem} \label{Thm-Nt-Scal}
For $t\in[0,T)$, we have
that
\begin{align} \label{Nt-Scal}
    \sup\limits_{t\leq s < T} D(s)
\leq C (\sum_{j=1}^{K} \sup\limits_{t\leq s < T} \frac{Scal_j(s)}{\lbb_j^2(s)}
+ \int_{t}^{T}\sum_{j=1}^{K}\frac{Scal_j(s)}{\lbb_j^3(s)}  + \ve^* \frac{D(s)}{T-s}ds).
\end{align}
\end{theorem}

{\bf Proof.}
By Lemma \ref{Lem-wtI-Nt} and Theorem \ref{Thm-wtI-mono},
for $t<\wt t <T$,
\begin{align*}
     C_1D(t)
&\leq  \wt {I}(t) +C_2\sum_{j=1}^{K}\frac{Scal_j(t)}{\lbb_j^2(t)}
    =  \wt I(\wt{t})   +C_2\sum_{j=1}^{K}\frac{Scal_j(t)}{\lbb_j^2(t)}
-\int_{t}^{\wt{t}}\frac{d \wt{I}}{ds}(s) ds \nonumber \\
&\leq  C( D(\wt t) + \sum_{j=1}^{K}\frac{Scal_j(t)}{\lbb_j^2(t)}
       + \int_{t}^{\wt {t}} D(s) + \sum_{j=1}^{K}\frac{Scal_j(s)}{\lbb_j^3(s)}
       + \ve^* \frac{D(s)}{T-s}ds),
\end{align*}
which yields that
\begin{align*}
   \sup\limits_{t\leq s\leq \wt t}D(s)  \leq C (D(\wt t)
             + \sum_{j=1}^{K} \sup\limits_{t\leq s\leq \wt t} \frac{Scal_j(s)}{\lbb_j^2(s)}
              + (\wt t- t)  \sup\limits_{t\leq s\leq \wt t}D(s)
         +\int_t^{\wt t} \sum\limits_{j=1}^K \frac{Scal_j(s)}{\lbb_j^3(s)}
         + \ve^* \frac{D(s)}{T-s}ds).
\end{align*}
By \eqref{w-Tt-Uniq},
$D(\wt t) \to 0$ as $\wt t \to T$.
Hence,
letting $\wt t \to T$ and
taking $T$ sufficiently small
we obtain \eqref{Nt-Scal}
and finish the proof.
\hfill $\square$

\subsection{Control of the null space}

In view of Theorem \ref{Thm-Nt-Scal},
the last step is to control the scalar products in $Scal_j$,
that is,  the
growth along six unstable directions
in the null space.
The key result is formulated in Theorem \ref{Thm-Scal} below.
Then,
at the end of this subsection,
we  finish the proof of the main uniqueness result.

\begin{theorem}\label{Thm-Scal}
For $T$ small enough
and for $1\leq j\leq K$,
there exist $C>0, \zeta\in (0,1)$
such that
\begin{align} \label{Scalj-N}
Scal_j(t)\leq C(T-t)^{2+\zeta} \sup\limits_{t\leq s<T}D(s).
\end{align}
\end{theorem}

To begin with,
we first treat the estimate of the scalar product ${\rm Re}\<U_j, w\>$.

\begin{proposition} \label{Prop-U-w}
For $1\leq j\leq K$, we have
\begin{align}  \label{x-Uj-w}
|{\rm Re}\int  U_{j}(t)\ol{w}(t)dx|
   \leq C(T-t)^{4+\zeta} \sup\limits_{t\leq s<T}\sqrt{D(s)}.
\end{align}
\end{proposition}

\begin{remark}
One may also use equation \eqref{equa-ej} below to obtain the bound
\begin{align} \label{Q-ej}
  |{\rm Re}\int  U_j(t) \ol{w}(t)dx|
 =|{\rm Re}\int  Q e_j(t)dy|
     \leq C (T-t)^{3+\zeta} \sup\limits_{t\leq s<T}\sqrt{D(s)},
\end{align}
where $e_j$ is defined in \eqref{w-wtej-ej-def} below.
Estimate \eqref{x-Uj-w} improves \eqref{Q-ej}
by a factor $(T-t)$,
by exploring the conservation law of mass.
\end{remark}

{\bf Proof of Proposition \ref{Prop-U-w}.}
We first note from \eqref{w-wj}, $v=u+w$, that
\begin{align*}
\int |v(t)|^2 \Phi_jdx
= \int |u(t)|^2 \Phi_jdx +\int |w(t)|^2 \Phi_jdx +2{\rm Re} \< w_j, {u}\>.
\end{align*}
Taking into account $u= U+R$,
Lemma \ref{Lem-inter-est}, \eqref{R-Tt-Uniq} and \eqref{w-Tt-Uniq}
we get,
for $T$ small enough,
\begin{align*}
{\rm Re} \< w_j, {u} \>
=& {\rm Re} \< w, {U}_j\> +{\rm Re} \< w_j, {R}\>
  + \calo(e^{-\frac {\delta}{T- t}} \|w\|_{L^2}) \\
=& {\rm Re} \<U_j, w\> + \calo((T-t)^{\kappa +1}\|w\|_{L^2}).
\end{align*}
Hence,
taking into account \eqref{w-Tt-Uniq}
and $\kappa \geq 3$
we obtain
\begin{align} \label{Uj-w-esti}
{\rm Re} \langle U_{j}(t),w(t) \rangle
  =& \frac{1}{2}(\int |v(t)|^2\Phi_jdx
      -\int |u(t)|^2\Phi_jdx) \nonumber \\
   & -\half\int |w(t)|^2 \Phi_jdx
    + \calo((T-t)^{\kappa+1} \|w(t)\|_{L^2}) \nonumber  \\
  =& \frac{1}{2}(\int |v(t)|^2\Phi_jdx
      -\int |u(t)|^2 \Phi_jdx)
      + \calo((T-t)^{3+\zeta}\|w(t)\|_{L^2}).
\end{align}

In order to estimate the first term on the right-hand side above,
we note that for $\wt t \in (t,T)$,
\begin{align} \label{v-u-L2-diff}
  \int |v(t)|^2 \Phi_jdx- \int |u(t)|^2 \Phi_jdx
=&\int_{\tilde{t}}^{t}(\frac{d}{ds}\int |v|^2 \Phi_j dx-\frac{d}{ds}\int |u|^2\Phi_jdx)ds  \nonumber \\
&+(\int |v(\tilde{t})|^2 \Phi_jdx- \int |u(\tilde{t})|^2 \Phi_jdx).
\end{align}
Similarly to \eqref{du2-bc}, we have
\begin{align*}
   \frac{d}{dt} \int |v|^2 \Phi_j dx
   = {\rm Im} \int (2\ol v \na v + b|v|^2)\cdot \na \Phi_j dx.
\end{align*}
Similar equation also holds for $u$.
Then, taking into account $v=u+w$ and $u=U+R$
we get
\begin{align}
&\frac{d}{dt}\int |v|^2 \Phi_j dx-\frac{d}{dt}\int |u|^2 \Phi_jdx  \nonumber \\
=&{\rm Im}\int (2(\ol{w}\nabla u+\ol{u}\nabla w)+2 \ol{w}\nabla w+b( \ol{w} u+ \ol{u} w+|w|^2))\cdot  \nabla\Phi_j dx  \nonumber \\
=&{\rm Im} \int_{ |x-x_l|\geq 4\sigma,1\leq l\leq K}(2( \ol{w}\nabla U+ \ol{U}\nabla w)+b( \ol{w} U+ \ol{U} w))\cdot  \nabla\Phi_j dx  \nonumber \\
&+{\rm Im}\int (2(\ol{w}\nabla R+ \ol{R}\nabla w)+2 \ol{w}\nabla w+b(\ol{w} R+\ol{R} w+|w|^2))\cdot  \nabla\Phi_j dx,
\end{align}
which along with Lemma \ref{Lem-inter-est},
integration by parts formula,
H\"older's inequality and \eqref{R-Tt-Uniq}
yields that
\begin{align}
|\frac{d}{dt}\int |v |^2 \Phi_jdx-\frac{d}{dt}\int |u |^2\Phi_jdx|
\leq& C(\|R\|_{H^1}+\| w\|_{H^1}+e^{-\frac {\delta}{T- t}}) \|w\|_{L^2}  \nonumber \\
\leq& C(T-t)^{3+\zeta} \|w\|_{L^2},
\end{align}
and thus
\begin{align}  \label{uni-4}
|\int_t^{\tilde{t}}(\frac{d}{ds}\int |v|^2 \Phi_jdx-\frac{d}{ds}\int |u|^2 \Phi_jdx)ds|
\leq C(T-t)^{4+\zeta} \sup\limits_{t\leq s\leq \wt t} \|w(s)\|_{L^2}.
\end{align}

Moreover,
taking into account \eqref{v-Sj-t} we   have
\begin{align} \label{v-u-L2-0}
\lim_{\wt t\rightarrow T}|\int |v(\wt t)|^2 \Phi_jdx- \int |u(\wt t)|^2 \Phi_jdx|=0.
\end{align}

Therefore,
plugging \eqref{uni-4} and \eqref{v-u-L2-0}  into \eqref{v-u-L2-diff}
and  passing to the limit $\tilde{t}\rightarrow T$
we arrive at
\begin{align}  \label{uni-5}
|\int |v(t)|^2 \Phi_jdx- \int |u(t)|^2 \Phi_jdx|
  \leq C(T-t)^{4+\zeta}  \sup\limits_{t\leq s <T} \|w(s)\|_{L^2},
\end{align}
which along with \eqref{w-Dt} and \eqref{Uj-w-esti}
yields \eqref{x-Uj-w},
thereby finishing the proof.
\hfill $\square$

Below we estimate the
growth in the remaining five unstable directions
associated to the null space of the operator $L$.

We define the renormalized variables $\wt e_j$ and $e_j$ by
\begin{align} \label{w-wtej-ej-def}
w(t,x)=\lbb_{j}(t)^{-\frac d2} \wt {e}_{j} (t,\frac{x-\a_{j}(t)}{\lbb_{j}(t)}) e^{i\theta_{j}(t)}, \ \
with\ \wt e_j(t,y)=e_j(t,y)e^{i(\b_j(t)\cdot y- \frac 14 \g_j(t) |y|^2)},
\end{align}
and let $e_{j,1} :={\rm Re} e_{j}$ and
$e_{j,2} := {\rm Im} e_{j}$,
where $1\leq j\leq K$.
Note that,
the renormalized variable $e_j$ is
different from the previous one $\epsilon_j$  in \eqref{wj-vej-def}.
The advantage to introduce $e_j$ and $\wt e_j$
can be seen in Proposition \ref{Prop-ej-Lkernel} below,
where the estimates of the unstable directions
can be diagonalized in some sense.

Using the Taylor expansions \eqref{g-gz-expan}, \eqref{g-gzz-expan}
we have
\begin{align}  \label{f-G1}
f(u+w)-f(u)=f^{\prime}(U) \cdot w+G_1,
\end{align}
where
\begin{align} \label{G1-equaw}
   G_1:=
   (\partial_z f)'(U, R)\cdot R w
    +  (\partial_{\ol z} f)'(U, R)\cdot R \ol w
    + f''(u, w)\cdot w^2.
\end{align}
We further split $f'(U; w)$ into three parts below
\begin{align} \label{f'w-G2G3}
   f'(U)\cdot w
   &= f'(U_j)\cdot w
      + \sum\limits_{l\not = j} f'(U_l)\cdot w
      + [f'(U) \cdot w - \sum\limits_{l=1}^K f'(U_l) \cdot w] \nonumber \\
   &=: f'(U_j)\cdot w + G_2 + G_3.
\end{align}
Note that,
$G_2$ contains the blow-up profiles different from $U_j$,
and $G_3$ contains the interactions between different blow-up profiles.

Moreover, let $G_4$ denote the lower order perturbations
\begin{align} \label{G4-bc}
   G_4 :=b\cdot \na w+cw ,
\end{align}
where $b,c$ are given by \eqref{b}, \eqref{c}, respectively.

Then, plugging \eqref{f-G1}, \eqref{f'w-G2G3} and \eqref{G4-bc} into
equation \eqref{equa-w}
we reformulate the equation of $w$ as follows
\begin{align} \label{equa-w-reform}
    i\partial_tw+\Delta w+f^{\prime}(U_j) \cdot w=-\sum_{l=1}^{4}G_l.
\end{align}

The equation of renomalized variable $e_j$
is contained in Lemma \ref{Lem-equa-ej} below,
which is, actually, a consequence of several algebraic cancellations.

\begin{lemma} \label{Lem-equa-ej}
For every $1\leq j\leq K$,
$e_j$ satisfies the equation
\begin{align} \label{equa-ej}
      i\lbb_j^2\partial_t e_j+\Delta e_j-e_j+(1+\frac{2}{d})Q^{\frac{4}{d}} e_j
        +\frac{2}{d}Q^{\frac{4}{d}} \ol {e_j}
    = -\sum_{l=1}^{4}H_l
      + \calo((\<y\>^2 |\wt e_j|+  \<y\> |\na \wt e_j|) Mod_j),
\end{align}
where
\begin{align} \label{Hl-Gl}
   H_l(t,y)
   = \lbb_j^{2+\frac d2} e^{-i\theta_j} e^{-i(\beta_j\cdot y - \frac{1}{4}\g_j|y|^2)} G_l(t,\lbb_j y + \a_j),
   \ \ 1\leq l\leq 4.
\end{align}
\end{lemma}

{\bf Proof.}
Using the identity
\begin{align*}
   \partial_t \wt e_j
   = \partial_t e_j e^{i(\beta_j\cdot y - \frac{1}{4}\g_j|y|^2)}
      + i (\dot \beta_j \cdot y - \frac 14 \dot \g_j |y|^2) \wt e_j,
\end{align*}
we infer from \eqref{w-wtej-ej-def} that,
if $y := \frac{x-\a_j}{\lbb_j}$,
\begin{align*}
   \partial_t w
   =& \lbb_j^{-2-\frac d2} e^{i\theta_j}
     \big(-\frac d2 \lbb_j \dot \lbb_j \wt e_j + \lbb_j^2 \partial_t e_j e^{i(\beta_j\cdot y - \frac{1}{4}\g_j|y|^2)}
       + i\lbb_j^2 \dot \beta_j\cdot y \wt e_j - \frac 14 i \lbb_j^2 \dot \g_j |y|^2 \wt e_j \\
    &\qquad \qquad \ \ -\lbb_j \dot \a_j \cdot \na\wt e_j
            -\lbb_j \dot \lbb_j y\cdot \na \wt e_j
            +i \lbb_j^2 \dot \theta_j \wt e_j  \big),
\end{align*}
which along with \eqref{Mod-def} yields that
\begin{align*}
   i \pa_t w
   =& \lbb_j^{-2-\frac d2} e^{i\theta_j}
     \big\{i\g_j \Lambda \wt e_j
       + i \lbb_j^2 \pa_t e_j e^{i(\beta_j\cdot y - \frac{1}{4}\g_j|y|^2)}
       + \g_j \beta_j \cdot y \wt e_j
       - \frac 14 \g_j^2 |y|^2 \wt e_j  \nonumber \\
    &\qquad \qquad \ \  - 2i \beta_j \cdot \na \wt e_j
       - \wt e_j
       - |\beta_j|^2 \wt e_j
     + \calo((\<y\>^2|\wt e_j| + \<y\>|\na \wt e_j|) Mod_j)\big\} .
\end{align*}
Then,
taking into account the identities,
similarly to \eqref{LaQj-LaQ} and \eqref{naQj-naQ},
\begin{align}
   & \Lambda \wt e_j = (\Lambda e_j + i(\beta_j \cdot y - \frac 12 \g_j |y|^2) e_j) e^{i(\beta_j\cdot y - \frac 14 \g_j |y|^2)}, \label{Lawtej-Laej} \\
   & \na \wt e_j = (\na e_j + i (\beta_j - \frac 12 \g_j y) e_j) e^{i(\beta_j\cdot y - \frac 14 \g_j |y|^2)}, \label{nawtej-naej}
\end{align}
we come to
\begin{align} \label{ptw}
  i\pa_t w
  =& \lbb_j^{-2-\frac d2} e^{i\theta_j} e^{i(\beta_j\cdot y - \frac{1}{4}\g_j|y|^2)}
    \big\{i \lbb_j^2 \pa_t e_j
          + i \g_j \Lambda e_j
          + |\beta_j - \frac 12 \g_j y|^2 e_j
          - 2 i \beta_j \cdot \na e_j  - e_j \nonumber \\
   & \qquad \qquad \qquad \qquad \qquad + \calo(\<y\>^2 |\wt e_j| + \<y\> |\na \wt e_j|) Mod_j \}.
\end{align}

Moreover,
by \eqref{w-wtej-ej-def},
direct computations show that
\begin{align} \label{Deltaw}
   \Delta \wt e_j
   = (\Delta e_j  - i \g_j \Lambda e_j
      - |\beta_j - \frac 12 \g_j y|^2 e_j
      + 2 i \beta_j \cdot \na e_j)  e^{i(\beta_j\cdot y - \frac{1}{4}\g_j|y|^2)} .
\end{align}

Thus, combing \eqref{ptw} and \eqref{Deltaw} altogether
we obtain that,
after algebraic cancellations,
\begin{align*}
   i \pa_t w + \Delta w
   = \lbb_j^{-2-\frac d2}
     e^{i\theta_j} e^{i(\beta_j\cdot y - \frac{1}{4}\g_j|y|^2)}
     \big\{ i\lbb_j^2 \pa_t e_j - e_j + \Delta e_j
           + \calo(\<y\>^2|\wt e_j| + \<y\>|\na \wt e_j|)Mod_j \big\},
\end{align*}
which along with \eqref{equa-w-reform} yields \eqref{equa-ej}
and thus  finishes the proof.
\hfill $\square$

The contributions of the error terms
$H_l$, $1\leq l\leq 4$,
are actually negligible,
which is the content of Lemma \ref{Lem-Hl} below.

\begin{lemma} \label{Lem-Hl}
Let $E$ belong to the generalized kernels of the linearized operator $L$,
i.e., $E\in\{Q, yQ, |y|^2 Q, \na Q, \Lambda Q, \rho\}$.
Then, there exists $C>0, \zeta,\delta\in (0,1)$ such that
\begin{align}
   & \int | H_1(t,y)||E(y)| dy\leq C(T-t)^{3+\zeta} \|w\|_{L^2},  \label{H1-E}   \\
   & \int ( |H_2(t,y)|+|H_3(t,y)|) |E(y)|dy \leq Ce^{-\frac{\delta}{T-t}}\|w\|_{L^2},  \label{H2H3-E}  \\
   & |\int H_4(t,y)E(y)dy| \leq C(T-t)^{\nu_*+1}\|w\|_{L^2}, \label{H4-E}
\end{align}
where $\nu_*$ is the flatness index of the spatial functions of noise in Assumption $(A1)$.
\end{lemma}

{\bf Proof.}
We first see that, by \eqref{Hl-Gl},
\begin{align}
   \int |H_1(t,y)| |E(y) | dy
   \leq C (T-t)^{2-\frac d2} \int |G_1(t,x) E(\frac{x-\a_j}{\lbb_j})| dx.
\end{align}
Since by \eqref{f'R}, \eqref{f''R2} and \eqref{G1-equaw},
\begin{align*}
   |G_1 |
   \leq C (|U|^{\frac 4d -1} + |R|^{\frac 4d -1}) |R| |w|
          + C (|U|^{\frac 4d -1} + |R|^{\frac 4d -1} + |w|^{\frac 4d -1}) |w|^2,
\end{align*}
taking into account $\|E\|_{L^\9}<\9$ and
using H\"older's inequality and \eqref{fLp-fH1}
we obtain
\begin{align*}
   \int|H_1(t,y)| |E(y)| dy
   \leq& C(T-t)^{2-\frac d2} \|w\|_{L^2}
         \big((T-t)^{-2+\frac d2}\|R\|_{L^2} + \|R\|_{H^1}^{\frac 4d} \\
      &\qquad \qquad  + (T-t)^{-2+\frac d2}\|w\|_{L^2}
          + \|R\|_{H^1}^{\frac 4d-1}\|w\|_{H^1} + \|w\|_{H^1}^\frac 4d \big),
\end{align*}
which along with \eqref{R-Tt-Uniq} and \eqref{w-Tt-Uniq} yields \eqref{H1-E}.

Moreover,
since by \eqref{f-linear}, \eqref{w-wtej-ej-def} and \eqref{f'w-G2G3},
\begin{align*}
   |G_2(t,x)|
   \leq \sum\limits_{l\not=j} |f'(U_l) \cdot w|
   \leq (T-t)^{-2-\frac d2} Q^{\frac 4d}(\frac{x-\a_l}{\lbb_l}) |e_j(t,\frac{x-\a_j}{\lbb_j})|,
\end{align*}
we see from \eqref{Hl-Gl} that
\begin{align}
   |H_2(t,y)| \leq C \sum\limits_{l\not =j} Q^\frac 4d( \frac{\lbb_j}{\lbb_l}y + \frac{\a_j-\a_l}{\lbb_l}) |e_j(t,y)| .
\end{align}
This yields that
\begin{align} \label{H2-E-esti}
    \int |H_2(t,y)| |E(y)| dy
   \leq& C \sum\limits_{l\not=j} \int Q^\frac 4d(\frac{\lbb_j}{\lbb_l}y + \frac{\a_j-\a_l}{\lbb_l}) |e_j (y)| |E(y)| dy \nonumber \\
   \leq& C \sum\limits_{l\not=j} e^{-\frac{\delta}{T-t}} \|e_j\|_{L^2}
   \leq C   e^{-\frac{\delta}{T-t}}  \|w\|_{L^2},
\end{align}
where $\delta\in(0,1)$,
and the second inequality  is due to the exponential decay  of
the ground state $Q$ and $E$.

Similarly,
by the definition of $G_3$
and \eqref{U-Tt},
$1\leq l\leq K$,
\begin{align*}
   |G_3(t,x)|
   \leq C (T-t)^{-2-\frac d2}
        \sum\limits_{l\not=h} Q(\frac{x-\a_l}{\lbb_l}) Q(\frac{x-\a_h}{\lbb_h}) |e_j(t,\frac{x-\a_j}{\lbb_j})|,
\end{align*}
which along with \eqref{Hl-Gl} yields that
\begin{align*}
   |H_3(t,y)|
    \leq& C \sum\limits_{l\not= j} Q(\frac{\lbb_j}{\lbb_l}y + \frac{\a_j- \a_l}{\lbb_l}) |e_j(t,y)|.
\end{align*}
Hence, similarly to \eqref{H2-E-esti},
we get that for some $\delta \in (0,1)$,
\begin{align*}
  \int  |H_3(t,y)||E(y) | dy
   \leq C e^{- \frac{\delta}{T-t}} \|w\|_{L^2}
\end{align*}
and thus \eqref{H2H3-E} follows.

Regarding $H_4$,
in view of \eqref{G4-bc}, \eqref{Hl-Gl} and \eqref{nawtej-naej},
we obtain
\begin{align*}
   H_4(t,y)
   =& \lbb_j \wt b \cdot (\na e_j + i(\beta_j - \frac 12 \g_j y)e_j)
     +  \lbb_j^2 \wt c e_j,
\end{align*}
where $\wt b(t,y) = b(t,\lbb_j y + \a_j)$,
$\wt c(t,y) =  c(t,\lbb_j y + \a_j)$.
Using integration by parts formula we have
\begin{align*}
   \int \wt b\cdot \na e_j E dy
   = - \int div \wt b\  E  e_j dy
     -\int \wt b \cdot \na E e_j dy,
\end{align*}
then using H\"older's inequality,
\eqref{b}, \eqref{c} and \eqref{phik-Taylor}
we obtain \eqref{H4-E} and finish the proof.
\hfill $\square$

By virtue of
the identities in \eqref{Q-kernel}
and Lemmas \ref{Lem-equa-ej} and \ref{Lem-Hl}
we are now able to control
the growth of $e_j$ along the remaining five unstable directions
as stated in Proposition \ref{Prop-ej-Lkernel} below.

\begin{proposition} \label{Prop-ej-Lkernel}
Let $e_j$ be as in \eqref{w-wtej-ej-def}
and $e_{j,1} :={\rm Re} e_{j}$,
$e_{j,2} := {\rm Im} e_{j}$,
Then, for $1\leq j\leq K$ and for $T$ small enough
we have
\begin{align}
&\frac{d}{dt}  \<e_{j,2}, \Lambda Q\>
     =  2\lbb_j^{-2}  \<e_{j,1}, Q\>
       + \calo((T-t)^{2+\zeta}\sqrt{D(t)}),  \label{LaQ-ej} \\
&\frac{d}{dt} \<e_{j,1}, |y|^2 Q\>
    = -4\lbb_j^{-2}  \<e_{j,2}, \Lambda Q\>
        + \calo((T-t)^{2+\zeta}\sqrt{D(t)}), \label{y2Q-ej} \\
&\frac{d}{dt} \<e_{j,2}, \rho\>
   = \lbb_j^{-2} \<e_{j,1}, |y|^2 Q\>
      + \calo((T-t)^{2+\zeta}\sqrt{D(t)}), \label{rho-ej} \\
&\frac{d}{dt}  \<e_{j,2}, \na Q\>
     =    \calo((T-t)^{2+\zeta}\sqrt{D(t)}).  \label{naQ-ej}  \\
&\frac{d}{dt}  \<e_{j,1}, yQ\>
     =  -2\lbb_j^{-2} \<e_{j,2}, \na Q\>
         + \calo((T-t)^{2+\zeta}\sqrt{D(t)}), \label{yQ-ej}
\end{align}
\end{proposition}

{\bf Proof.}
Let us take
$\frac{d}{dt} \<e_{j,2}, \Lambda Q \>$
as an example to illustrate the arguments.
Using equation \eqref{equa-ej}
we have
\begin{align} \label{rho-ej-equa}
   \frac{d}{dt}\<e_{j,2},  \Lambda Q \>
   =& \lbb_j^{-2} {\rm Re}
      \int  \Lambda Q (\Delta  e_j -  e_j + (1+\frac 2d)Q^\frac 4d  e_j + \frac 2d Q^{\frac 4d} \ol {e_j}) dy   \nonumber \\
    & + \lbb_j^{-2}\sum\limits_{l=1}^4 \calo(\int  \Lambda Q H_l dy)
      + \lbb_j^{-2} Mod_j \calo(\int \Lambda Q  (\<y\>^2|\wt  e_j|+ \<y\> |\na \wt  e_j|) dy).
\end{align}

Note that, by the definition \eqref{L+-L-} of the operator $L_+$,
the integration by parts formula and the identity $L_+  \Lambda Q  =- 2 Q$ in \eqref{Q-kernel},
\begin{align} \label{rho-ej-equa.1}
    {\rm Re}
      \int  \Lambda Q  (\Delta e_j - e_j + (1+\frac 2d)Q^\frac 4d e_j + \frac 2d Q^{\frac 4d} \ol {e_j}) dy
   =  -  \int  \Lambda Q  L_+ e_{j,1} dy
     =  2 \int  Q e_{j,1} dy.
\end{align}

Moreover,
using Lemma \ref{Lem-Hl}
and \eqref{w-Dt}
we infer that for each $1\leq l\leq 4$,
\begin{align}  \label{rho-ej-equa.2}
   \lbb_j^{-2} |\int  \Lambda Q   H_l dy|
   \leq  C (T-t)^{1+\zeta} \|w\|_{L^2}
   \leq C (T-t)^{2+\zeta} \sqrt{D(t)}.
\end{align}

We also note from Lemma \ref{Lem-Mod-Uniq} that
$Mod(t)\leq C (T-t)^{\kappa+3} \leq C(T-t)^{3+\zeta}$.
Then, using H\"older's inequality,
the boundedness of $\|\<y\>^2 \Lambda Q \|_{L^2}$
and
\begin{align*}
   \|\wt e_j\|_{L^2} + \|\na \wt e_j\|_{L^2}
   \leq C (\|w\|_{L^2} + \lbb_j \|\na w\|_{L^2})
   \leq C(T-t) \sqrt{D(t)},
\end{align*}
we get
\begin{align}  \label{rho-ej-equa.3}
   \lbb_j^{-2} Mod_j  |\int  \Lambda Q  (\<y\>^2|\wt e_j|+ \<y\> |\na \wt e_j|) dy|
   \leq& C (T-t)^{1+\zeta} \|\<y\>^2 \Lambda Q\|_{L^2} ( \|\wt e_j\|_{L^2} + \|\na \wt e_j\|_{L^2} ) \nonumber \\
   \leq& C (T-t)^{2+\zeta} \sqrt{D(t)}.
\end{align}

Thus, plugging \eqref{rho-ej-equa.1}-\eqref{rho-ej-equa.3} into \eqref{rho-ej-equa}
yields immediately \eqref{rho-ej}.

Similar arguments apply also to
the remaining four estimates in Proposition \ref{Prop-ej-Lkernel}.
For simplicity,
the details are omitted here.
\hfill $\square$

We are now ready to prove the key estimate \eqref{Scalj-N} in Theorem \ref{Thm-Scal}.

{\bf Proof of Theorem \ref{Thm-Scal}.}
Similarly to \eqref{Scaj-def},
we define the growth quantity associated to $e_j$ by
\begin{align}  \label{wtScalj-def}
   \wt{Scal}_j(t):= \<e_{j,1},Q\>^2+\<e_{j,1},yQ\>^2+\<e_{j,1},|y|^2Q\>^2
           +\<e_{j,2},\nabla Q\>^2+\<e_{j,2},\Lambda Q\>^2+\<e_{j,2},\rho\>^2.
\end{align}
As we shall see below, the two renomalized variables $\epsilon_j$ and $e_j$
defined in \eqref{wj-vej-def} and \eqref{w-wtej-ej-def} respectively
are almost the same
up to the negligible error of order $\calo(P + e^{-\frac{\delta}{T-t}})$,
and thus the two quantities $Scal_j$ and $\wt{Scal_j}$
should be close to each other.

We first claim that for some $\zeta>0$,
\begin{align} \label{wtScalj-N}
 \wt{Scal}_j(t)\leq C(T-t)^{2+ \zeta}  \sup\limits_{t\leq s<T} D(s).
\end{align}

To this end, we use Proposition \ref{Prop-U-w}
and the change of variables to get
\begin{align} \label{Q-ej-Tt}
  | \<e_{j,1}, Q\>| \leq C(T-t)^{4+\zeta} \sup\limits_{t\leq s<T} \sqrt{D(s)},
\end{align}
which along with \eqref{LaQ-ej} yields that
\begin{align*}
   |\frac{d}{dt} \<e_{j,2}, \Lambda Q\>| \leq C(T-t)^{2+\zeta} \sup\limits_{t\leq s<T} \sqrt{D(s)}.
\end{align*}
Since  by \eqref{w-Tt-Uniq}, $\lim_{t\to T} \|w(t)\|_{H^1} =0$,
we infer that
\begin{align}
   \lim\limits_{t \to T} \<e_{j,1}(t), \Lambda  Q\> =0.
\end{align}
This  yields  that
\begin{align} \label{LaQ-ej-Tt}
   |\<e_{j,2}, \Lambda Q\>|
   \leq \int_t^T |\frac{d}{ds}\<e_{j,2}, \Lambda Q\>| ds
   \leq C(T-t)^{3+\zeta} \sup\limits_{t\leq s<T} \sqrt{D(s)}.
\end{align}
Then, plugging \eqref{LaQ-ej-Tt} into \eqref{y2Q-ej} yields
\begin{align} \label{y2Q-ej-Tt}
   |\<e_{j,1}, |y|^2 Q\>| \leq C(T-t)^{2+\zeta} \sup\limits_{t\leq s<T} \sqrt{D(s)},
\end{align}
which along with \eqref{rho-ej} yields
\begin{align} \label{rho-ej-Tt}
   |\<e_{j,2}, \rho\>| \leq C(T-t)^{1+\zeta} \sup\limits_{t\leq s<T} \sqrt{D(s)}.
\end{align}
We also see from \eqref{naQ-ej} that
\begin{align} \label{naQ-ej-Tt}
   |\<e_{j,2}, \na Q\>|
 \leq C(T-t)^{3+\zeta}  \sup\limits_{t\leq s<T}  \sqrt{D(s)},
\end{align}
which along with \eqref{yQ-ej} yields that
\begin{align} \label{yQ-ej-Tt}
   |\<e_{j,1}, yQ\>| \leq C(T-t)^{2+\zeta} \sup\limits_{t\leq s<T} \sqrt{D(s)}.
\end{align}
Thus,
combining estimates \eqref{Q-ej-Tt}-\eqref{yQ-ej-Tt} altogether
we obtain \eqref{wtScalj-N}, as claimed.

Next we show that there exist $C,\delta>0$ such that
\begin{align} \label{Scalj-wtScalj}
   Scal_j(t) = \wt{Scal}_j(t) + \calo( P + e^{-\frac{\delta}{T-t}}) \|w(t)\|^2_{L^2}.
\end{align}
Then, taking into account \eqref{lbbj-gj-P-Tt-Uniq},
\eqref{w-Dt}
and \eqref{wtScalj-N}
we obtain \eqref{Scalj-N}.

It remains to prove \eqref{Scalj-wtScalj}.
For this purpose,
by \eqref{U-Uj-Uniq}, \eqref{w-wj}, \eqref{w-wtej-ej-def}
and the change of variables,
\begin{align*}
  {\rm Re} \<e_j, Q\>
  = {\rm Re} \< w, U_j\>
  =  {\rm Re}  \< w_j, U_j\>
     + \sum\limits_{l\not =j} {\rm Re}  \< w_l, U_j\>.
\end{align*}
Using Lemma \ref{Lem-inter-est} to decouple
$w_l$ and $U_j$, $l\not =j$,
and then using \eqref{wj-vej-def} we get
\begin{align*}
   {\rm Re} \<e_j, Q\>
   =  {\rm Re}   \< w_j, U_j\>
      + \calo(e^{-\frac{\delta}{T-t}} \|w\|_{L^2})
   = {\rm Re}  \<\epsilon_j, Q_j\>
      + \calo(e^{-\frac{\delta}{T-t}} \|w\|_{L^2}),
\end{align*}
where $Q_j$ is given by \eqref{Qj-Q-Uniq}.
Then, using the fact
\begin{align} \label{Qj-Q-P}
   Q_j = Q+ \calo(P \<y\>^2 Q),
\end{align}
we arrive at
\begin{align} \label{epj-vej-Q-p}
   {\rm Re} \<e_j, Q\>
   = {\rm Re}  \<\epsilon_j, Q\>
      + \calo(P  +e^{-\frac{\delta}{T-t}}) \|w\|_{L^2}.
\end{align}

Similar arguments with slight modifications also yield that
\begin{align}
    {\rm Re} \<e_j, y Q\>
  =& {\rm Re} \<\epsilon_j, y Q\>
    + \calo(P  +e^{-\frac{\delta}{T-t}}) \|w\|_{L^2},  \label{ejyQ-epjyQ}  \\
  {\rm Re} \<e_j, |y|^2 Q\>
      =&  {\rm Re} \<\epsilon_j, |y|^2 Q\>
        + \calo(P  +e^{-\frac{\delta}{T-t}}) \|w\|_{L^2},  \label{ejy2Q-epjy2Q}
\end{align}
and
taking into account $\rho_j = \rho + \calo(P\<y\>^2 \rho)$
we also have
\begin{align} \label{epj-vej-rho-p}
   {\rm Re} \<e_j, \rho\>
   = {\rm Re}  \<\epsilon_j, \rho\>
      + \calo(P  +e^{-\frac{\delta}{T-t}}) \|w\|_{L^2}.
\end{align}

Regarding ${\rm Im} \<\Lambda Q, e_j\>$,
using the change of variables
and the identity  \eqref{LaQj-LaQ}
we have
\begin{align*}
   {\rm Im} \<e_j,\Lambda Q\>
   =& {\rm Im} \int \wt e_j \ol{\Lambda Q e^{i(\beta_j\cdot y - \frac{1}{4}\g_j|y|^2)}} dy \\
   =& {\rm Im} \int w (\Lambda \ol {U_j} + i(\beta_j\cdot (\frac{x-\a_j}{\lbb_j}) - \frac 12 \g_j |\frac{x-\a_j}{\lbb_j}|^2) \ol {U_j}) dx \\
   =& {\rm Im} \int w_j (\Lambda  \ol {U_j} + i(\beta_j\cdot (\frac{x-\a_j}{\lbb_j}) - \frac 12 \g_j |\frac{x-\a_j}{\lbb_j}|^2) \ol{U_j} ) dx
      + \calo(e^{-\frac{\delta}{T-t}} \|w\|_{L^2}),
\end{align*}
where in the last step we also used Lemma \ref{Lem-inter-est} to
decouple the different profiles $w_l$ and $U_j$, $l\not= j$,
which merely contribute the exponentially small error.
Then, using again  \eqref{LaQj-LaQ}, \eqref{wj-vej-def},
the change of variables
and \eqref{Qj-Q-P}
we obtain
\begin{align}  \label{epj-vej-LaQ-p}
   {\rm Im} \<e_j,\Lambda Q\>
   =& {\rm Im} \int  \epsilon_j \ol{\Lambda Q e^{i(\beta_j\cdot y - \frac{1}{4}\g_j|y|^2)}} dy  + \calo( e^{-\frac{\delta}{T-t}} \|w\|_{L^2})  \nonumber   \\
    =& {\rm Im}\<  \epsilon_j, \Lambda Q \>+ \calo(P + e^{-\frac{\delta}{T-t}}) \|w\|_{L^2}.
\end{align}

Using similar arguments,
but with the identity  \eqref{naQj-naQ} instead,
we also have
\begin{align}  \label{epj-vej-naQ-p}
   {\rm Im} \<e_j, \na Q\>
    =& {\rm Im} \<  \epsilon_j, \na Q \> + \calo(P + e^{-\frac{\delta}{T-t}}) \|w\|_{L^2}.
\end{align}

Therefore, combining estimates \eqref{epj-vej-Q-p},
\eqref{ejyQ-epjyQ}, \eqref{ejy2Q-epjy2Q}, \eqref{epj-vej-rho-p},
\eqref{epj-vej-LaQ-p}
and \eqref{epj-vej-naQ-p} altogether
we obtain \eqref{Scalj-wtScalj}
and thus finish the proof of Theorem \ref{Thm-Scal}.
\hfill $\square$

Now, we are ready to prove the main uniqueness result, i.e., Theorem \ref{Thm-Uniq-Blowup-RNLS}.

{\bf Proof of Theorem \ref{Thm-Uniq-Blowup}.}
First take any $\ve\in (0,\ve^*]$,
where $\ve^*$ is sufficiently small and is to be specified later.
Then,
let $T$ be sufficiently small and satisfy \eqref{T-M}.
We shall use iteration arguments to show that $D\equiv 0$.

By Theorem \ref{Thm-Nt-Scal}, we have for any $t\in [0,T)$,
\begin{align} \label{D-iter}
   \sup\limits_{t\leq s<T} D(s)
\leq C_1 \sup_{t\leq s<T}\sum_{j=1}^{K}\frac{Scal_j(s)}{\lbb_j^2(s)}
+C_1 \int_{t}^{T}\sum_{j=1}^{K}\frac{Scal_j(s)}{\lbb_j^3(s)}
   + \ve^*  \frac{D(s)}{T-s} ds.
\end{align}
Then, in view of Theorem \ref{Thm-Scal} and \eqref{w-Tt-Uniq}
we obtain that for some $\zeta>0$,
\begin{align} \label{D-iter*}
   \sup\limits_{t\leq s<T} D(s)
\leq C_2 (T-t)^{\zeta} \sup\limits_{t\leq s<T} D(s)
      +  C_2 \ve^* \int_t^T \frac{D(s)}{T-s} ds,
\end{align}
where $C_2$ is independent of $\ve^*$.
This yields that for $T$ even smaller
such that $C_2 T^{\zeta} \leq \frac 12$,
\begin{align} \label{D-iter.0}
   \sup\limits_{t\leq s<T} D(s)
\leq 2 C_2 \ve_* \int_t^T \frac{D(s)}{T-s} ds.
\end{align}

We also infer from \eqref{w-Tt-Uniq} and \eqref{D-def} that
\begin{align*}
   D(s) \leq C_3 (T-s)^{4+\zeta},
\end{align*}
where we may take $C_3\geq 1$, independent of $\ve^*$.

Then, plugging this into \eqref{D-iter.0} we get
\begin{align}  \label{D-iter.1}
   \sup\limits_{t\leq s<T} D(s)
   \leq 2 C_2 \ve^* \int_t^T C_3 (T-s)^{3+\zeta} ds
   \leq (\frac{2C_2C_3\ve^*}{4+\zeta}) (T-t)^{4+\zeta}.
\end{align}
But, plugging \eqref{D-iter.1} into \eqref{D-iter.0} again
we can obtain the refined estimate
\begin{align} \label{D-iter.2}
   \sup\limits_{t\leq s<T} D(s)
\leq  (\frac{2 C_2C_3 \ve^*}{4+\zeta})^2 (T-t)^{4+\zeta}.
\end{align}
Thus, iterating similar arguments
we infer that for any $k \geq 1$,
\begin{align} \label{D-iter.k}
   \sup\limits_{t\leq s<T} D(s)
\leq  (\frac{2 C_2C_3 \ve^*}{4+\zeta})^k (T-t)^{4+\zeta}
\leq  (\frac{2 C_2C_3 \ve^*}{4+\zeta})^k,
\end{align}
where $C_2, C_3>0$ are independent of $\ve^*$ and $k$.

Therefore,
taking $\ve^*$ small enough
such that $\frac{2 C_2C_3 \ve^*}{4+\zeta} <1$
and using the arbitrariness of $k$
we infer that
$D(t) =0$ for any $t\in [0,T)$,
which yields immediately that $w\equiv 0$,
and thus $v\equiv u$.
The proof of Theorem \ref{Thm-Uniq-Blowup} is complete.
\hfill $\square$

\section{Appendix} \label{Sec-App}

{\bf Proof of Lemma \ref{Lem-inter-est}.}
First,
using \eqref{aj-xj}
and the separation of $\{x_j\}_{j=1}^K$
we have that
for any $t^*\leq t<T_*$,
\begin{align}
   |(\alpha_j(t)-\alpha_l(t))\cdot {\bf v_1}|\geq 10\sigma, \ \ j\neq l, \label{aj-al-sep}
\end{align}
where $\sigma$ is given by \eqref{sep-xj-0}.
Note that
\begin{align}  \label{lbbj-lbbl-bdd}
\max_{1\leq j\neq l\leq K} \frac{\lambda_{j}(t)}{\lambda_{l}(t)}
\leq c_*:=4 \max_{1\leq j\neq l\leq K} \frac{\omega_{j}}{\omega_{l}},\ \
   \max_{1\leq j\neq l\leq K} |\a_j- \a_l| \leq 2 M := 2 (1+\max_{1\leq j\leq K} |x_j|).
\end{align}

In order to prove \eqref{Gj-Gl-decoup}, we note that
\begin{align*}
      & |\int\limits_{\bbr^d} |x-\a_l|^n |\partial^\nu G_l(t)| |x-\a_j|^m |G_j(t)|dx   \\
  \leq& C (T-t)^{-|\nu|+m+n-d}
         \int_{\R^d} |\frac{x-\a_l}{\lbb_l}|^n |\frac{x-\a_j}{\lbb_j}|^m
             |\pa^\nu g_l(t,\frac{x-\a_l}{\lbb_l})|
             |g_l(t,\frac{x-\a_j}{\lbb_j})| dx    \\
  \leq& C (T-t)^{-|\nu|+m+n}
          \int |\frac{\lbb_j y + \a_j - \a_l}{\lbb_l}|^n |y|^m
               | \partial^\nu g_l (t, \frac{\lbb_j}{\lbb_l} y + \frac{\a_j- \a_l}{\lbb_l})|
               | g(y)| dy.
\end{align*}
Using \eqref{T-M} we have
\begin{align*}
    |\frac{\lbb_j y + \a_j - \a_l}{\lbb_l}|^n
    \leq C(|y|^n + (T-t)^{-n} M^n)
    \leq C (T-t)^{-2n} \<y\>^n.
\end{align*}
It follows that
\begin{align*}
       &  |\int\limits_{\bbr^d} |x-\a_l|^n |\partial^\nu G_l(t)| |x-\a_j|^m |G_j(t)|dx  \\
   \leq& C (T-t)^{-|\nu|+m-n}
         ( \int_{|y\cdot {\bf v_1}|\leq\frac{5\sigma}{c_*\lambda_j}}
           + \int_{|y\cdot {\bf v_1}|\geq\frac{5\sigma}{c_*\lambda_j}} )
        \<y\>^{m+n}
        |\partial^\nu g_l(t, \frac{\lbb_j}{\lbb_l} y + \frac{\a_j-\a_l}{\lbb_l})|
        |g(y)| dy.
\end{align*}
On one hand,
in the region $\{y\in \bbr^d: |y\cdot {\bf v_1}|\leq \frac{5\sigma}{c_*\lbb_j}\}$,
by \eqref{aj-al-sep} and \eqref{lbbj-lbbl-bdd},
\begin{align*}
   |\frac{\lambda_{j}}{\lambda_{l}}y +\frac{\alpha_{j}-\alpha_l}{\lambda_{l}}|
   \geq& |\frac{(\alpha_{j}-\alpha_l)\cdot {\bf v_1}}{\lambda_{l}}| - |\frac{\lambda_{j}}{\lambda_{l}}(y\cdot {\bf v_1})|
    \geq  \frac{5\sigma}{\lbb_l} \to \9,\ \ as\ T\to 0.
\end{align*}
Then, in view of the exponential decay of $\partial_\nu g$,
we obtain that
\begin{align*}
   &(T-t)^{-|\nu|+m-n}
   \int_{|y\cdot {\bf v_1}|\leq\frac{5\sigma}{c_*\lambda_l}}
     \<y\>^{m+n} |\partial^\nu g_l (t,\frac{\lambda_{j}}{\lambda_{l}}y +\frac{\alpha_{j}-\alpha_l}{\lambda_{l}}) |g(y)| dy \\
   \leq& C (T-t)^{-|\nu|+m-n} e^{-\frac{\delta_1}{T-t}}
          \|\<y\>^{m+n} g\|_{L^1}
   \leq C e^{-\frac{\delta_1}{2(T-t)}},
\end{align*}
where $\delta_1>0$,
and we also used $\sup_{r>0} r^{-|\nu|+m-n} e^{-\frac{\delta_1}{2r}} <\9$ in the last step.
On the other hand,
in the region $\{y: |y\cdot {\bf v_1}|\geq \frac{5\sigma}{c_*\lbb_l}\}$,
$|y| \geq |y\cdot {\bf v_1}| \geq \frac{5\sigma}{c_*\lbb_l} \to \9$, as $t\to T$.
Using the change of variables
and the exponential decay of $g$
we get
\begin{align*}
   & (T-t)^{-|\nu|+m-n}
   \int_{|y\cdot {\bf v_1}|\geq\frac{5\sigma}{c_*\lambda_l}}
    \<y\>^{m+n} |\partial^\nu g_l(t,,\frac{\lambda_{j}}{\lambda_{l}}y +\frac{\alpha_{j}-\alpha_l}{\lambda_{l}})|
    |g(y)| dy  \\
   \leq& C(T-t)^{-|\nu|+m-n} e^{-\frac{\delta_2}{T-t}}
         \|\<y\>^{m+n} |g(y)|^\frac 12\|_{L^\9}
         \|\partial^\nu g_l\|_{L^1}
   \leq C e^{-\frac{\delta_2}{2(T-t)}},
\end{align*}
where $\delta_2>0$.
Thus, combing two estimates above we obtain \eqref{Gj-Gl-decoup}.

Regarding \eqref{Gj-hl-decoup},
since for $j\not = l$, $|(x-x_l)\cdot {\bf v_1}|\geq 4\sigma$
on the support of $\Phi_j$,
we infer from \eqref{aj-xj} that
$|(x-\a_l)\cdot {\bf v_1}| \geq 3\sigma$.
Then,
using the change of variables
and the inequality
\begin{align*}
    |\frac{\lbb_l y + \a_l - \a_j}{\lbb_j}|^m
    \leq C M (T-t)^{-m} \<y\>^m
    \leq C (T-t)^{-2m} \<y\>^m
\end{align*}
we get
\begin{align} \label{Gl-hPhij-esti}
    &  |\int_{\bbr^d} |x-\a_l|^n | \partial^\nu G_l(t)| |x-\a_j|^m |h| \Phi_j dx |  \nonumber \\
\leq&  C (T-t)^{-|\nu|+m+n-\frac d2}
       \int_{|(x-x_l)\cdot {\bf v_1}|\geq 4\sigma}
       |\frac{x-\a_l}{\lbb_l}|^n
       |\partial^\nu g_l(t, \frac{x-\a_l}{\lbb_l})| |\frac{x-\a_j}{\lbb_j}|^m
       |h(x)|  dx  \nonumber \\
\leq&C (T-t)^{-|\nu|-m+n+\frac d2}
         \int_{|y\cdot {\bf v_1}|\geq \frac{3\sigma}{\lambda_l}}
         \<y\>^{m+n} |\partial^\nu g_l(t,y) |
         |h(\lambda_{l}y+\alpha_l)|dy.
\end{align}
Then,
using Cauchy's inequality and the exponential decay of $\partial^\nu g$ again
we obtain that for $t$ close to $T$,
the right-hand-side of \eqref{Gl-hPhij-esti} is bounded by
\begin{align*}
   & C(T-t)^{-|\nu|-m+n+\frac d2}
         e^{-\frac{\delta_3}{T-t}}
         \|\<y\>^{m+n}| \partial^\nu g_l|^\frac 12 \|_{L^2}
         (\int|h(\lbb_l y + \a_l)|^2 dy)^\frac 12  \\
   \leq& C(T-t)^{-|\nu|-m+n}
         e^{-\frac{\delta_3}{T-t}}
         \|\<y\>^{m+n}| \partial^\nu g_l|^\frac 12 \|_{L^2}
         \|h\|_{L^2}
   \leq C e^{-\frac{\delta_3}{2(T-t)}} \|h\|_{L^2},
\end{align*}
where $\delta_3>0$
and in the last step we also used $\sup_{r>0}r^{-|\nu|-m+n}
         e^{-\frac{\delta_3}{2r}}<\9$.
One can also bound the right-hand-side of \eqref{Gl-hPhij-esti}
by
\begin{align*}
   & C(T-t)^{-|\nu|+n+\frac d2}
     e^{-\frac{\delta_3}{T-t}}
    \|\<y\>^{m+n}| \partial^\nu g_l|^\frac 12 \|_{L^\9}
      \int|h(\lbb_l y + \a_l)|  dy \\
   \leq&  C(T-t)^{-|\nu|-m+n-\frac d2}
     e^{-\frac{\delta_3}{T-t}}
    \|\<y\>^{m+n}| \partial^\nu g_l|^\frac 12 \|_{L^\9} \|h\|_{L^1}
   \leq C e^{-\frac{\delta_3}{2(T-t)}} \|h\|_{L^1}.
\end{align*}
Hence, combining the two estimate together
we obtain \eqref{Gj-hl-decoup},
thereby finishing the proof of Lemma \ref{Lem-inter-est}.
\hfill $\square$

{\bf Proof of Corollary \ref{Cor-coer-f-local}.}
Let $\wt {f} :=f\phi_A^{\frac{1}{2}}$.
Since
$\nabla f\phi_A^{\frac{1}{2}}=\nabla \wt {f}-\frac{\nabla\phi_A}{2\phi_A} \wt {f}$,
we have
\begin{align}\label{a4}
&\int (|\nabla f|^2 +|f|^2)\phi_A -(1+\frac{4}{d})Q^{\frac4d}f_1^2-Q^{\frac4d}f_2^2dx \nonumber \\
=&\int|\nabla \wt {f}|^2+|\wt {f}|^2-(1+\frac{4}{d})Q^{\frac4d} \wt {f}_1^2-Q^{\frac4d} \wt {f}_2^2dx \nonumber \\
&- \int(1-\phi_A)((1+\frac{4}{d})Q^{\frac4d}f_1^2+Q^{\frac4d}f_2^2)dx \nonumber \\
&+\frac{1}{4}\int |\frac{\nabla\phi_A}{\phi_A} |^2|\tilde{f}|^2dx
-{\rm Re} \int\frac{\nabla\phi_A}{\phi_A}\cdot\nabla \wt {f} \ol{\wt {f}}dx
=: \sum\limits_{i=1}^4 K_i.
\end{align}

Using Corollary \ref{Cor-coer-f-H1} we have
\begin{align} \label{K1-fH1-Scal.0}
 K_1 \geq C'_1\|\wt {f}\|_{H^1}^2- C'_2Scal(\tilde{f}),
\end{align}
where $C'_1, C'_2>0$.
Moreover, we claim that there exist $C, \delta>0$ such that
\begin{align}   \label{Scal-f-wtf}
Scal(\wt {f})\leq Scal(f) + C e^{-\delta A}\|\wt f\|_{L^2}^2.
\end{align}
To this end,
we see that
\begin{align*}
   \<\wt f_1, Q\> = \< f_1, Q\> + \<\wt f_1\phi_A^{-\frac 12}(\phi_A^\frac 12 -1), Q\>.
\end{align*}
Since on the support of $\phi_A^\frac 12 -1$,
$|x|\geq A$,
by the exponential decay of $Q$ we have that for some $\delta'>0$
\begin{align*}
   |\phi_A^{-\frac 12} (x) Q^\frac 12(x)|
   \leq C e^{-\frac 12 |x|(\delta' -A^{-1})}
   \leq C e^{-\frac 14 \delta' A}, \ \ |x|\geq A,
\end{align*}
where we also took $A$ large enough such that $\delta' -A^{-1} \geq \frac {1}{2}\delta'$.
It follows that
\begin{align} \label{phiA-Q-exp}
    | \<\wt f_1\phi_A^{-\frac 12}(\phi_A^\frac 12 -1), Q\>|
    \leq C \|\wt f\|_{L^2} \|\phi_A^{-\frac 12} Q^\frac 12\|_{L^\9(|x|\geq A)} \|Q^\frac 12\|_{L^2}
    \leq C e^{-\frac 14 \delta' A} \|\wt f\|_{L^2},
\end{align}
which yields that
\begin{align*}
   | \<\wt f_1, Q\> - \< f_1, Q\>|
   \leq  C  e^{-\frac 14 \delta' A} \|\wt f\|_{L^2}.
\end{align*}
Similar arguments apply also to the remaining five scalar products in $Scal(f)$,
and thus we obtain \eqref{Scal-f-wtf}, as claimed.

Hence, we infer from \eqref{K1-fH1-Scal.0} and \eqref{Scal-f-wtf} that
for $C_1, C_2>0$,
\begin{align} \label{K1-fH1-Scal}
 K_1 \geq C_1\|\tilde{f}\|_{H^1}^2- C_2Scal(f).
\end{align}

Regarding the second term $K_2$, we see that
\begin{align}
 K_2= \int (1-\phi_A)\phi_A^{-1}((1+\frac 4d)Q^2 \wt{f}_1^2+Q^2 \wt{f}^2_2)dx
    \leq  C \|\phi_A^{-1} Q^2 \|_{L^\9(|x|\geq A)} \|\wt f\|_{L^2}.
\end{align}
Similar arguments as in the proof of \eqref{phiA-Q-exp}
yield that for $A$ large enough
$ \|\phi_A^{-1} Q^2 \|_{L^\9(|x|\geq A)}  \leq C  e^{- \delta A}$,
where $C,\delta>0$.
This implies that for $A$ large enough
\begin{align} \label{K2-fH1-Scal}
    K_2 \leq  C  e^{- \delta A}\|\wt f\|_{L^2}.
\end{align}

We also note that,
since $|\frac{\nabla\phi_A}{\phi_A}|\leq C A^{-1}$,
by  H\"{o}lder's inequality,
\begin{align} \label{K3K4-fH1-Scal}
 K_3 + K_4 \leq \frac C A \|\tilde{f}\|_{H^1}^2.
\end{align}

Thus, plugging (\ref{K1-fH1-Scal}), (\ref{K2-fH1-Scal}) and (\ref{K3K4-fH1-Scal}) into (\ref{a4}),
we obtain  that for $A$ possibly larger
\begin{align} \label{esti-wtf}
  \int (|\nabla f|^2 +|f|^2)\phi_A -(1+\frac{4}{d})Q^{\frac4d}f_1^2-Q^{\frac4d}f_2^2dx
  \geq \frac{C_1}{2} \|\wt f\|_{H^1}^2-C_2Scal(f).
\end{align}

It remains to show that for $A$ large enough,
\begin{align} \label{wtf-fphiA}
     \frac{C_1}{2} \|\wt f\|_{H^1}^2
     \geq \frac{C_1}{8} \int (|\nabla f|^2 +|f|^2)\phi_A dx.
\end{align}

To this end, since
$\nabla\tilde{f} = \nabla f\phi_A^{\frac{1}{2}}+ \frac{\nabla\phi_A}{2\phi_A}\tilde{f}$
we have
\begin{align} \label{esti-wtf.1}
   \frac{C_1}{2}\|\wt f\|_{H^1}^2
   =&  \frac{C_1}{2} \int |f|^2 \phi_A dx
     + \frac{C_1}{2} \int |\nabla f\phi_A^{\frac{1}{2}}+ \frac{\nabla\phi_A}{2\phi_A}\tilde{f}|^2 dx \nonumber \\
   =&   \frac{C_1}{2} \int (|f|^2 + |\na f|^2) \phi_A dx
      +  \frac{C_1}{2} \int {\rm Re} (\na f \phi^\frac 12_A \frac{\na \phi_A}{\phi_A} \ol{\wt f}) dx
      +  \frac{C_1}{8} \int |\frac{\na \phi_A}{\phi_A} \wt f|^2 dx.
\end{align}
Note that,
by H\"older's inequality and $|\frac{\nabla\phi_A}{\phi_A}|\leq C A^{-1}$,
\begin{align} \label{esti-wtf.1}
   |\frac{C_1}{2} \int {\rm Re} (\na f \phi^\frac 12_A \frac{\na \phi_A}{\phi_A} \ol{\wt f}) dx |
   \leq& \frac{C_1C}{2A} (\int |\na f|^2 \phi_A dx)^{\frac 12}
        (\int | f|^2 \phi_Adx)^\frac 12 \nonumber  \\
   \leq& \frac{C_1}{4} \int |\na f|^2 \phi_A dx
         + \frac{4C_1 C^2}{A^2} \int |f|^2 \phi_A dx.
\end{align}
We also have
\begin{align} \label{esti-wtf.2}
   \frac{C_1}{8} \int |\frac{\na \phi_A}{\phi_A} \wt f|^2 dx
  \leq \frac{C_1 C^2}{2A^2} \int |f|^2 \phi_A dx.
\end{align}
Therefore,
combining \eqref{esti-wtf}-\eqref{esti-wtf.2}
and taking $A$ large enough
such that
$\frac{4C^2}{A^2} + \frac{C^2}{2A^2} \leq \frac 18$
we obtain \eqref{wtf-fphiA}, as claimed.
This along with \eqref{esti-wtf} yields \eqref{coer-f-local}.
\hfill $\square$

{\bf Proof of Lemma \ref{Lem-eeta}.}
Using  \eqref{equa-Ut} and \eqref{etan-Rn} we have
\begin{align} \label{eta-esti.1}
   \|\partial^\nu \eta\|_{L^2}
   \leq C(T-t)^{-2-|\nu|} Mod
        + \sum\limits_{j=1}^K  \|\partial^\nu(b\cdot \na U_j + c U_j)\|_{L^2}
        + \|\partial^\nu(f(U)- \sum\limits_{j=1}^K f(U_j))\|_{L^2}.
\end{align}
Note that,
\begin{align*}
   \|\partial^\nu (b\cdot \na U_j) \|_{L^2}
   =& \sum\limits_{\nu_1+\nu_2} \|\partial^{\nu_1} b \cdot \partial^{\nu_2} \na U_j \|_{L^2}   \\
   \leq& C (T-t)^{-|\nu_2|-1}
         \sum\limits_{j=1}^K ( \int |(\partial^{\nu_1} b) (\lbb_j y +\a_j) \partial^{\nu_2} \na Q_j(y)|^2 dy )^\frac 12.
\end{align*}
Since by \eqref{b} and \eqref{phik-Taylor},
\begin{align*}
   |(\partial^{\nu_1} b)(\lbb_jy +\a_j)|
   \leq C (T-t)^{\nu_*-|\nu_1|} \<y\>^{\nu_*+1}.
\end{align*}
This yields that
\begin{align} \label{eta-b-esti}
   \|\partial^\nu (b\cdot \na U_j) \|_{L^2}
   \leq C (T-t)^{\nu_* - |\nu| -1}.
\end{align}
Similarly, using \eqref{c} and \eqref{phik-Taylor} again
we also have
\begin{align} \label{eta-c-esti}
   \|\partial^\nu  (cU_j) \|_{L^2}
   \leq C (T-t)^{\nu_* - |\nu| -1}.
\end{align}

Moreover, we have that for $T$ small enough,
\begin{align}  \label{pxfU-fUj}
  \|\partial^\nu(f(U)- \sum\limits_{j=1}^K f(U_j))\|_{L^2}
  \leq& C(T-t)^{-\frac 8d (2+\frac d2)} (\sum\limits_{k\not =l} \int |U_k||U_l| dx)^{\frac 12}  \nonumber \\
  \leq& C(T-t)^{-\frac 8d (2+\frac d2)} e^{-\frac{\delta}{T-t}}
  \leq C(T-t)^{\nu_* - |\nu| -1}.
\end{align}

Therefore,
plugging \eqref{eta-b-esti}, \eqref{eta-c-esti} and \eqref{pxfU-fUj}  into \eqref{eta-esti.1}
we arrive at
\begin{align*}
   \| \partial^\nu \eta\|_{L^2}
   \leq C (T-t)^{-2-|\nu|} Mod + C (T-t)^{\nu_* - |\nu| -1},
\end{align*}
which, via \eqref{Mod-w-lbb} and $\nu_*=\kappa+3$,
yields \eqref{eta-L2} immediately.
The proof is complete.
\hfill $\square$

\section*{Acknowledgements}
The authors would like to thank Prof. Daomin Cao and Prof. Michael R\"ockner for valuable comments to improve the manuscript.
Y. Su is supported by NSFC (No. 11601482) and D. Zhang  is supported by NSFC (No. 11871337).
Financial support by the Key Laboratory of Scientific and Engineering Computing
(Ministry of Education) is also gratefully acknowledged.

\end{document}